\let\chooseClass2

\ifx\chooseClass1
\IfFileExists{extarticle.cls}{
\documentclass[14pt]{extarticle}
\emergencystretch 7 pt
\setlength{\oddsidemargin}{-17 mm} %-17     %-19    %-19     %-19
\setlength{\textwidth}{192 mm}     %192     %196    %196     %196
\setlength{\textheight}{246 mm}    %246     %246    %270     %266
\setlength{\topmargin}{-31 mm}     %-31     %-31    %-31     %-31
}{
\documentclass[12pt]{article}
\textheight = 8.5in
\textwidth 6.3in
\setlength{\oddsidemargin}{0mm}
\setlength{\topmargin}{-20 mm}
}{}
\fi%\chooseClass1

\ifx\chooseClass2
\documentclass[12pt]{article}
\fi

\usepackage{    amsmath,
                amsfonts,
                amssymb,
                amsthm
}

\IfFileExists{url.sty}{\usepackage{url}}{}
\providecommand{\url}[1]{{\tt #1}}

\usepackage{ifpdf}
\ifpdf
\usepackage[pdftex]{graphics}
\IfFileExists{hyperref.sty}{\usepackage[pdftex]{hyperref}}{}
\else
\usepackage[dvips]{graphics}
\IfFileExists{hyperref.sty}{\usepackage[hypertex]{hyperref}}{}
\fi

\IfFileExists{euscript.sty}{\usepackage[mathcal]{euscript}}{}
\IfFileExists{mathrsfs.sty}{\usepackage{mathrsfs}}{\let\mathscr\mathfrak}

\InputIfFileExists{diagrams.tex}{}{%
\IfFileExists{diagrams.sty}{\usepackage{diagrams}}
{\message{Without diagrams you can not process this file!}}}
\diagramstyle[height=2em,balance,righteqno,PostScript=dvips,nohug]

\IfFileExists{dsfont.sty}%
{\usepackage[sans]{dsfont}\newcommand\1{{\mathds 1}}}%
{\newcommand\1{{1\mkern-5mu {\mathrm I}}}}

\makeatletter
\def\@seccntformat#1{\csname the#1\endcsname.\quad}
\renewcommand\section{\@startsection {section}{1}{\z@}%
                                   {-3.5ex \@plus -1ex \@minus -.2ex}%
                                   {2.3ex \@plus.2ex}%
                                   {\normalfont\large\bfseries}}
\renewcommand\subsection{\@startsection{subsection}{2}{\z@}%
                        {3.25ex plus 1ex minus .2ex}{-.5em}%
                        {\normalfont\normalsize\bfseries}}
\@addtoreset{equation}{section}
\numberwithin{equation}{subsection}
\makeatother

\newtheoremstyle{boldhead}%     name
{\topsep}%                      abovespace
{\topsep}%                      belowspace
{\slshape}%                     bodyfont
{}%                             indentation=noindent
{\bfseries}%                    headfont
{.}%                            headpunctuation
{ }%                            headspace=interword space
{\thmname{#1}\thmnumber{ #2}\thmnote{ (#3)}}%   custom head specification

\newtheoremstyle{boldremark}%   name
{\topsep}%                      abovespace
{\topsep}%                      belowspace
{\upshape}%                     bodyfont
{}%                             indentation=noindent
{\bfseries}%                    headfont
{.}%                            headpunctuation
{ }%                            headspace=interword space
{\thmname{#1}\thmnumber{ #2}\thmnote{ (#3)}}%   custom head specification

\swapnumbers

\theoremstyle{boldhead}

\newtheorem{corollary}[subsection]{Corollary}
\newtheorem{lemma}[subsection]{Lemma}
\newtheorem{proposition}[subsection]{Proposition}

\theoremstyle{boldremark}

\newtheorem{definition}[subsection]{Definition}
\newtheorem{example}[subsection]{Example}

\newtheorem{remark}[subsection]{Remark}

\newtheorem*{acknowledgement}{Acknowledgements}

\def\rhaha{\raise.24ex\hbox{$\rightharpoonup$}\kern-1em\lower.24ex\hbox{$\rightharpoondown$}}%
\def\lhaha{\raise.24ex\hbox{$\leftharpoonup$}\kern-1em\lower.24ex\hbox{$\leftharpoondown$}}%
\def\dhaha{\downharpoonleft\kern-.22em\downharpoonright\kern.02em}%
\def\uhaha{\upharpoonleft\kern-.22em\upharpoonright\kern.02em}% -.21em
\newarrowhead{twoharpoons}\rhaha\lhaha\dhaha\uhaha

\newarrow{MapsTo}{mapsto}---{->}
\newarrow{Mono}{boldhook}---{->}
\newarrow{TTo}----{->}
\newarrow{Twoar}===={=>}

\newcommand\ZZ{{\mathbb Z}}

\newcommand{\ca}{{\mathcal A}}
\newcommand{\cb}{{\mathcal B}}
\newcommand{\cc}{{\mathcal C}}
\newcommand{\cd}{{\mathcal D}}
\newcommand{\ce}{{\mathcal E}}
\newcommand{\cf}{{\mathcal F}}

\newcommand{\ci}{{\mathcal I}}
\newcommand{\cj}{{\mathcal J}}
\newcommand{\ck}{{\mathcal K}}
\newcommand{\co}{{\mathcal O}}
\newcommand{\cx}{{\mathcal X}}

\newcommand{\fu}{{\mathscr U}}

\newcommand{\bull}{{\scriptscriptstyle\bullet}}

\newcommand{\tensor}[1]{\overset{#1}{\otimes}}
\newcommand{\uni}{{\mathbf i}}
\newcommand{\iu}{\bar{\imath}}
\newcommand\ju{{\overline{\jmath}}}
\newcommand\bm{{\boldsymbol\mu}}
\newcommand{\ovmi}{\overline{\;\text-\;}}
\newcommand{\unmi}{\underline{\;\text-\;}}

\newcommand{\sS}[2]{\vphantom{#2}#1 #2}
\newcommand{\tcolimit}{{\displaystyle\lim_{\longrightarrow}}\raisebox{-1.8mm}{\vphantom{l}}}

\newcommand{\n}[1]{\nobreakdash-\hspace{0pt}}
\newcommand{\ainf}[1]{$A_\infty$\nobreakdash-\hspace{0pt}}

\newcommand{\plquot}{/\mkern-4mu_p\mkern2mu}
\newcommand{\Cat}{{\mathcal C}at}
\newcommand{\Acyc}{\mathsf A}
\newcommand{\AcIn}{\mathsf{AI}}
\newcommand{\Com}{\mathsf C}
\newcommand{\Dc}{\mathsf D}
\newcommand{\Dr}{\mathsf D}
\newcommand{\Inj}{\mathsf I}
\newcommand{\Kht}{\mathsf K}
\newcommand{\quo}{\mathsf q}
\newcommand\cQuiver{{\mathscr Q}}

\let\kk\Bbbk

\let\eps\varepsilon
\let\epsilon\varepsilon
\let\ge\geqslant
\let\le\leqslant

\let\tens\otimes

\DeclareMathOperator\ad{ad}

\DeclareMathOperator\Cone{Cone}

\DeclareMathOperator\Hom{Hom}

\DeclareMathOperator\id{id}

\DeclareMathOperator\inj{in}
\DeclareMathOperator\im{Im}
\DeclareMathOperator\Ker{Ker}
\DeclareMathOperator\modul{-mod}

\DeclareMathOperator\Ob{Ob}
\DeclareMathOperator\pr{pr}

\newcommand{\corref}[1]{Corollary~\ref{#1}}
\newcommand{\defref}[1]{Definition~\ref{#1}}
\newcommand{\exaref}[1]{Example~\ref{#1}}

\newcommand{\lemref}[1]{Lemma~\ref{#1}}
\newcommand{\propref}[1]{Proposition~\ref{#1}}
\newcommand{\remref}[1]{Remark~\ref{#1}}
\newcommand{\secref}[1]{Section~\ref{#1}}

\begin{document}
\bibliographystyle{amsalpha}
\title{A construction of quotient $A_\infty$-categories}
\author{Volodymyr Lyubashenko%
\thanks{Institute of Mathematics,
National Academy of Sciences of Ukraine,
3 Tereshchenkivska st.,
Kyiv-4, 01601 MSP,
Ukraine,
lub@imath.kiev.ua}
\thanks{The research of V.~L. was supported in part by grant 01.07/132
of State Fund for Fundamental Research of Ukraine}
\ and Sergiy Ovsienko%
\thanks{Department of Algebra,
Faculty of Mechanics and Mathematics,
Kyiv Taras Shevchenko University,
64 Volodymyrska st.,
Kyiv, 01033,
Ukraine,
ovsienko@zeos.net}
}
\maketitle

\begin{abstract}
We construct an \ainf-category $\Dr(\cc|\cb)$ from a given
\ainf-category $\cc$ and its full subcategory $\cb$. The construction
is similar to a particular case of Drinfeld's construction of the
quotient of differential graded categories \cite{Drinf:DGquot}. We use
$\Dr(\cc|\cb)$ to construct an \ainf-functor of K\n-injective
resolutions of a complex, when the ground ring is a field. The
conventional derived category is obtained as the 0\n-th cohomology of
the quotient of the differential graded category of complexes over
acyclic complexes. This result follows also from Drinfeld's theory of
quotients of differential graded categories \cite{Drinf:DGquot}.
\end{abstract}

\allowdisplaybreaks[1]

In \cite{Drinf:DGquot} Drinfeld reviews and develops Keller's
construction of the quotient of differential graded categories
\cite{MR99m:18012} and gives a new construction of the quotient. This
construction consists of two parts. The first part replaces given pair
$\cb\subset\cc$ of a differential graded category $\cc$ and its full
subcategory $\cb$ with another such pair
$\widetilde{\cb}\subset\widetilde{\cc}$, where $\widetilde{\cc}$ is
homotopically flat over the ground ring $\kk$ (K\n-flat)
\cite[Section~3.3]{Drinf:DGquot}, and there is a quasi-equivalence
$\widetilde{\cc}\to\cc$ \cite[Section~2.3]{Drinf:DGquot}. The first
step is not needed, when $\cc$ is already homotopically flat, for
instance, when $\kk$ is a field. In the second part a new differential
graded category $\cc/\cb$ is produced from a given pair
$\cb\subset\cc$, by adding to $\cc$ new morphisms $\epsilon_U:U\to U$
of degree $-1$ for every object $U$ of $\cb$, such that
$d(\epsilon_U)=\id_U$.

In the present article we study an \ainf-analogue of the second part of
Drinfeld's construction. Namely, to a given pair $\cb\subset\cc$ of an
\ainf-category $\cc$ and its full subcategory $\cb$ we associate
another \ainf-category $\Dr(\cc|\cb)$ via a construction related to the
bar resolution of $\cc$. The \ainf-category $\Dr(\cc|\cb)$ plays the
role of the quotient of $\cc$ over $\cb$ in some cases, for instance,
when $\kk$ is a field. When $\cc$ is a differential graded category,
$\Dr(\cc|\cb)$ is precisely the category $\cc/\cb$ constructed by
Drinfeld \cite[Section~3.1]{Drinf:DGquot}.

There exists another construction of quotient \ainf-category
\(\quo(\cc|\cb)\), see \cite{math.CT/0306018}. It enjoys some universal
property, and is significantly bigger in size than $\Dr(\cc|\cb)$.
However, when $\cc$ is unital, the two quotient constructions turn out
to be \ainf-equivalent. When $\cc$ is strictly unital, so is
$\Dr(\cc|\cb)$, while \(\quo(\cc|\cb)\) is unital, but not strictly
unital.

We apply our construction to the case of practical interest:
$\cc=\Com(\ca)$ is the differential graded category of complexes in an
Abelian category $\ca$, and $\cb=\Acyc(\ca)$ is the full subcategory of
acyclic complexes. When the full subcategory $\ci\subset\cc$ of
K\n-injective complexes is big enough (every complex has a right
K\n-injective resolution) and $\kk$ is a field, we obtain an
\ainf-functor $i:\cc\to\ci$, which assigns a K\n-injective resolution
to a complex. Using it we get a new proof of the already known result:
$H^0\bigl(\Dr(\Com(\ca)|\Acyc(\ca))\bigr)$ is equivalent to the derived
category $\Dc(\ca)$ of $\ca$.

\subsection*{Outline of the article with comments.}
In the first section we describe conventions and notations used in the
article. In particular, we recall some conventions and useful formulas
from \cite{Lyu-AinfCat}.

In the second section we describe a construction of the quotient
\ainf-category $\Dr(\cc|\cb)$, departing from an \ainf-category $\cb$
fully embedded into an \ainf-category $\cc$. The underlying quiver of
$\Dr(\cc|\cb)$ is described in \defref{def-Dr-cc-cb}. Its particular
case $\Dr(\cc|\cc)$ is $s^{-1}T^+s\cc=\oplus_{n>0}s^{-1}T^ns\cc$, where
$s\cc=\cc[1]$ stands for the suspended quiver $\cc$. We introduce two
\ainf-category structures for $s^{-1}T^+s\cc$. The first,
$\underline{\cc}=(s^{-1}T^+s\cc,\underline{b})$ uses the differential
$\underline{b}$, whose components all vanish except
$\underline{b}_1=b:T^+s\cc\to T^+s\cc$, which is the \ainf-structure of
$\cc$. The second, $\overline{\cc}=(s^{-1}T^+s\cc,\overline{b})$ is
isomorphic to the first via a coalgebra automorphism
$\bm:T(T^+s\cc)\to T(T^+s\cc)$, whose components are multiplications in the
tensor algebra $T^+s\cc$. The resulting differential
$\bar{b}=\bm\underline{b}\bm^{-1}:T(T^+s\cc)\to T(T^+s\cc)$ is
described componentwise in \propref{pro-b-bar-diff}. The subquiver
$\Dr(\cc|\cb)\subset\overline{\cc}$ turns out to be an
\ainf-subcategory (\propref{pro-b-bar-diff}).
$\underline{\cc}$ and $\overline{\cc}$ are, in a sense, trivial (they
are contractible if $\cc$ is unital), but $\Dr(\cc|\cb)$, in general,
is not trivial. If $\cc$ is strictly unital, then so are
$\overline{\cc}$ and $\Dr(\cc|\cb)$ and their units are identified
(\secref{sec-str-unit}). Notice that $\underline{\cc}$ is never
strictly unital except when $\cc=0$. Nevertheless, $\underline{\cc}$
can be unital. When $\cb$, $\cc$ are differential graded categories, we
show in \secref{sec-dg-cat}, that $\Dr(\cc|\cb)$ coincides with the
category $\cc/\cb$ defined by Drinfeld
\cite[Section~3.1]{Drinf:DGquot}.

In the third section we construct functors between the obtained
\ainf-categories. When $\cb\subset\cc$, $\cj\subset\ci$ are full
\ainf-subcategories, and $i:\cc\to\ci$ is an \ainf-functor which maps
objects of $\cb$ into objects of $\cj$, we construct a strict
\ainf-functor $\underline{i}:\underline{\cc}\to\underline{\ci}$, whose
first component $\underline{i}_1=i:T^+s\cc\to T^+s\ci$ is given by $i$
itself. The components of the conjugate \ainf-functor
$\iu=\bm\underline{i}\bm^{-1}:\overline{\cc}\to\overline{\ci}$ are
described in \propref{prop-In-ce-cf}. It turns out that $\iu$ restricts
to an \ainf-functor $\Dr(i)=\iu:\Dr(\cc|\cb)\to\Dr(\ci|\cj)$
(\propref{prop-In-ce-cf}). If $\cc$ is strictly unital, then
$\underline{\cc}$ is unital (and contractible) and its unit
transformation is computed in \secref{sec-iso-strict-uni}.

In the fourth section we construct \ainf-transformations between
functors obtained in the third section. When $\cb\hookrightarrow\cc$
and $\cj\hookrightarrow\ci$ are full \ainf-subcategories,
$f,g:\cc\to\ci$ are \ainf-functors, which map objects of $\cb$ into
objects of $\cj$, and $r:f\to g:\cc\to\ci$ is an \ainf-transformation,
we construct an \ainf-transformation
$\underline{r}:\underline{f}\to\underline{g}:
\underline{\cc}\to\underline{\ci}$,
whose only non-trivial components are $\underline{r}_0=r_0$ and
$\underline{r}_1=r\big|_{T^+s\cc}$. The components of the conjugate
\ainf-transformation
\(\overline{r}=\bm\underline{r}\bm^{-1}:\overline{f}\to\overline{g}:
\overline{\cc}\to\overline{\ci}\)
are computed in \propref{prop-over-r-compon}. It turns out that
$\overline{r}$ restricts to an \ainf-transformation
$\Dr(r):\Dr(f)\to\Dr(g):\Dr(\cc|\cb)\to\Dr(\ci|\cj)$
(\propref{prop-over-r-compon}). Thus, $\unmi$ and $\ovmi$ are defined
as maps on objects, 1\n-morphisms and 2\n-morphisms of the 2\n-category
$A_\infty$ of \ainf-categories. Actually, they are strict 2\n-functors
$A_\infty\to A_\infty$ (\corref{cor-2fun-unmi-A-A},
\corref{cor-2fun-ovmi-A-A}). We prove more: they are strict
$\ck$-2-functors $\ck A_\infty\to\ck A_\infty$, where the 2\n-category
$\ck A_\infty$ is enriched in $\ck$ --  the category of differential
graded $\kk$\n-modules, whose morphisms are chain maps modulo homotopy
(\propref{pro-K2fun-unmi}, \corref{cor-K2fun-ovmi}). Compatibility of
$\ovmi$ with the composition of 2\n-morphisms is expressed via explicit
homotopy \eqref{eq-ovrpB2-ovrovpovB2}. Components of this homotopy are
found in \propref{pro-homo-rpR-explicit}. It turns out that this
homotopy restricts to subcategories $\Dr(\text-|\text-)$
(\propref{pro-homo-rpR-explicit}). Therefore, $\Dr$ is a
$\ck$-2-functor from the non-2-unital $\ck$-2-category of pairs
(\ainf-category, full \ainf-subcategory) to $\ck A_\infty$
(\corref{cor-D-K2-functor}). It can be viewed also as a 2\n-functor
$\Dr$ from the non-2-unital 2\n-category of pairs (\ainf-category, full
\ainf-subcategory) to $A_\infty$ (\corref{cor-D-2-functor-Ap-A}).

In the fifth section we consider unital \ainf-categories and prove that
some of our \ainf-categories are contractible. If $\cb$ is a full
subcategory of a unital \ainf-category $\cc$, then $\Dr(\cc|\cb)$ is
unital as well, and $\Dr(\uni^\cc)$ is its unit transformation
(\propref{pro-C-unital-D(C|B)-unital}). In particular, for a unital
\ainf-category $\cc$, both $\overline\cc$ and $\underline\cc$ are
unital with the unit transformation $\overline{\uni^\cc}$ (resp.
$\underline{\uni^\cc}$) (Corollaries \ref{cor-overline-unital},
\ref{cor-underline-unital}). If $i:\cc\to\ci$ is a unital
\ainf-functor, then $\iu:\overline{\cc}\to\overline{\ci}$,
$\underline{i}:\underline{\cc}\to\underline{\ci}$ and (whenever
defined) $\Dr(i):\Dr(\cc|\cb)\to\Dr(\ci|\cj)$ are unital as well
(Corollaries \ref{cor-overline-underline-functor-unital},
\ref{cor-D-functor-unital}). When we restrict $\unmi$, $\ovmi$ or $\Dr$
to unital \ainf-categories (and unital \ainf-functors), we get strict
2\n-functors of (usual 1-2-unital) ($\ck$\n-)2-categories.

In the sixth section we consider contractible \ainf-categories and
\ainf-functors. A unital \ainf-functor $f:\ca\to\cb$ is called
\emph{contractible} if many equivalent conditions hold, including
contractibility of complexes $(s\cb(Xf,Y),b_1)$, $(s\cb(Y,Xf),b_1)$ for
all $X\in\Ob\ca$, $Y\in\Ob\cb$ (Propositions
\ref{pro-C1-C5-contract}--\ref{pro-C10-C11-unital-f},
\defref{def-contractible-functor-category}). A unital \ainf-category
$\ca$ is called \emph{contractible} if several equivalent conditions
hold, including contractibility of complexes $(s\ca(X,Y),b_1)$ for all
objects $X$, $Y$ of $\ca$ (\defref{def-contractible-functor-category},
\propref{pro-C0-C10'-contractible-category}). If $\cc$ is a unital
\ainf-category, then $\underline\cc$, $\overline\cc$ are contractible
(\exaref{exam-undc-overc-contractible}). Nevertheless, in general, the
subcategories $\Dr(\cc|\cb)\subset\overline\cc$ are not contractible.
Contractible \ainf-categories $\cb$ may be considered as trivial,
because in this case any natural \ainf-transformation
$r:f\to g:\ca\to\cb$ is equivalent to 0
(\corref{cor-contractible-r-equivalent-0}). Moreover, non-empty
contractible categories are equivalent to the 1-object-1-morphism
\ainf-category $\1$, such that $\Ob\1=\{*\}$ and $\1(*,*)=0$
(\propref{pro-C0-C10'-contractible-category},
\remref{rem-O-0-equivalent-1}).

In the seventh section we consider the case of a contractible full
subcategory $\cf$ of a unital \ainf-category $\ce$. In this case the
canonical strict embedding $\ce\to\Dr(\ce|\cf)$ is an equivalence
(\propref{pro-D(E|contractible-F)-E}).

In the eighth section we prepare to construct the K-injective
resolution \ainf-functor. This concrete construction is deferred until
the next section. In the eighth section we consider an abstract version
of it. Given an \ainf-functor $f:\cb\to\cc$, a map $g:\Ob\cb\to\Ob\cc$
and cycles $r_X\in\cc^0(Xf,Xg)$, $X\in\Ob\cb$, producing certain
quasi-isomorphisms, we make $g$ into an \ainf-functor $g:\cb\to\cc$ and
$r_X$ into 0\n-th component $\sS{_X}r_0s^{-1}$ of a natural
\ainf-transformation $r:f\to g:\cb\to\cc$
(\propref{prop-ainf-fr-g-tr-r}). Next we prove the uniqueness of so
constructed $g$ and $r$. Assuming that the initial data
$(g:\Ob\cb\to\Ob\cc,(r_X)_{x\in\Ob\cb})$ give rise to two
\ainf-functors $g,g':\cb\to\cc$ and two natural \ainf-transformations
$r:f\to g:\cb\to\cc$, $r':f\to g':\cb\to\cc$, we construct another
natural \ainf-transformation $p:g\to g':\cb\to\cc$, such that $r'$ is
the composition $(f \rTTo^r g \rTTo^p g')$ in the 2\n-category
$A_\infty$ (\propref{pro-p-g-g-exists}). Moreover, such $p$ is unique
up to an equivalence (\propref{prop-p-unique}). If, in addition, $\cc$
is unital, then the constructed $p$ is invertible
(\corref{cor-p-invertible}). If $f$ is unital, then the constructed
\ainf-functor $g$ is unital as well (\propref{prop-g-unital}).

In the ninth section we consider categories of complexes. Let $\kk$ be
a field, let $\ca$ be an Abelian $\kk$\n-linear category, and let
$\cc=\Com(\ca)$ be the differential graded category of complexes in $\ca$.
Let $\cb=\Acyc(\ca)$ be its full subcategory of acyclic complexes,
$\ci=\Inj(\ca)$ denotes K\n-injective complexes, $\cj=\AcIn(\ca)$
denotes acyclic K\n-injective complexes. We assume that each complex
$X\in\Ob\cc$ has a right K\n-injective resolution $r_X:X\to Xi$ (a
quasi-isomorphism with K\n-injective $Xi\in\Ob\ci$). We notice that
quasi-isomorphisms from $\cc$ become ``invertible modulo boundary'' in
the differential graded category $\Dr(\cc|\cb)$
(\secref{sec-invert-qis}). From the identity functor
$f=\id_\cc:\cc\to\cc$, a map $g:\Ob\cc\to\Ob\cc$, $X\mapsto Xi$ and
quasi-isomorphisms $r_X$ we produce an \ainf-functor $g:\cc\to\cc$,
which factors as $g=\bigl(\cc \rTTo^i \ci \rMono^e \cc\bigr)$, into
``K-injective resolution'' unital \ainf-functor $i$ and the full
embedding $e$ (\secref{sec-K-injective-complexes}). The unital
\ainf-functor $\iu:\Dr(\cc|\cb)\to\Dr(\ci|\cj)$ and the faithful
differential graded functor $\overline{e}:\Dr(\ci|\cj)\to\Dr(\cc|\cb)$
are \ainf-equivalences quasi-inverse to each other. Due to
contractibility of $\cj$ the natural embedding $\ci\to\Dr(\ci|\cj)$ is
an equivalence. Since $\Dr(\cc|\cb)$ and $\ci$ are \ainf-equivalent,
their 0\n-th cohomology categories are equivalent as usual
$\kk$\n-linear categories. That is,
$H^0(\Dr(\cc|\cb))=H^0\bigl(\Dr(\Com(\ca)|\Acyc(\ca))\bigr)$ is
equivalent to $H^0(\ci)=H^0(\Inj(\ca))$ -- homotopy category of
K-injective complexes, which is equivalent to the derived category
$\Dc(\ca)$ of $\ca$. This result (\secref{sec-K-injective-complexes})
motivated our studies. It follows also from Drinfeld's theory of
quotients of differential graded categories \cite{Drinf:DGquot}. This
agrees with Bondal and Kapranov's proposal to produce triangulated
categories as homotopy categories of some differential graded
categories \cite{BondalKapranov:FramedTR}.

\section{Conventions}\label{sec-convent-nota}
We keep the notations and conventions of \cite{Lyu-AinfCat}, sometimes
without explicit mentioning. Some of the conventions are recalled here.

We assume as in \cite{Lyu-AinfCat} that quivers, \ainf-categories, etc.
are small with respect to some universe $\fu$.

The ground ring $\kk\in\fu$ is a unital associative commutative ring.

We use the right operators: the composition of two maps (or morphisms)
$f:X\to Y$ and $g:Y\to Z$ is denoted $fg:X\to Z$; a map is written on
elements as $f:x\mapsto xf=(x)f$. However, these conventions are not
used systematically, and $f(x)$ might be used instead.

If $C$ is a $\ZZ$\n-graded $\kk$\n-module, then $sC=C[1]$ denotes the
same $\kk$\n-module with the grading $(sC)^d=C^{d+1}$, the suspension
of $C$. The shift ``identity'' map $C\to sC$ of degree $-1$ is also
denoted $s$. Getzler and Jones demonstrated in
\cite{GetzlerJones:A-infty} that the suspension $s$ and the shift map
$s$ are useful in the theory of \ainf-algebras. We follow the Koszul
sign convention:
\[ (x\tens y)(f\tens g)=(-)^{yf}xf\tens yg
=(-1)^{\deg y\cdot\deg f}xf\tens yg. \]

A chain complex is called \emph{contractible} if its identity
endomorphism is homotopic to zero.

The category $\cQuiver/S$ of $\fu$\n-small graded $\kk$\n-linear
quivers with fixed set of objects $S$ admits a monoidal
structure with the tensor product $\ca\times\cb\mapsto\ca\tens\cb$,
$(\ca\tens\cb)(X,Y)=\oplus_{Z\in S}\ca(X,Z)\tens_\kk\cb(Z,Y)$. Thus, we
have tensor powers $T^n\ca=\ca^{\tens n}$ of a given graded
$\kk$\n-quiver $\ca$, such that $\Ob T^n\ca=\Ob\ca$. Explicitly,
\[ T^n\ca(X,Y) = \bigoplus_{X_1,\dots,X_{n-1}\in\Ob\ca}
\ca(X_0,X_1)\tens_\kk\ca(X_1,X_2)\tens_\kk\dots\tens_\kk\ca(X_{n-1},X_n),
\]
where $X_0=X$ and $X_n=Y$. In particular, \(T^0\ca\) denotes the unit
object $\kk S$, where $\kk S(X,Y)=\kk$ if $X=Y$ and vanishes otherwise.

As in any monoidal category, there is a notion of coassociative
coalgebras \((\cb,\overline\Delta:\cb\to\cb\tens\cb)\) in $\cQuiver/S$
(in general, without counit). For \ainf-category theory, we need only
coalgebras \((\cb,\overline\Delta)\) in $\cQuiver/S$ that satisfy an
additional requirement: for all \(X,Y\in S\)
\[ \cb(X,Y) = \cup_{k=2}^\infty
\Ker\bigl(\overline\Delta^{(k)}:\cb(X,Y)\to\cb^{\tens k}(X,Y)\bigr),
\]
where $\overline\Delta^{(2)}=\overline\Delta$,
 $\overline\Delta^{(3)}=\overline\Delta(1\tens\overline\Delta)
 =\overline\Delta(\overline\Delta\tens1):\cb\to\cb^{\tens3}$,
etc. Such coalgebras are named \emph{cocomplete cocategories} by
Keller~\cite{math.RT/0510508}. A counital coassociative coalgebra
$(\ca=T^0\cb\oplus\cb,\Delta:\ca\to\ca\tens\ca,\eps:\ca\to T^0\ca)$ in
$\cQuiver/S$ is associated with \((\cb,\overline\Delta)\), namely:
\begin{align*}
\Delta\big|_{T^0\cb} &= \bigl(T^0\cb \rTTo^\sim T^0\cb\tens T^0\cb
\rMono^{\hspace*{.3em}\inj_0\tens\inj_0\hspace*{.2em}}
\ca\tens\ca\bigr),
\\
\Delta\big|_\cb &= \bigl(\cb \rTTo^\sim \cb\tens T^0\cb
\rMono^{\hspace*{.3em}\inj_1\tens\inj_0\hspace*{.2em}}
\ca\tens\ca\bigr)
\\
&+ \bigl(\cb \rTTo^{\overline\Delta} \cb\tens\cb
\rMono^{\hspace*{.3em}\inj_1\tens\inj_1\hspace*{.2em}}
\ca\tens\ca\bigr)
\\
&+ \bigl(\cb \rTTo^\sim T^0\cb\tens\cb
\rMono^{\hspace*{.3em}\inj_0\tens\inj_1\hspace*{.2em}}
\ca\tens\ca\bigr),
\end{align*}
or simply \(f\Delta=f\tens1+f\overline\Delta+1\tens f\) for
\(f\in\cb(X,Y)\), and \(\eps=\pr_0:\ca\to T^0\cb=T^0\ca\). The triple
\((\ca,\Delta,\epsilon)\) is a \emph{cocategory} in the sense of
\cite{Lyu-AinfCat}. In the present article, we shall use only one kind
of cocategory associated with quivers $\cc$, namely, the cocomplete
cocategory \(\cb=T^+\cc=\oplus_{n=1}^\infty T^n\cc\), equipped with the
comultiplication \(\overline\Delta:\cb\to\cb\tens\cb\),
 $(h_1\tens h_2\tens\dots\tens h_n)\overline\Delta=\sum_{k=1}^{n-1}
 h_1\tens\dots\tens h_k\bigotimes h_{k+1}\tens\dots\tens h_n$,
gives rise to the cocategory \(\ca=T\cc=\oplus_{n=0}^\infty T^n\cc\),
equipped with the cut comultiplication \(\Delta:\ca\to\ca\tens\ca\),
 $(h_1\tens h_2\tens\dots\tens h_n)\Delta=\sum_{k=0}^n
 h_1\tens\dots\tens h_k\bigotimes h_{k+1}\tens\dots\tens h_n$,
and with the counit \(\eps=\pr_0:\ca\to T^0\cc=T^0\ca\).

By definition, cocategory homomorphisms (in particular, \ainf-functors)
respect the cut comultiplication $\Delta$, and \ainf-transformations
are coderivations with respect to $\Delta$ (see e.g.
\cite{Lyu-AinfCat}).

We use the following standard equations for a differential $b$ in an
\ainf-category
\begin{equation}
\sum_{r+n+t=k} (1^{\tens r}\tens b_n\tens1^{\tens t})b_{r+1+t} = 0
: T^ks\ca \to s\ca.
\label{eq-b-b-0}
\end{equation}
Since $b$ is a differential and a coderivation, it may be called a
\emph{codifferential}. Commutation relation $fb=bf$ for an
\ainf-functor $f:\ca\to\cb$ expands to the following
\begin{equation}
\sum_{l>0;i_1+\dots+i_l=k}(f_{i_1}\tens f_{i_2}\tens\dots\tens f_{i_l})b_l
= \sum_{r+n+t=k} (1^{\tens r}\tens b_n\tens1^{\tens t}) f_{r+1+t}
: T^ks\ca \to s\cb.
\label{eq-fff-b-b-f}
\end{equation}

Given \ainf-functors $f,g,h:\cb\to\cc$ and coderivations
$f\rTTo^r g\rTTo^p h:\cb\to\cc$ of arbitrary degree we construct a map
$\theta:Ts\cb\to Ts\cc$ as in Section~3 of \cite{Lyu-AinfCat}. We view
$\theta$ as a bilinear function $(r\tens p)\theta$ of $r$, $p$. Its
components $\theta_{kl}=\theta\big|_{T^ks\cb}\pr_l:T^ks\cb\to T^ls\cc$
are given by formula (3.0.1) of \cite{Lyu-AinfCat}
\begin{equation}
\theta_{kl} = \sum f_{a_1}\tens\dots\tens f_{a_\alpha}\tens
r_j\tens g_{c_1}\tens\dots\tens g_{c_\beta}\tens
p_t\tens h_{e_1}\tens\dots\tens h_{e_\gamma},
\label{eq-theta-kl-intro}
\end{equation}
where the summation is taken over all terms with
\[ \alpha+\beta+\gamma+2=l, \qquad
a_1+\dots+a_\alpha+j+c_1+\dots+c_\beta+t+e_1+\dots+e_\gamma = k. \]
The same formula can be presented as
\begin{equation}
\theta_{kl} =
\sum_{\substack{\alpha+\beta+\gamma+2=l\\a+j+c+t+e=k}}
f_{a\alpha}\tens r_j\tens g_{c\beta}\tens p_t\tens h_{e\gamma},
\label{eq-theta-kl-short}
\end{equation}
where $f_{a\alpha}:T^a\ca\to T^\alpha\cb$ are matrix elements of $f$
and similarly for $g$, $h$. By Proposition~3.1 of \cite{Lyu-AinfCat}
the map $\theta$ satisfies the equation
\[ \theta\Delta = \Delta(f\tens\theta + r\tens p + \theta\tens h). \]

Given \ainf-categories $\ca$ and $\cb$, one constructs an
\ainf-category $A_\infty(\ca,\cb)$ of \ainf-functors $\ca\to\cb$,
equipped with a differential $B$
\cite{Fukaya:FloerMirror-II,math.RT/0510508,math.RA/0606241,KonSoi-book,Lefevre-Ainfty-these},
\cite[Section~5]{Lyu-AinfCat}.

The category of graded $\kk$\n-linear quivers admits a symmetric
monoidal structure with the tensor product
$\ca\times\cb\mapsto\ca\boxtimes\cb$, where
$\Ob\ca\boxtimes\cb=\Ob\ca\times\Ob\cb$ and
$(\ca\boxtimes\cb)\bigl((X,U),(Y,V)\bigr)=\ca(X,Y)\tens_\kk\cb(U,V)$.
The same tensor product was denoted $\tens$ in \cite{Lyu-AinfCat}.
Given \ainf-categories $\ca$, $\cb$, $\cc$, there is a graded
cocategory morphism of degree 0
\[ M: TsA_\infty(\ca,\cb)\boxtimes TsA_\infty(\cb,\cc)
\to TsA_\infty(\ca,\cc), \]
which satisfies equation $(1\boxtimes B+B\boxtimes1)M=MB$
\cite[Section~6]{Lyu-AinfCat}.

\section{\texorpdfstring{An $A_\infty$-category}{An A8-category}}
Let $\cb\hookrightarrow\cc$ be a full \ainf-subcategory. It means that
$\Ob\cb\subset\Ob\cc$, $\cb(X,Y)=\cc(X,Y)$ for all $X,Y\in\Ob\cb$, and
the operations for $\cb$ coincide with those for $\cc$. Let us define
another \ainf-category $\Dr(\cc|\cb)$. If $\cb$, $\cc$ are differential
graded categories, then $\Dr(\cc|\cb)$ is differential graded as well
and it coincides with the category $\cc/\cb$ defined by Drinfeld in
\cite[Section~3.1]{Drinf:DGquot}.

\begin{definition}\label{def-Dr-cc-cb}
Let $T^+s\cc=\oplus_{n>0}T^ns\cc$ and $\ce=\Dr(\cc|\cb)$ be the
following graded $\kk$\n-quivers: the class of objects is
$\Ob T^+s\cc=\Ob\ce=\Ob\cc$, the morphisms for $X,Y\in\Ob\ce$ are
\begin{align*}
T^+s\cc(X,Y) &= \oplus_{C_1,\dots,C_{n-1}\in\cc} s\cc(X,C_1)\tens
s\cc(C_1,C_2)\tens\dots\tens s\cc(C_{n-2},C_{n-1})\tens
s\cc(C_{n-1},Y), \\
s\ce(X,Y) &= \oplus_{C_1,\dots,C_{n-1}\in\cb} s\cc(X,C_1)\tens
s\cc(C_1,C_2)\tens\dots\tens s\cc(C_{n-2},C_{n-1})\tens s\cc(C_{n-1},Y),
\end{align*}
where in the second case summation extends over all sequences of
objects $(C_1,\dots,C_{n-1})$ of $\cb$. To the empty sequence ($n=1$)
corresponds the summand $s\cc(X,Y)$.
\end{definition}

Let us endow $s^{-1}T^+s\cc$ with a structure of \ainf-category, given
by $\underline{b}:T(T^+s\cc)\to T(T^+s\cc)$, with the components
$\underline{b}_0=0$, $\underline{b}_1=b:T^+s\cc\to T^+s\cc$,
$\underline{b}_k=0$ for $k>1$. This \ainf-category is denoted
$\underline{\cc}=(s^{-1}T^+s\cc,\underline{b})$. There is an
\ainf-functor $\underline{j}:\cc\to(s^{-1}T^+s\cc,\underline{b})$,
specified by its components $\underline{j}_k:T^ks\cc\to T^+s\cc$,
$k\ge1$, where $\underline{j}_k$ is the canonical embedding of the
direct summand. The property
$b\underline{j}=\underline{j}\underline{b}$, or
\[ \sum_{r+k+t=n}(1^{\tens r}\tens b_k\tens 1^{\tens t})
\underline{j}_{r+1+t} = \underline{j}_nb:T^ns\cc\to T^+s\cc,
\]
is clear -- this is just the expression of $b$ in terms of its components.

There is a coalgebra automorphism $\bm:TT^+s\cc\to TT^+s\cc$, specified
by its components $\bm_k=\mu^{(k)}:T^kT^+s\cc\to T^+s\cc$, $k\ge1$,
where $\mu:T^+s\cc\tens T^+s\cc\to T^+s\cc$ is the multiplication in
the tensor algebra, $\mu^{(k)}=0$ for $k\le0$,
$\mu^{(1)}=1:T^+s\cc\to T^+s\cc$, $\mu^{(2)}=\mu$,
$\mu^{(3)}=(\mu\tens1)\mu:(T^+s\cc)^{\tens3}\to T^+s\cc$ and so on. Its
inverse is the coalgebra automorphism
$\bm^{-1}=\bm^-:TT^+s\cc\to TT^+s\cc$, specified by its components
$\bm^-_k=(-)^{k-1}\mu^{(k)}:T^kT^+s\cc\to T^+s\cc$. The fact that $\bm$
and $\bm^-$ are inverse to each other is proven as follows:
\begin{multline*}
(\bm\bm^-)_n
= \sum_{l_1+\dots+l_k=n}(\bm_{l_1}\tens\dots\tens\bm_{l_k})\bm^-_k
= \sum_{l_1+\dots+l_k=n} (-)^{k-1}\mu^{(n)} \\
= \mu^{(n)}\sum_{k=1}^n (-1)^{k-1}\binom{n-1}{k-1}
= \mu^{(n)}(1-1)^{n-1},
\end{multline*}
which equals $\id$ for $n=1$ and vanishes for $n>1$. Similarly,
$\bm^-\bm=\id$.

\begin{proposition}\label{pro-b-bar-diff}
The conjugate codifferential
$\bar{b}=\bm\underline{b}\bm^{-1}:T(T^+s\cc)\to T(T^+s\cc)$ has the
following components: $\overline{b}_0=0$, $\overline{b}_1=b$ and for
$n\ge2$
\begin{align}
\bar{b}_n &= \mu^{(n)}
\sum_{m;q<k;t<l} 1^{\tens q}\tens b_m\tens1^{\tens t}
:T^ks\cc\tens(T^+s\cc)^{\tens n-2}\tens T^ls\cc \to T^+s\cc,
\label{eq-Bn-sum-1bm1}
\\
\bar{b}_n &= \mu^{(n)}b -(1\tens\mu^{(n-1)}b)\mu -(\mu^{(n-1)}b\tens1)\mu
+ (1\tens\mu^{(n-2)}b\tens1)\mu^{(3)} :(T^+s\cc)^{\tens n}\to T^+s\cc,
\label{eq-Bn-mub-4terms}
\end{align}
for all $n\ge0$. The operations $\bar{b}_n$ restrict to maps
$s\ce^{\tens n}\to s\ce$ via the natural embedding
$s\ce\subset T^+s\cc$ of graded $\kk$\n-quivers. Hence, $\bar{b}$ turns
$\ce$ and $\overline{\cc}\overset{\text{def}}=(s^{-1}T^+s\cc,\bar{b})$
into an \ainf-category.
\end{proposition}

\begin{proof}
Let us define a $(1,1)$-coderivation $\bar{b}$ of degree 1 by its
components~\eqref{eq-Bn-mub-4terms}. Substituting the definition of $b$
via its components, we get formula~\eqref{eq-Bn-sum-1bm1}. Clearly,
$\bm\underline{b}\bm^{-1}$ is also a $(1,1)$-coderivation of degree 1.
Let us show that it coincides with $\bar{b}$, that is,
$\bar{b}\bm=\bm\underline{b}$. This equation expands to
\[ \sum_{r+k+t=n}(1^{\tens r}\tens\bar{b}_k\tens1^{\tens t})\mu^{(r+1+t)}
= \mu^{(n)}b = \mu^{(n)}\underline{b}_1,
\]
which follows immediately from \eqref{eq-Bn-sum-1bm1} and from the
standard expression of $b$ via its components.

Clearly, $\bar{b}^2=\bm\underline{b}^2\bm^{-1}=0$, hence,
$\overline{\cc}=(s^{-1}T^+s\cc,\bar{b})$ is an \ainf-category.
Map~\eqref{eq-Bn-sum-1bm1} is a sum of maps of the form
\begin{multline*}
1^{\tens q}\tens b_m\tens1^{\tens t}:
s\cc(X,C_1)\tens\dots\tens s\cc(C_{q-1},C_q)\tens s\cc(C_q,C_{q+1})
\tens\dots\tens(T^+s\cc)^{\tens n-2}\tens \\
\dots\tens s\cc(D_{l-t-1},D_{l-t})\tens
s\cc(D_{l-t},D_{l-t+1})\tens\dots\tens s\cc(D_{l-1},Y) \\
\to s\cc(X,C_1)\tens\dots\tens s\cc(C_{q-1},C_q)\tens
s\cc(C_q,D_{l-t})\tens s\cc(D_{l-t},D_{l-t+1})\tens\dots\tens
s\cc(D_{l-1},Y),
\end{multline*}
where $C_0=X$ for $q=0$ and $D_l=Y$ for $t=0$. If the source is
contained in $(s\ce)^{\tens n}(X,Y)$, then $C_i$, $D_j$ are in $\Ob\cb$
for all $0<i<k$, $0<j<l$. Therefore, the target is a direct summand of
$s\ce(X,Y)$. Thus the maps $\bar{b}_n$ restrict to maps
$\bar{b}_n:(s\ce)^{\tens n}(X,Y)\to s\ce(X,Y)$. The obtained
(1,1)-coderivation $\bar{b}:Ts\ce\to Ts\ce$ also satisfies
$\bar{b}^2=0$. Thus it makes $\ce$ into an \ainf-category.
\end{proof}

In particular, \eqref{eq-Bn-mub-4terms} gives
\begin{align*}
\bar{b}_2 &= \mu b - (1\tens b + b\tens1)\mu, \\
\bar{b}_3 &= \mu^{(3)}b - (1\tens\mu b)\mu - (\mu b\tens1)\mu
+ (1\tens b\tens1)\mu^{(3)}, \\
\bar{b}_4 &= \mu^{(4)}b - (1\tens\mu^{(3)}b)\mu - (\mu^{(3)}b\tens1)\mu
+ (1\tens\mu b\tens1)\mu^{(3)}.
\end{align*}

\begin{remark}\label{rem-move-diff}
Let $\ca$ be an \ainf-category, defined by a codifferential
$b:Ts\ca\to Ts\ca$, let $\cb$ be a graded $\kk$\n-quiver and let
$f:Ts\ca\to Ts\cb$ (resp. $g:Ts\cb\to Ts\ca$) be an isomorphism of
graded cocategories. Then the codifferential $f^{-1}bf$ (resp.
$gbg^{-1}$) is the unique codifferential on $Ts\cb$, which turns $f$
(resp. $g$) into an invertible \ainf-functor between $\ca$ and $\cb$.
\end{remark}

\begin{corollary}\label{cor-mu-ainf-functor}
The coalgebra isomorphism
 $\bm^{-1}:\underline{\cc}=(s^{-1}T^+s\cc,\underline{b})
 \to\overline{\cc}=(s^{-1}T^+s\cc,\bar{b})$
is an \ainf-functor. Its composition with $\underline{j}$ is a strict
\ainf-functor $\ju=\underline{j}\bm^{-1}:\cc\to\Dr(\cc|\cb)$,
$X\mapsto X$, whose components are the direct summand embedding
$\ju_1:s\cc(X,Y)=T^1s\cc(X,Y)\hookrightarrow s\ce(X,Y)$ and $\ju_n=0$
for $n>1$.
\end{corollary}

Indeed,
\[ \ju_n = \sum_{l_1+\dots+l_k=n}
(\underline{j}_{l_1}\tens\dots\tens\underline{j}_{l_k})(-1)^{k-1}\mu^{(k)}
= \sum_{k=1}^n (-1)^{k-1}\binom{n-1}{k-1}\underline{j}_n
= (1-1)^{n-1}\underline{j}_n, \]
which equals $\underline{j}_1$ for $n=1$ and vanishes for $n>1$.

\subsection{Strict unitality.}\label{sec-str-unit}
Assume that \ainf-category $\cc$ is strictly unital. It means that for
each object $X$ of $\cc$ there is an element $1_X\in\cc^0(X,X)$, such
that the map $\uni^\cc_0:\kk\to(s\cc)^{-1}(X,X)$, $1\mapsto1_Xs$ of
degree $-1$ satisfies equations
$(1\tens \uni^\cc_0)b_2=1:s\cc(Y,X)\to s\cc(Y,X)$ and
$(\uni^\cc_0\tens1)b_2=-1:s\cc(X,Z)\to s\cc(X,Z)$ for all
$Y,Z\in\Ob\cc$, and $(\dots\tens1_Xs\tens\dots)b_n=0$ if $n\ne2$. Since
$\cc$ is strictly unital, its full \ainf-subcategory $\cb$ is strictly
unital as well.

Let us show that in these assumptions $\ce=\Dr(\cc|\cb)$ is also
strictly unital. We take the same elements
$1_X\in\cc^0(X,X)\subset\ce^0(X,X)$ as strict units of $\ce$. We have
$1_Xs\bar{b}_1=1_Xsb=1_Xsb_1=0$. Explicit formulas give
$(\dots\tens1_Xs\tens\dots)\bar{b}_n=0$ for $n>2$. The map
$\bar{b}_2:T^ks\cc(Y,X)\tens s\cc(X,X)\to T^+s\cc(Y,X)$ is the sum of maps
\begin{multline*}
1^{\tens k-t}\tens b_{t+1}: s\cc(Y,C_1)\tens\dots\tens
s\cc(C_{k-t},C_{k-t+1})\tens\dots\tens s\cc(C_{k-1},X)\tens s\cc(X,X) \\
\to s\cc(Y,C_1)\tens\dots\tens s\cc(C_{k-t},X)
\end{multline*}
over $t>0$. Therefore, the map $\uni^\ce_0:\kk\to(s\ce)^{-1}(X,X)$,
$1\mapsto1_Xs$ satisfies equations
$(1\tens\uni^\ce_0)\bar{b}_2=(1^{\tens k}\tens\uni^\cc_0)
(1^{\tens k-1}\tens b_2)=1^{\tens k-1}\tens1=1$.
Similarly, for $\bar{b}_2:s\cc(X,X)\tens T^ks\cc(X,Z)\to T^+s\cc(X,Z)$
we have
$(\uni^\ce_0\tens1)\bar{b}_2=(\uni^\cc_0\tens1^{\tens k})
(b_2\tens1^{\tens k-1})=-1\tens1^{\tens k-1}=-1$.
Therefore, $\ce$ and $\overline{\cc}$ are strictly unital with the unit
$\uni^\ce$.

\subsection{Differential graded categories.}\label{sec-dg-cat}
If $b_k=0$ for $k>2$, then explicit formulae in the case of $\ce$ show
that we also have $\bar{b}_k=0$ for $k>2$. Combining this fact with the
above unitality considerations, we see that if $\cc$ is a differential
graded category, then so is $\Dr(\cc|\cb)$. The differential graded
category $\ce=\Dr(\cc|\cb)=\cc/\cb$ was constructed by Drinfeld
\cite[Section~3.1]{Drinf:DGquot}. This construction was a starting
point of the present article. Let us describe it in detail.

Write down elements of $\ce(X,Y)$ as sequences
$f_1\epsilon_{C_1}f_2\dots\epsilon_{C_{n-1}}f_n$, where
$f_i\in\cc(C_{i-1},C_i)$, $C_0=X$, $C_n=Y$, and $C_i\in\Ob\cb$ for
$0<i<n$. The symbol $\epsilon_C$ for $C\in\Ob\cb$ is assigned degree
$-1$. Its differential is set equal to $\epsilon_Cd=1_C$. The graded
Leibniz rule gives
\begin{multline*}
(f_1\epsilon_{C_1}f_2\dots\epsilon_{C_{n-1}}f_n)d =
\hspace*{-2mm} \sum_{q+1+t=n} \hspace*{-2mm}(-)^{f_{n-t+1}+\dots+f_n-t}
f_1\epsilon_{C_1}f_2\dots\epsilon_{C_q}(f_{q+1}m_1)
\epsilon_{C_{q+1}}f_{q+2}\dots\epsilon_{C_{n-1}}f_n \\
+ \sum_{q+2+t=n} (-)^{f_{n-t}+\dots+f_n-t}
f_1\epsilon_{C_1}f_2\dots\epsilon_{C_q}(f_{q+1}\cdot f_{q+2})
\epsilon_{C_{q+2}}f_{q+3}\dots\epsilon_{C_{n-1}}f_n,
\end{multline*}
where $f_{q+1}\cdot f_{q+2}=(f_{q+1}\tens f_{q+2})m_2$ is the
composition. Introduce a degree $-1$ map
\[ s:\ce\to s\ce\subset T^+s\cc, \quad
f_1\epsilon_{C_1}f_2\dots\epsilon_{C_{n-1}}f_n
\mapsto f_1s\tens f_2s\tens\dots\tens f_ns. \]
One can check that $ds=s\bar{b}_1$, where, naturally,
$\bar{b}_1=b=\sum_{q+1+t=n}1^{\tens q}\tens b_1\tens1^{\tens t}
+\sum_{q+2+t=n}1^{\tens q}\tens b_2\tens1^{\tens t}$.

The composition $\bar{m}_2$ in $\ce$ consists of the concatenation and
the composition $m_2$ in $\cc$:
\begin{multline*}
(f_1\epsilon_{C_1}\dots f_{n-1}\epsilon_{C_{n-1}}f_n\tens
g_1\epsilon_{D_1}g_2\dots\epsilon_{D_{m-1}}g_m)\bar{m}_2 \\
= f_1\epsilon_{C_1}\dots f_{n-1}\epsilon_{C_{n-1}}
(f_n\cdot g_1)\epsilon_{D_1}g_2\dots\epsilon_{D_{m-1}}g_m.
\end{multline*}
One can check that $\bar{m}_2s=(s\tens s)\bar{b}_2$; here
$\bar{b}_2=1^{\tens n-1}\tens b_2\tens1^{\tens m-1}$.

Specifically this construction applies to the case of the differential
graded category $\cc=\Com(\ca)$ of complexes of objects of an abelian
category $\ca$. One may take for $\cb$ the subcategory of acyclic
complexes $\cb=\Acyc(\ca)$.

\section{\texorpdfstring{An $A_\infty$-functor}{An A8-functor}}
Let $\cb\hookrightarrow\cc$, $\cj\hookrightarrow\ci$ be full
\ainf-subcategories. Let $i:\cc\to\ci$ be an \ainf-functor, such that
$Xi\in\Ob\cj$ for $X\in\Ob\cb$. Then it restricts to an \ainf-functor
$\cb\to\cj$, denoted by $i'$. We are going to construct an extension of
this functor to the \ainf-categories $\ce=\Dr(\cc|\cb)$ and
$\cf=\Dr(\ci|\cj)$.

Let us begin with a strict \ainf-functor
$\underline{i}:\underline{\cc}\to\underline{\ci}$, given by its
components $\underline{i}_1=i:T^+s\cc\to T^+s\ci$ and
$\underline{i}_k=0$ for $k>1$. The equation
$\underline{i}\,\underline{b}=\underline{b}\,\underline{i}$ reduces to
familiar $ib=bi$. Therefore,
$\iu\overset{\text{def}}=\bm\underline{i}\bm^{-1}
:\overline{\cc}\to\overline{\ci}$
is an \ainf-functor as well.

The following diagram of \ainf-functors commutes
\begin{diagram}
\cb & \rMono & \cc & \rTTo^{\underline{j}^\cc} & \underline{\cc}
& \rTTo^{\bm^{-1}} & \overline{\cc} \\
\dTTo<{i'} && \dTTo>i && \dTTo>{\underline{i}} && \dTTo>{\iu} \\
\cj & \rMono & \ci & \rTTo^{\underline{j}^\ci} & \underline{\ci}
& \rTTo^{\bm^{-1}} & \overline{\ci}
\end{diagram}
Indeed, $\underline{j}^\cc\underline{i}=i\underline{j}^\cj$ expands to
\[ \underline{j}^\cc_ni = \sum_{l_1+\dots+l_k=n}
(i_{l_1}\tens\dots\tens i_{l_k})\underline{j}^\cj_k:T^ns\cc\to T^+s\ci,
\]
which expresses $i$ in terms of its components.

\begin{proposition}\label{prop-In-ce-cf}
The \ainf-functor $\iu$ has the following components:
\begin{equation}
\iu_n = \sum_{l_1+\dots+l_k=n} (-)^{k-1}
\bigl(\mu^{(l_1)}\tens\dots\tens\mu^{(l_k)}\bigr)i^{\tens k}\mu^{(k)}
:(T^+s\cc)^{\tens n}\to T^+s\ci.
\label{eq-In-mu-i-mu}
\end{equation}
The restriction of this map to $T^{k_1}s\cc\tens\dots\tens T^{k_n}s\cc$
is
\begin{equation}
\iu_n = \mu^{(n)} \sum_{(l_1,\dots,l_t)\in L(k_1,\dots,k_n)}
(i_{l_1}\tens\dots\tens i_{l_t}):
T^{k_1}s\cc\tens\dots\tens T^{k_n}s\cc \to T^+s\ci,
\label{eq-In-mun-iii}
\end{equation}
\begin{multline*}
L(k_1,\dots,k_n) = \cup_{t>0} \bigl\{ (l_1,\dots,l_t)\in\ZZ^t_{>0} \mid
\forall q,s\in\ZZ_{>0},q\le t,s\le n \\
l_1+\dots+l_q=k_1+\dots+k_s \Longleftrightarrow q=t,s=n \bigr\}.
\end{multline*}
These maps restrict to maps $\iu_n:T^ns\Dr(\cc|\cb)\to s\Dr(\ci|\cj)$,
which are components of an \ainf-functor
$\Dr(i)=\iu:\Dr(\cc|\cb)\to\Dr(\ci|\cj)$. The restriction of $\iu_n$
to $T^ns\cc \rMono^{\ju_1^{\tens n}} T^ns\Dr(\cc|\cb)$ equals
$T^ns\cc \rTTo^{i_n} s\ci \rMono^{\ju_1} s\Dr(\ci|\cj)$.
\end{proposition}

\begin{proof}
Let us prove \eqref{eq-In-mun-iii}. Let $\iu$ denote the cocategory
homomorphism $\iu:\overline{\cc}\to\overline{\ci}$, defined by its
components~\eqref{eq-In-mun-iii}. We are going to prove that it
satisfies $\iu\bm=\bm\underline{i}$. Indeed, this equation expands to
the following
\begin{equation}
\sum_{n_1+\dots+n_a=p}(\iu_{n_1}\tens\dots\tens\iu_{n_a})\mu^{(a)}
= \mu^{(p)}i: T^{c_1}s\cc\tens\dots\tens T^{c_p}s\cc \to T^+s\ci,
\label{eq-ibar-ibar-mu-mu-i}
\end{equation}
which has to be proven for all $p\ge1$. The right hand side is the sum
of terms $i_{m_1}\tens\dots\tens i_{m_t}$ such that
$m_1+\dots+m_t=c_1+\dots+c_p$. Consider a set of positive integers
 \[ N = \{m_1,m_1+m_2,\dots,m_1+\dots+m_t\}
\cap \{c_1,c_1+c_2,\dots,c_1+\dots+c_p\}. \]
It contains $c_1+\dots+c_p$. Clearly, $i_{m_1}\tens\dots\tens i_{m_t}$
will appear in the term $\iu_{n_1}\tens\dots\tens\iu_{n_a}$ if and only
if $N=\{n_1,n_1+n_2,\dots,n_1+\dots+n_a\}$. Since any finite subset
$N\subset\ZZ_{>0}$ has a unique presentation of this form via $n_1$,
\dots, $n_a$, equation~\eqref{eq-ibar-ibar-mu-mu-i} holds.

Let $X$, $Z_j$, $Y$ be objects of $\cc$ and let $C_j^i$ be objects of
$\cb$. When $\iu_n$ is applied to the $\kk$\n-module
\begin{multline}
s\cc(X,C_1^1)\tens\dots\tens s\cc(C^1_{k_1-1},Z_1)\bigotimes
s\cc(Z_1,C^2_1)\tens\dots\tens s\cc(C^2_{k_2-1},Z_2)\bigotimes\dots \\
\bigotimes s\cc(Z_{n-2},C_1^{n-1})\tens\dots\tens
s\cc(C^{n-1}_{k_{n-1}-1},Z_{n-1})\bigotimes
s\cc(Z_{n-1},C^n_1)\tens\dots\tens s\cc(C^n_{k_n-1},Y),
\label{eq-XCCZ-ZCCZ-ZCCY}
\end{multline}
the target space for the term $i_{m_1}\tens\dots\tens i_{m_t}$ has the
form
\[ s\ci(Xi,C^\bull_\bull i)\tens\dots\tens
s\ci(C^\bull_\bull i,C^\bull_\bull i)
\tens\dots\tens s\ci(C^\bull_\bull i,Yi), \]
where $C^\bull_\bull$ are objects of $\cb$ (no $Z_j$ will appear!).
Since $Xi,Yi\in\Ob\ci$ and $C^\bull_\bull i\in\Ob\cj$, the above space
is a direct summand of $s\Dr(\ci|\cj)(Xi,Yi)$. Therefore, the required
map $\iu_n:T^ns\Dr(\cc|\cb)\to s\Dr(\ci|\cj)$ is constructed.

The last statement is a particular case of \eqref{eq-In-mun-iii}.
Indeed, if $k_1=\dots=k_n=1$, then $L(1,\dots,1)$ consists of only one
sequence $(n)$ of the length $t=1$.

Since $\underline{i}$ is strict and $\underline{i}_1=i$,
equation~\eqref{eq-In-mu-i-mu} is the expansion of the definition
$\iu=\bm\underline{i}\bm^{-1}$.
\end{proof}

For example,
\begin{align*}
\iu_1 &= i, \\
\iu_2 &= \mu i - (i\tens i)\mu, \\
\iu_3 &= \mu^{(3)}i - (i\tens\mu i)\mu - (\mu i\tens i)\mu
+ (i\tens i\tens i)\mu^{(3)}, \\
\iu_4 &= \mu^{(4)}i - (i\tens\mu^{(3)}i)\mu
- (\mu i\tens\mu i)\mu - (\mu^{(3)}i\tens i)\mu \\
&\quad+ (i\tens i\tens\mu i)\mu^{(3)} + (i\tens\mu i\tens i)\mu^{(3)}
+ (\mu i\tens i\tens i)\mu^{(3)} - (i\tens i\tens i\tens i)\mu^{(4)}.
\end{align*}

\begin{corollary}
We have a commutative diagram of \ainf-functors
\begin{diagram}
\cb & \rMono & \cc & \rTTo^{\ju^\cc} & \Dr(\cc|\cb) \\
\dTTo<i && \dTTo>i && \dTTo>{\iu} \\
\cj & \rMono & \ci & \rTTo^{\ju^\ci} & \Dr(\ci|\cj)
\end{diagram}
\end{corollary}

When $\ca \rTTo^f \cb \rTTo^g \cc$ are \ainf-functors, then
$\underline{f}\,\underline{g}=\underline{fg}:
\underline{\ca}\to\underline{\cc}$.
This implies
$\overline{f}\,\overline{g}=\overline{fg}:
\overline{\ca}\to\overline{\cc}$.
Assume that $\ca'\hookrightarrow\ca$, $\cb'\hookrightarrow\cb$,
$\cc'\hookrightarrow\cc$ are full \ainf-subcategories such that
$(\Ob\ca')f\subset\Ob\cb'$, $(\Ob\cb')g\subset\Ob\cc'$. Denote
$f'=f\big|_{\ca'}$, $g'=g\big|_{\cb'}$. Since $\Dr(f)$ and
$\Dr(g)$ are just the restrictions of $\overline{f}$ and
$\overline{g}$, we conclude that
\begin{equation}
\Dr(f)\Dr(g)=\Dr(fg): \Dr(\ca|\ca') \to \Dr(\cc|\cc').
\label{eq-Dr-Dr-Dr}
\end{equation}

\subsection{Strict unitality.}\label{sec-iso-strict-uni}
Assume that the \ainf-category $\cc$ is strictly unital. As we know
from \secref{sec-str-unit}, $\Dr(\cc|\cb)$ and $\overline{\cc}$ are
strictly unital with the unit transformation $\uni^{\overline{\cc}}$.
Since $\bm^{-1}:\underline{\cc}\to\overline{\cc}$ is an invertible
\ainf-functor, $\underline{\cc}$ is unital (see
\cite[Section~8.12]{Lyu-AinfCat}). Notice that $\underline{\cc}$ is
never strictly unital except when $\cc=0$, because
\(\underline{b}_2=0\). The transformation
 $\uni^{\underline{\cc}}=\bm^{-1}\uni^{\overline{\cc}}\bm:
 \id_{\underline{\cc}}\to\id_{\underline{\cc}}:
 \underline{\cc}\to\underline{\cc}$,
whose components are
\begin{align*}
\uni^{\underline{\cc}}_0 &= \uni^{\overline{\cc}}_0, \\
\uni^{\underline{\cc}}_1
&= (\uni^{\overline{\cc}}_0\tens1 + 1\tens\uni^{\overline{\cc}}_0)\mu, \\
\uni^{\underline{\cc}}_2
&= (1\tens\uni^{\overline{\cc}}_0\tens1)\mu^{(3)}, \\
\uni^{\underline{\cc}}_k &= 0 \qquad\text{for } k>2,
\end{align*}
is a unit transformation of $\underline{\cc}$. Indeed, let us define
$\uni^{\underline{\cc}}$ by the above components and let us prove that
$\bm\uni^{\underline{\cc}}=\uni^{\overline{\cc}}\bm$. Clearly,
 $(\bm\uni^{\underline{\cc}})_0=\uni^{\overline{\cc}}_0
 =(\uni^{\overline{\cc}}\bm)_0$.
For $n>0$ we have
\begin{align*}
(\bm\uni^{\underline{\cc}})_n &= \mu^{(n)}\uni^{\underline{\cc}}_1
+ \sum_{k+l=n}(\mu^{(k)}\tens\mu^{(l)})\uni^{\underline{\cc}}_2 \\
&= (\uni^{\overline{\cc}}_0\tens1^{\tens n})\mu^{(n+1)}
+ (1^{\tens n}\tens\uni^{\overline{\cc}}_0)\mu^{(n+1)}
+ \sum_{k,l>0;k+l=n}
(1^{\tens k}\tens\uni^{\overline{\cc}}_0\tens1^{\tens l})\mu^{(n+1)} \\
&= \sum_{k,l\ge0;k+l=n}
(1^{\tens k}\tens\uni^{\overline{\cc}}_0\tens1^{\tens l})\mu^{(n+1)}
= (\uni^{\overline{\cc}}\bm)_n.
\end{align*}

\section{\texorpdfstring{An $A_\infty$-transformation}{An A8-transformation}}
Let $\cb\hookrightarrow\cc$ and $\cj\hookrightarrow\ci$ be full
\ainf-subcategories. Let $f,g:\cc\to\ci$ be two \ainf-functors such
that $(\Ob\cb)f\subset\Ob\cj$, $(\Ob\cb)g\subset\Ob\cj$, and let
$r:f\to g:\cc\to\ci$ be an \ainf-transformation. Denote by
$r':f'\to g':\cb\to\cj$ the restriction of $r$ to $\cb$. We already
have $\underline{\cc}$ and $\overline{\cc}$ for \ainf-categories $\cc$,
$\underline{f}$ and $\overline{f}$ for \ainf-functors $f$. Now let us
proceed with \ainf-transformations.

Let us define an \ainf-transformation
 $\underline{r}:\underline{f}\to\underline{g}:
 \underline{\cc}\to\underline{\ci}$
via its components
\begin{alignat}2
\underline{r}_0 &= r_0\underline{j}_1, \qquad &&
\underline{r}_0 = \bigl[\kk \rTTo^{r_0} (s\ci)(Xf,Xg)
\rMono^{\underline{j}_1} (s\underline{\ci})(Xf,Xg) \bigr]; \notag \\
\underline{r}_1 &= r, \qquad && \underline{r}_1 = r\big|_{T^+s\cc}:
T^+s\cc=s\underline{\cc} \to T^+s\ci=s\underline{\ci};
\label{eq-ur0-ur1-urk} \\
\underline{r}_k &= 0 \qquad && \text{for } k>1. \notag
\end{alignat}
Let us check that $\underline{\text{ - }}$ maps the $\omega$\n-globular
set $A_\omega$ \cite[Definition~6.4]{Lyu-AinfCat} into itself (so that
sources and targets are preserved). It suffices to notice that the
correspondence $r\mapsto\underline{r}$ is additive, and if $r=[v,b]$,
then $\underline{r}=[\underline{v},\underline{b}]$. Indeed,
\begin{align}
[\underline{v},\underline{b}]_0 &= \underline{v}_0\underline{b}_1 =
v_0\underline{j}_1b = v_0b_1\underline{j}_1 = r_0\underline{j}_1 =
\underline{r}_0, \notag \\
[\underline{v},\underline{b}]_1 &= \underline{v}_1\underline{b}_1 -(-)^v
\underline{b}_1\underline{v}_1 = vb -(-)^vbv = r = \underline{r}_1,
\label{eq-v-b-01k} \\
[\underline{v},\underline{b}]_k &= \underline{v}_k\underline{b}_1 -(-)^v
\underline{b}_k\underline{v}_1 = 0\pm0 = 0 = \underline{r}_k
\qquad\text{for } k>1. \notag
\end{align}
In particular, a natural \ainf-transformation $r:f\to g:\cc\to\ci$ goes
to the natural \ainf-transformation
$\underline{r}:\underline{f}\to\underline{g}:
\underline{\cc}\to\underline{\ci}$,
and equivalent natural \ainf-transformations $r$, $p$ go to equivalent
$\underline{r}$, $\underline{p}$.

We claim that
\[ r\underline{j}^\ci = \underline{j}^\cc\underline{r}:
f\underline{j} = \underline{j}\,\,\underline{f} \to
g\underline{j} = \underline{j}\,\,\underline{g}:\cc\to \underline{\ci}.
\]
Indeed,
$(r\underline{j})_0-(\underline{j}\underline{r})_0
=r_0\underline{j}_1-\underline{r}_0=0$,
and for $n>0$
\begin{multline*}
(r\underline{j})_n-(\underline{j}\underline{r})_n
= \sum_{a_1+\dots+a_l+k+c_1+\dots+c_m=n}
(f_{a_1}\tens\dots\tens f_{a_l}\tens r_k\tens
g_{c_1}\tens\dots\tens g_{c_m})\underline{j}_{l+1+m}
- \underline{j}_n\underline{r}_1 \\
= r\big|_{T^ns\cc} - r\big|_{T^ns\cc} = 0.
\end{multline*}

We define also the \ainf-transformation conjugate to $\underline{r}$
\[ \overline{r} = \bm\underline{r}\bm^{-1}:
\overline{f} = \bm\underline{f}\bm^{-1} \to
\overline{g} = \bm\underline{g}\bm^{-1}:
\overline{\cc} \to \overline{\ci} \]
(not necessarily natural). Summing up, we have a commutative cylinder
\begin{diagram}[width=4em,LaTeXeqno]
\cb & \rMono & \cc & \rTTo^{\underline{j}^\cc} & \underline{\cc}
& \rTTo^{\bm^{-1}} & \overline{\cc} \\
\dTTo<{f'} \overset{r'}\implies \dTTo>{g'}
&& \dTTo<f \overset{r}\implies \dTTo>g &&
\dTTo<{\underline{f}} \overset{\underline{r}}\implies \dTTo>{\underline{g}}
&&
\dTTo<{\overline{f}} \overset{\overline{r}}\implies \dTTo>{\overline{g}}
\label{dia-underj-bm-1} \\
\cj & \rMono & \ci & \rTTo^{\underline{j}^\ci} & \underline{\ci}
& \rTTo^{\bm^{-1}} & \overline{\ci}
\end{diagram}
The correspondence $\overline{\text{ - }}$ also maps the
$\omega$\n-globular set $A_\omega$ into itself. Indeed, if $r=[v,b]$,
then
\[ \overline{r} = \bm\underline{r}\bm^{-1}
= \bm[\underline{v},\underline{b}]\bm^{-1}
= [\bm\underline{v}\bm^{-1},\bm\underline{b}\bm^{-1}]
= [\overline{v},\overline{b}]. \]

\begin{proposition}\label{prop-over-r-compon}
The \ainf-transformation $\overline{r}$ has the following components
\begin{multline}
\overline{r}_n = \sum^{0\le q\le t}_{l_1+\dots+l_t=n}
(-)^t(\mu^{(l_1)}f\tens\dots\tens\mu^{(l_q)}f\tens r_0\underline{j}_1
\tens\mu^{(l_{q+1})}g\tens\dots\tens\mu^{(l_t)}g)\mu^{(t+1)} \\
+ \sum^{1\le q\le t}_{l_1+\dots+l_t=n}
(-)^{t-1}(\mu^{(l_1)}f\tens\dots\tens\mu^{(l_{q-1})}f\tens\mu^{(l_q)}r
\tens\mu^{(l_{q+1})}g\tens\dots\tens\mu^{(l_t)}g)\mu^{(t)}.
\label{eq-over-rn-r0-r-mu}
\end{multline}
Explicitly, $\overline{r}_0=\underline{r}_0=r_0\underline{j}_1$ and for
$n>0$ the restriction of $\overline{r}_n$ to
$T^{k_1}s\cc\tens\dots\tens T^{k_n}s\cc$ is
\begin{multline}
\overline{r}_n = \mu^{(n)}
\sum_{(a_1,\dots,a_\alpha;k;c_1,\dots,c_\beta)\in P(k_1,\dots,k_n)}
(f_{a_1}\tens\dots\tens f_{a_\alpha}\tens r_k\tens
g_{c_1}\tens\dots\tens g_{c_\beta})\underline{j}_{\alpha+1+\beta}: \\
T^{k_1}s\cc\tens\dots\tens T^{k_n}s\cc \to T^+s\ci,
\label{eq-rn-mun-ffrgg}
\end{multline}
\begin{multline*}
P(k_1,\dots,k_n) = \sqcup_{\alpha,\beta\ge0} \bigl\{
(l_1,\dots,l_\alpha;l_{\alpha+1};l_{\alpha+2},\dots,l_{\alpha+1+\beta})
\in\ZZ^\alpha_{>0}\times\ZZ_{\ge0}\times\ZZ^\beta_{>0} \mid \\
\forall q\in\ZZ_{>0},q\le\alpha+1+\beta\; \forall s\in\ZZ_{\ge0},s\le n
\;\, l_1+\dots+l_q=k_1+\dots+k_s \Leftrightarrow
q=\alpha+1+\beta,s=n \bigr\}.
\end{multline*}
These maps restrict to maps
$\overline{r}_n:T^ns\Dr(\cc|\cb)\to s\Dr(\ci|\cj)$, which are
components of an \ainf-transformation
$\Dr(r)=\overline{r}:\overline{f}\to\overline{g}:
\Dr(\cc|\cb)\to\Dr(\ci|\cj)$.
The restriction of $\overline{r}_n$ to
$T^ns\cc \rMono^{\ju_1^{\tens n}} T^ns\Dr(\cc|\cb)$ equals
$T^ns\cc \rTTo^{r_n} s\ci \rMono^{\ju_1} s\Dr(\ci|\cj)$.
\end{proposition}

\begin{proof}
Similarly to the case of \ainf-functors, discussed in
\propref{prop-In-ce-cf}, let us define an \ainf-transformation
$\overline{r}:\overline{f}\to\overline{g}:\overline{\cc}\to\overline{\ci}$
by its components~\eqref{eq-rn-mun-ffrgg} and prove that the equation
$\overline{r}\bm=\bm\underline{r}$ holds. Clearly,
$(\overline{r}\bm)_0=\overline{r}_0=\underline{r}_0=(\bm\underline{r})_0$.
We have to prove that for $n>0$
\begin{equation}
\sum_{i_1+\dots+i_t=n}
(\overline{f}_{i_1}\tens\dots\tens\overline{f}_{i_{q-1}}
\tens\overline{r}_{i_q}\tens\overline{g}_{i_{q+1}}
\tens\dots\tens\overline{g}_{i_t})\mu^{(t)} = \mu^{(n)}r:
T^{k_1}s\cc\tens\dots\tens T^{k_n}s\cc \to T^+s\ci.
\label{eq-of-of-or-og-og}
\end{equation}
The right hand side is the sum of terms
$f_{a_1}\tens\dots\tens f_{a_\alpha}\tens r_k\tens
g_{c_1}\tens\dots\tens g_{c_\beta}$,
such that $a_1+\dots+a_\alpha+k+c_1+\dots+c_\beta=k_1+\dots+k_n$.
Denote
$(l_1,\dots,l_{\alpha+1+\beta})=(a_1,\dots,a_\alpha,k,c_1,\dots,c_\beta)$.
Consider the subsequence $N$ of the sequence
$L=(0,l_1,l_1+l_2,\dots,l_1+\dots+l_{\alpha+1+\beta})$ consisting of
all elements which belong to the set
$\{0,k_1,k_1+k_2,\dots,k_1+\dots+k_n\}$. The term
$f_{a_1}\tens\dots\tens f_{a_\alpha}\tens r_k\tens
g_{c_1}\tens\dots\tens g_{c_\beta}$
will appear as a summand of the term
$\overline{f}_{i_1}\tens\dots\tens\overline{f}_{i_{q-1}}
\tens\overline{r}_{i_q}\tens\overline{g}_{i_{q+1}}
\tens\dots\tens\overline{g}_{i_t}$
if and only if
\begin{gather}
N=(0,i_1,i_1+i_2,\dots,i_1+\dots+i_t), \label{eq-N-0lllll} \\
i_1+\dots+i_{q-1}\le a_1+\dots+a_\alpha, \qquad
a_1+\dots+a_\alpha+k\le i_1+\dots+i_q, \\
\forall 1\le y\le t \qquad i_y=0 \implies y=q.
\label{eq-li0-qi}
\end{gather}

Let us prove that for a given sequence
$(a_1,\dots,a_\alpha;k;c_1,\dots,c_\beta)
\in\ZZ^\alpha_{>0}\times\ZZ_{\ge0}\times\ZZ^\beta_{>0}$
there exists exactly one sequence
$(i_1,\dots,i_t;q)\in\ZZ^t_{\ge0}\times\ZZ_{>0}$ such that conditions
\eqref{eq-N-0lllll}--\eqref{eq-li0-qi} are satisfied. Indeed, the
sequence $N$ determines uniquely a sequence $(i_1,\dots,i_t)$ of
non-negative integers such that \eqref{eq-N-0lllll} holds. If $k>0$ or
$a_1+\dots+a_\alpha$ does not belong to $N$, then all $i_y$ are
positive. This implies that the interval
$[a_1+\dots+a_\alpha,a_1+\dots+a_\alpha+k]$ is contained in a unique
interval of the form $[i_1+\dots+i_{q-1},i_1+\dots+i_q]$.

If $k=0$ and $a_1+\dots+a_\alpha$ belongs to $N$, then
$a_1+\dots+a_\alpha=a_1+\dots+a_\alpha+k$ is repeated in $N$. Hence,
there exists $q>0$ such that $a_1+\dots+a_\alpha=i_1+\dots+i_{q-1}$,
$i_q=0$. Since $M$ and $N$ may contain no more than one repeated
element, $i_y$ is positive for $y\ne q$. Therefore, conditions
\eqref{eq-N-0lllll}--\eqref{eq-li0-qi} are satisfied.

We conclude that \eqref{eq-of-of-or-og-og} holds.
Formula~\eqref{eq-over-rn-r0-r-mu} is the expansion of the proven
property $\overline{r} = \bm\underline{r}\bm^{-1}$.

The target space for map~\eqref{eq-rn-mun-ffrgg} applied to
$\kk$\n-module~\eqref{eq-XCCZ-ZCCZ-ZCCY} has the form
\[ s\ci(Xf,C^\bull_\bull f)\tens\dots\tens
s\ci(C^\bull_\bull f,C^\bull_\bull g)
\tens\dots\tens s\ci(C^\bull_\bull g,Yg), \]
where $C^\bull_\bull$ are objects of $\cb$ (and no $Z_j$ will appear).
Since objects $C^\bull_\bull f$ and $C^\bull_\bull g$ belong to $\cj$,
the above space is a direct summand of $s\Dr(\ci|\cj)(Xf,Yg)$.
Therefore, the required map
$\overline{r}_n:T^ns\Dr(\cc|\cb)\to s\Dr(\ci|\cj)$ is constructed.

The last statement is a particular case of \eqref{eq-rn-mun-ffrgg}.
Indeed, if $k_1=\dots=k_n=1$, then $P(1,\dots,1)$ consists of only one
element $(;n;)$, that is, $\alpha=\beta=0$, $i_1=n\in\ZZ_{\ge0}$.
\end{proof}

In particular, the correspondence $r\mapsto\Dr(r)=\overline{r}$ maps
natural \ainf-transformations to natural ones, and equivalent
$r,p:f\to g:\cb\to\cc$ are mapped to equivalent
\[ \Dr(r), \Dr(p): \overline{f} \to \overline{g}:
\Dr(\cc|\cb) \to \Dr(\ci|\cj). \]

For example,
\begin{align*}
\overline{r}_1 &= r - (f\tens r_0 + r_0\tens g)\mu, \\
\overline{r}_2 &= \mu r - (f\tens r + r\tens g)\mu
- \mu(f\tens r_0 + r_0\tens g)\mu \\
&\qquad+ (f\tens f\tens r_0 + f\tens r_0\tens g + r_0\tens g\tens g)\mu^{(3)}.
\end{align*}

\begin{corollary}\label{cor-ju-cylinder}
We have a commutative cylinder
\begin{diagram}[width=4em]
\cb & \rMono & \cc & \rTTo^{\ju^\cc} & \Dr(\cc|\cb) \\
\dTTo<{f'} \overset{r'}\implies \dTTo>{g'}
&& \dTTo<f \overset{r}\implies \dTTo>g &&
\dTTo<{\overline{f}} \overset{\overline{r}}\implies \dTTo>{\overline{g}} \\
\cj & \rMono & \ci & \rTTo^{\ju^\ci} & \Dr(\ci|\cj)
\end{diagram}
\end{corollary}

\subsection{\texorpdfstring{$\ck$-2-categories and $\ck$-2-functors.}
    {K-2-categories and K-2-functors.}}
Let $\ck$ denote the category $\Kht(\kk\modul)=H^0(\Com(\kk\modul))$ of
differential graded complexes of $\kk$\n-modules, whose morphisms are
chain maps modulo homotopy. A 1\n-unital, non-2-unital $\ck$-2-category
$\ck A_\infty$ of \ainf-categories is described in
\cite[Proposition~7.1]{Lyu-AinfCat}. Instead of the complex of
2\n-morphisms $(A_\infty(\ca,\cb)(f,g),m_1)$, $m_1=sB_1s^{-1}$, we work
with the shifted complex $(sA_\infty(\ca,\cb)(f,g),B_1)$. There is an
obvious notion of a strict $\ck$-2-functor between such
$\ck$-2-categories -- a map of objects, maps of 1\n-morphisms and chain
maps of 2\n-morphisms, which preserve all operations. The operations
involving 2\n-morphisms are subject to equations in $\ck$, which mean
equations between chain maps up to homotopy.

We have applied the ``underline'' construction $\unmi$ to three kinds
of arguments $\text-$ : \ainf-categories, \ainf-functors and
\ainf-transformations. Let us summarize the properties of this
construction.

\begin{proposition}\label{pro-K2fun-unmi}
The following assignment defines a strict $\ck$-2-functor
$\unmi:\ck A_\infty\to\ck A_\infty$ : an \ainf-category $\ca$ is mapped
to $\underline{\ca}$, an \ainf-functor $f:\ca\to\cb$ is mapped to
$\underline{f}:\underline{\ca}\to\underline{\cb}$, and the chain map of
complexes of 2\n-morphisms is
\[ \unmi: (sA_\infty(\ca,\cb)(f,g),B_1) \to
(sA_\infty(\underline{\ca},\underline{\cb})
(\underline{f},\underline{g}),\underline{B}_1),
\qquad r\mapsto\underline{r},
\]
where $\underline{B}$ denotes the codifferential in
$TsA_\infty(\underline{\ca},\underline{\cb})$, in particular,
$v\underline{B}_1=[v,\underline{b}]$.
\end{proposition}

\begin{proof}
We have seen in \eqref{eq-v-b-01k} that
$\underline{rB_1}=\underline{[r,b]}=[\underline{r},\underline{b}]
=\underline{r}\,\,\underline{B}_1$,
thus, $r\mapsto\underline{r}$ is a chain map. The composition of
\ainf-functors is preserved,
$\underline{fg}=\underline{f}\,\,\underline{g}$. The right action of a
1\n-morphism $h$ on a 2\n-morphism $r$ is preserved, since
$\underline{rh}=(r_0h_1\underline{j}_1,rh,0,0,\dots)
=\underline{r}\,\,\underline{h}$.
The left action of a 1\n-morphism $e$ on a 2\n-morphism $r$ is
preserved, since
$\underline{er}=(\sS{_{\_e}}r_0\underline{j}_1,er,0,0,\dots)
=\underline{e}\,\,\underline{r}$.
The identity \ainf-functor $\id_\ca$ is mapped to the identity
\ainf-functor $\underline{\id_\ca}=\id_{\underline{\ca}}$.

It remains to prove that the vertical composition of 2\n-morphisms
\[ m_2 = (s\tens s)B_2s^{-1}:
A_\infty(\ca,\cb)(f,g)\tens A_\infty(\ca,\cb)(g,h) \to
A_\infty(\ca,\cb)(f,h)
\]
is preserved, that is, the diagram
\begin{diagram}
A_\infty(\ca,\cb)(f,g)\tens A_\infty(\ca,\cb)(g,h)
& \rTTo^{s\unmi s^{-1}\tens s\unmi s^{-1}} &
A_\infty(\underline{\ca},\underline{\cb})(\underline{f},\underline{g})\tens
A_\infty(\underline{\ca},\underline{\cb})(\underline{g},\underline{h}) \\
\dTTo<{m_2} && \dTTo>{\underline{m}_2} \\
A_\infty(\ca,\cb)(f,h) & \rTTo^{s\unmi s^{-1}} &
A_\infty(\underline{\ca},\underline{\cb})(\underline{f},\underline{h})
\end{diagram}
is commutative in $\ck$. Here \(s\unmi s^{-1}\) denotes the composition
\[ A_\infty(\ca,\cb)(f,h) \rTTo^s sA_\infty(\ca,\cb)(f,h)
\rTTo^\unmi
sA_\infty(\underline{\ca},\underline{\cb})(\underline{f},\underline{h})
\rTTo^{s^{-1}}
A_\infty(\underline{\ca},\underline{\cb})(\underline{f},\underline{h}).
\]
Since $\underline{b}_k=0$ for $k\ge2$, we have $\underline{B}_2=0$ due
to \cite[Equation~(5.1.3)]{Lyu-AinfCat}, hence, $\underline{m}_2=0$.
Let us prove that $m_2(s\unmi s^{-1})\sim0$. The homotopy is sought in
the form $(s\tens s)Hs^{-1}$, where
\[ H: sA_\infty(\ca,\cb)(f,g)\tens sA_\infty(\ca,\cb)(g,h) \to
sA_\infty(\underline{\ca},\underline{\cb})(\underline{f},\underline{h})
\]
is a $\kk$\n-linear map of degree 0. It has to satisfy the equation
\[ B_2\unmi = H\underline{B}_1 - (1\tens B_1 + B_1\tens1)H, \]
that is, for each $r\in sA_\infty(\ca,\cb)(f,g)$,
$p\in sA_\infty(\ca,\cb)(g,h)$
\begin{equation}
\underline{(r\tens p)B_2} = [(r\tens p)H,\underline{b}]
- [(r\tens p)(1\tens B_1+B_1\tens1)]H.
\label{eq-un-r-p-B2}
\end{equation}
A candidate for $H$ is chosen similarly to
definition~\eqref{eq-ur0-ur1-urk}. We choose the components of the
$(\underline{f},\underline{h})$-coderivation $(r\tens p)H$ as follows:
\begin{align*}
[(r\tens p)H]_0 & = (r_0\tens p_0)\underline{j}_2, \\
[(r\tens p)H]_1 & = (r\tens p)\theta: s\underline{\ca}(X,Y) \to
s\underline{\cb}(Xf,Yh), \\
[(r\tens p)H]_k & = 0 \qquad\qquad \text{ for } k>1.
\end{align*}
Let us verify equation~\eqref{eq-un-r-p-B2} for this $H$. Both
sides of \eqref{eq-un-r-p-B2} are
$(\underline{f},\underline{h})$-coderivations. It suffices to check
that all their components coincide. The 0\n-th component of the right
hand side of \eqref{eq-un-r-p-B2} is
\begin{align*}
& [(r\tens p)H]_0\underline{b}_1-[(r\tens[p,b]+(-)^p[r,b]\tens p)H]_0 \\
&= (r_0\tens p_0)\underline{j}_2b -
(r_0\tens p_0b_1 + (-)^pr_0b_1\tens p_0)\underline{j}_2 \\
&= (r_0\tens p_0)(b - 1\tens b_1 - b_1\tens1)
= (r_0\tens p_0)b_2\underline{j}_2,
\end{align*}
which equals $[\underline{(r\tens p)B_2}]_0$. Due to
\cite[Equation~(5.1.2)]{Lyu-AinfCat} the first component of the right
hand side of \eqref{eq-un-r-p-B2} equals
\begin{align*}
& [(r\tens p)H]_1\underline{b}_1
- (-)^{r+p}\underline{b}_1[(r\tens p)H]_1
- [(r\tens p)(1\tens B_1+B_1\tens1)H]_1 \\
&= (r\tens p)\theta b - (-)^{r+p}b(r\tens p)\theta
- [(r\tens p)(1\tens B_1+B_1\tens1)]\theta \\
&= (r\tens p)B_2,
\end{align*}
which is $[\underline{(r\tens p)B_2}]_1$. The $k$\n-th component of the
right hand side of \eqref{eq-un-r-p-B2} vanishes for $k>1$, and so does
$[\underline{(r\tens p)B_2}]_k=0$. Therefore, \eqref{eq-un-r-p-B2} and
the proposition are proven.
\end{proof}

\begin{corollary}\label{cor-2fun-unmi-A-A}
The same assignment $\ca\mapsto\underline{\ca}$,
$f\mapsto\underline{f}$, $r\mapsto\underline{r}$ as in
\propref{pro-K2fun-unmi} gives a strict 2\n-functor
$\unmi:A_\infty\to A_\infty$ of non-2-unital 2-categories.
\end{corollary}

This is obtained by taking the 0\n-th cohomology of $\ck A_\infty$ in
\propref{pro-K2fun-unmi}.

Similarly, the ``overline'' construction $\ovmi$, applied to three
kinds of arguments, \ainf-categories, \ainf-functors and
\ainf-transformations, gives a strict $\ck$-2-functor.

\begin{corollary}\label{cor-K2fun-ovmi}
The following assignment defines a strict $\ck$-2-functor
$\ovmi:\ck A_\infty\to\ck A_\infty$: an \ainf-category $\ca$ is mapped
to $\overline{\ca}$, an \ainf-functor $f:\ca\to\cb$ is mapped to
$\overline{f}:\overline{\ca}\to\overline{\cb}$, and the chain map of
complexes of 2\n-morphisms is
\[ \ovmi: (sA_\infty(\ca,\cb)(f,g),B_1) \to
(sA_\infty(\overline{\ca},\overline{\cb})
(\overline{f},\overline{g}),\overline{B}_1),
\qquad r\mapsto\overline{r}, \]
where $\overline{B}$ denotes the codifferential in
$TsA_\infty(\overline{\ca},\overline{\cb})$, in particular,
$v\overline{B}_1=[v,\overline{b}]$. There is an invertible strict
$\ck$\n-2-transformation $\bm:\ovmi\to\unmi$,
$\bm_\ca:\overline{\ca}\to\underline{\ca}$.
\end{corollary}

\begin{proof}
Starting with a $\ck$-2-functor $\unmi$, a mapping
$\Ob\ck A_\infty\to\Ob\ck A_\infty$, $\ca\mapsto\overline{\ca}$, and a
family of invertible \ainf-functors
$\bm_\ca:\overline{\ca}\to\underline{\ca}$, one may construct another
$\ck$\n-2-functor $\ovmi$, which maps an \ainf-category $\ca$ to
$\overline{\ca}$, so that $\bm$ is a strict $\ck$\n-2-transformation.
Since $\bm$ is strict and invertible, the values of $\overline{f}$ and
$\overline{r}$ are fixed by the requirements
$\overline{f}\bm=\bm\underline{f}$, $\overline{r}\bm=\bm\underline{r}$
for each \ainf-functor $f$ and \ainf-transformation $r$.
\end{proof}

The detailed definition of strict $\ck$\n-2-transformations is left to
the interested reader.

\begin{corollary}\label{cor-2fun-ovmi-A-A}
The same assignment $\ca\mapsto\overline{\ca}$, $f\mapsto\overline{f}$,
$r\mapsto\overline{r}$ as in \corref{cor-K2fun-ovmi} gives a strict
2\n-functor $\ovmi:A_\infty\to A_\infty$ of non-2-unital categories.
\end{corollary}

This is obtained by taking the 0\n-th cohomology of $\ck A_\infty$ in
\corref{cor-K2fun-ovmi}.

It is instructive to find the homotopy which forces $\ovmi$ to preserve
the vertical composition of 2\n-morphisms. Denote $\ad\bm$ the maps
$sA_\infty(\underline{\ca},\underline{\cb})(\underline{f},\underline{g})
\to sA_\infty(\overline{\ca},\overline{\cb})(\overline{f},\overline{g})$,
$v\mapsto\bm v\bm^{-1}$. The following diagram commutes modulo homotopy:
\begin{diagram}[LaTeXeqno]
\begin{array}{l}
sA_\infty(\ca,\cb)(f,g)\tens
\\
\quad\tens sA_\infty(\ca,\cb)(g,h)\quad\;
\end{array}
&& \rTTo^{\ovmi\tens\ovmi} &&
\begin{array}{l}
\; sA_\infty(\overline{\ca},\overline{\cb})(\overline{f},\overline{g})\tens
\\
\;\quad\tens sA_\infty(\overline{\ca},\overline{\cb})(\overline{g},\overline{h})\quad
\end{array}
\\
& \rdTTo_{\unmi\tens\unmi} & = & \ruTTo_{\ad\bm\tens\ad\bm} &
\\
&&
\begin{array}{l}
\;\quad sA_\infty(\underline{\ca},\underline{\cb})(\underline{f},\underline{g})\tens
\\
\;\quad\quad\tens sA_\infty(\underline{\ca},\underline{\cb})(\underline{g},\underline{h})\;
\end{array}
&&
\\
\dTTo<{B_2} & \sim & \dTTo<{0=}>{\underline{B}_2} & \sim
& \dTTo>{\overline{B}_2} \\
&&
sA_\infty(\underline{\ca},\underline{\cb})(\underline{f},\underline{h})
&& \\
& \ruTTo^{\unmi} & = & \rdTTo^{\ad\bm} & \\
sA_\infty(\ca,\cb)(f,h) && \rTTo_{\ovmi} &&
sA_\infty(\overline{\ca},\overline{\cb})(\overline{f},\overline{h})
\label{dia-ovov-ovB2-B2-ov}
\end{diagram}
The right homotopy commutative square is obtained from
\cite[Equation~(7.1.2)]{Lyu-AinfCat}:
\[ (\rho\tens\pi)\underline{B}_2\bm^{-1} -
(\rho\bm^{-1}\tens\pi\bm^{-1})B_2 = (\rho\tens\pi\mid\bm^{-1})M_{20}B_1-
[(\rho\tens\pi)(1\tens\underline{B}_1+\underline{B}_1\tens1)|\bm^{-1}]M_{20}
\]
for all $\rho\in
sA_\infty(\underline{\ca},\underline{\cb})(\underline{f},\underline{g})$,
$\pi\in
sA_\infty(\underline{\ca},\underline{\cb})(\underline{g},\underline{h})$.
Recall that $\underline{B}_2=0$ and compose with $\bm$ to get
\[ - (\bm\rho\bm^{-1}\tens\bm\pi\bm^{-1})\overline{B}_2
= [\bm(\rho\tens\pi\mid\bm^{-1})M_{20}]\overline{B}_1 -
\bm[(\rho\tens\pi)(1\tens\underline{B}_1+\underline{B}_1\tens1)
|\bm^{-1}]M_{20}. \]
In particular, for $\rho=\underline{r}$, $\pi=\underline{p}$ we have
\[ - (\overline{r}\tens\overline{p})\overline{B}_2 =
[\bm(\underline{r}\tens\underline{p}\mid\bm^{-1})M_{20}]\overline{B}_1
- \bm\{[(r\tens p)(1\tens B_1+B_1\tens1)](\unmi\tens\unmi)
|\bm^{-1}\}M_{20}. \]
The left homotopy commutative square, that is, \eqref{eq-un-r-p-B2}
composed with $\ad\bm$ gives
\[ \overline{(r\tens p)B_2} = [\bm(r\tens p)H\bm^{-1}]\overline{B}_1
- \bm[(r\tens p)(1\tens B_1+B_1\tens1)]H\bm^{-1}.
\]
We conclude that the exterior of diagram \eqref{dia-ovov-ovB2-B2-ov} is
commutative modulo homotopy
\begin{gather*}
R: sA_\infty(\ca,\cb)(f,g)\tens sA_\infty(\ca,\cb)(g,h) \to
sA_\infty(\overline{\ca},\overline{\cb})(\overline{f},\overline{h}), \\
(r\tens p)R = \bm(r\tens p)H\bm^{-1}
+ \bm(\underline{r}\tens\underline{p}\mid\bm^{-1})M_{20},
\end{gather*}
that is,
\begin{equation}
\overline{(r\tens p)B_2} -
(\overline{r}\tens\overline{p})\overline{B}_2
= (r\tens p)R\overline{B}_1 - [(r\tens p)(1\tens B_1+B_1\tens1)]R.
\label{eq-ovrpB2-ovrovpovB2}
\end{equation}

\begin{proposition}\label{pro-homo-rpR-explicit}
The $(\overline{f},\overline{h})$-transformation $(r\tens p)R$ has the
following components: $[(r\tens p)R]_0=0$, and for $n>0$ the restriction
of $[(r\tens p)R]_n$ to
$T^{k_1}s\ca\tens\dots\tens T^{k_n}s\ca$ is
\begin{multline}
[(r\tens p)R]_n = \mu^{(n)}
\sum_{(\bar{a};k;\bar{c};t;\bar{e})\in Q(k_1,\dots,k_n)}
(f_{a_1}\tens\dots\tens f_{a_\alpha}\tens r_k\tens
g_{c_1}\tens\dots\tens g_{c_\beta}\tens p_t\tens
h_{e_1}\tens\dots\tens h_{e_\gamma}) \\
\underline{j}_{\alpha+\beta+\gamma+2}:
T^{k_1}s\ca\tens\dots\tens T^{k_n}s\ca \to T^+s\cb,
\label{eq-rpRn-ffrggphh}
\end{multline}
where
$(\bar{a};k;\bar{c};t;\bar{e})
=(a_1,\dots,a_\alpha;k;c_1,\dots,c_\beta;t;e_1,\dots,e_\gamma)$
and
\begin{multline*}
Q(k_1,\dots,k_n) = \sqcup_{\alpha,\beta,\gamma\ge0} \bigl\{
(l_1,\dots,l_\alpha;l_{\alpha+1};l_{\alpha+2},\dots,l_{\alpha+\beta+1};
l_{\alpha+\beta+2};l_{\alpha+\beta+3},\dots,l_{\alpha+\beta+\gamma+2}) \\
\in\ZZ^\alpha_{>0}\times\ZZ_{\ge0}\times\ZZ^\beta_{>0}\times\ZZ_{\ge0}
\times\ZZ^\gamma_{>0} \mid
\forall q\in\ZZ_{>0},q\le\alpha+\beta+\gamma+2\;
\forall s\in\ZZ_{\ge0},s\le n \\
l_1+\dots+l_q=k_1+\dots+k_s \Leftrightarrow
q=\alpha+\beta+\gamma+2,s=n \bigr\}.
\end{multline*}
If $\ca'\subset\ca$, $\cb'\subset\cb$ are full \ainf-subcategories and
$(\Ob\ca')f\subset\Ob\cb'$, $(\Ob\ca')g\subset\Ob\cb'$,
$(\Ob\ca')h\subset\Ob\cb'$, then $[(r\tens p)R]_n$ restrict to maps
\[ [(r\tens p)R]_n: T^ns\Dr(\ca|\ca')(X,Y)\to s\Dr(\cb|\cb')(Xf,Yh), \]
which are components of an \ainf-transformation
\[ (r\tens p)R \in sA_\infty\bigl(\Dr(\ca|\ca'),\Dr(\cb|\cb')\bigr)
(\overline{f},\overline{h}). \]
\end{proposition}

\begin{proof}
Denote by $R'$ an $(\overline{f},\overline{h})$-coderivation, whose
components are $R'_0=0$ and $R'_n$ is given by the right hand side of
\eqref{eq-rpRn-ffrggphh}. We want to prove that $(r\tens p)R=R'$. This
is equivalent to the equation
\begin{equation}
R'\bm = \bm(r\tens p)H
+ \bm(\underline{r}\tens\underline{p}\mid\bm^{-1})M_{20}\bm.
\label{eq-Rm-mrpH-mrpmMm}
\end{equation}
Let us transform the last term. Applying the identity
$(1\boxtimes M)M\pr_1=(M\boxtimes1)M\pr_1$
\cite[Proposition~4.1]{Lyu-AinfCat} to an element
\[ 1\tens\underline{r}\tens\underline{p}\tens1 \in
T^0sA_\infty(\overline{\ca},\underline{\ca})(\bm,\bm)\tens
T^2sA_\infty(\underline{\ca},\underline{\cb})(\underline{f},\underline{h})
\tens T^0sA_\infty(\underline{\cb},\overline{\cb})(\bm^{-1},\bm^{-1}) \]
we find from
\begin{multline*}
(1\tens\underline{r}\tens\underline{p}\tens1)(1\tens M)M\pr_1 \\
= [1\tens(\underline{r}\tens\underline{p}\mid\bm^{-1})M_{20}
+ 1\tens\underline{r}\bm^{-1}\tens\underline{p}\bm^{-1}]M\pr_1
= \bm(\underline{r}\tens\underline{p}\mid\bm^{-1})M_{20},
\end{multline*}
\begin{equation*}
(1\tens\underline{r}\tens\underline{p}\tens1)(M\tens1)M\pr_1
= (\bm\underline{r}\tens\bm\underline{p}\tens1)M\pr_1
= (\bm\underline{r}\tens\bm\underline{p}\mid\bm^{-1})M_{20},
\end{equation*}
that
$\bm(\underline{r}\tens\underline{p}\mid\bm^{-1})M_{20}
=(\bm\underline{r}\tens\bm\underline{p}\mid\bm^{-1})M_{20}$.
Applying the same identity to an element
\[ \rho\tens\pi\tens1\tens1 \in
T^2sA_\infty(\overline{\ca},\underline{\cb})(\bm f,\bm h)\tens
T^0sA_\infty(\underline{\cb},\overline{\cb})(\lambda,\lambda)\tens
T^0sA_\infty(\overline{\cb},\underline{\cb})(\bm,\bm) \]
we find from
\begin{equation*}
(\rho\tens\pi\tens1\tens1)(1\tens M)M\pr_1
= (\rho\tens\pi\tens1)M\pr_1 = (\rho\tens\pi\mid\lambda\bm)M_{20},
\end{equation*}
\begin{multline*}
(\rho\tens\pi\tens1\tens1)(M\tens1)M\pr_1
= [(\rho\tens\pi\mid\lambda)M_{20}\tens1
+ \rho\lambda\tens\pi\lambda\tens1]M\pr_1 \\
= (\rho\tens\pi\mid\lambda)M_{20}\bm
+ (\rho\lambda\tens\pi\lambda\mid\bm)M_{20},
\end{multline*}
that
\[ (\rho\tens\pi\mid\lambda\bm)M_{20}
= (\rho\tens\pi\mid\lambda)M_{20}\bm
+ (\rho\lambda\tens\pi\lambda\mid\bm)M_{20}. \]
When $\lambda=\bm^{-1}$, the left hand side vanishes. Indeed, for each
$k\ge0$
\[ [(\rho\tens\pi\mid\id_{\underline{\cb}})M_{20}]_k
= \sum_{l\ge2} (\rho\tens\pi)\theta_{kl}\id_l = 0. \]
Hence,
\[ (\rho\tens\pi\mid\bm^{-1})M_{20}\bm
= - (\rho\bm^{-1}\tens\pi\bm^{-1}\mid\bm)M_{20}. \]
In particular, for $\rho=\bm\underline{r}$, $\pi=\bm\underline{p}$ we
have
\[ (\bm\underline{r}\tens\bm\underline{p}\mid\bm^{-1})M_{20}\bm
= - (\overline{r}\tens\overline{p}\mid\bm)M_{20}. \]
Therefore, equation~\eqref{eq-Rm-mrpH-mrpmMm}  can
be rewritten as follows:
\begin{equation}
\bm(r\tens p)H = R'\bm + (\overline{r}\tens\overline{p}\mid\bm)M_{20}.
\label{eq-mrpH-Rm-mrpmM}
\end{equation}
Both sides are $(\bm\underline{f},\bm\underline{h})$-coderivations,
or $(\overline{f}\bm,\overline{h}\bm)$-coderivations, which is the same
thing. Let us prove that all their components coincide.

The 0-th components coincide, since
\[ [(r\tens p)H]_0 = (r_0\tens p_0)\underline{j}_2
= (\overline{r}\tens\overline{p})\theta_{02}\bm_2
= [(\overline{r}\tens\overline{p}\mid\bm)M_{20}]_0. \]
For $n>0$ we have to verify the following equation for $n$\n-th
components
\begin{multline*}
\mu^{(n)}(r\tens p)\theta = \sum_{i_1+\dots+i_x=n}
\bigl(\overline{f}_{i_1}\tens\dots\tens\overline{f}_{i_{q-1}}\tens
R'_{i_q}\tens\overline{h}_{i_{q+1}}
\tens\dots\tens\overline{h}_{i_x}\bigr) \mu^{(x)} \\
+ \sum_x (\overline{r}\tens\overline{p})\theta_{nx}\mu^{(x)}:
T^{k_1}s\ca\tens\dots\tens T^{k_n}s\ca \to T^+s\cb.
\end{multline*}
The left hand side is
\begin{equation}
\sum_{a_1+\dots+a_\alpha+k+c_1+\dots+c_\beta+t+e_1+\dots+e_\gamma=n}
f_{a_1}\tens\dots\tens f_{a_\alpha}\tens r_k\tens
g_{c_1}\tens\dots\tens g_{c_\beta}\tens p_t\tens
h_{e_1}\tens\dots\tens h_{e_\gamma}.
\label{eq-LHS-ffrggphh}
\end{equation}
Both sums in the right hand side consist of some of the above
summands. Let us verify that each summand of \eqref{eq-LHS-ffrggphh}
will occur exactly once either in
\( \sum (\overline{f}_{i_1}\tens\dots\tens\overline{f}_{i_{q-1}}\tens
R'_{i_q}\tens\overline{h}_{i_{q+1}}
\tens\dots\tens\overline{h}_{i_x}) \mu^{(x)} \),
or in
\( \sum_x (\overline{r}\tens\overline{p})\theta_{nx}\mu^{(x)} \).
Let us rewrite the sequence
\( (a_1,\dots,a_\alpha;k;c_1,\dots,c_\beta;t;e_1,\dots,e_\gamma) \)
as
\[
(l_1,\dots,l_\alpha;l_{\alpha+1};l_{\alpha+2},\dots,l_{\alpha+\beta+1};
l_{\alpha+\beta+2};l_{\alpha+\beta+3},\dots,l_{\alpha+\beta+\gamma+2}).
\]
Consider the subsequence $N$ of the sequence
\( L=(0,l_1,l_1+l_2,\dots,l_1+\dots+l_{\alpha+\beta+\gamma+2}) \)
consisting of all elements which belong to the set
$\{0,k_1,k_1+k_2,\dots,k_1+\dots+k_n\}$. The term
\begin{equation}
f_{l_1}\tens\dots\tens f_{l_\alpha}\tens r_{l_{\alpha+1}}\tens
g_{l_{\alpha+2}}\tens\dots\tens g_{l_{\alpha+\beta+1}}\tens
p_{l_{\alpha+\beta+2}}\tens h_{l_{\alpha+\beta+3}}\tens\dots\tens
h_{l_{\alpha+\beta+\gamma+2}}
\label{eq-llllllll-ffrggphh}
\end{equation}
will appear as a summand of
\begin{equation}
(\overline{f}_{i_1}\tens\dots\tens\overline{f}_{i_{q-1}}\tens
R'_{i_q}\tens\overline{h}_{i_{q+1}}
\tens\dots\tens\overline{h}_{i_x}) \mu^{(x)},
\label{eq-of-ofRoh-oh-m}
\end{equation}
if and only if
\begin{gather}
N=(0,i_1,i_1+i_2,\dots,i_1+\dots+i_x), \label{eq-N-0iiiii} \\
i_1+\dots+i_{q-1} \le l_1+\dots+l_\alpha, \qquad
l_1+\dots+l_{\alpha+\beta+2} \le i_1+\dots+i_q, \label{eq-iill-llii} \\
\forall 1\le\gamma\le x \qquad i_\gamma=0 \implies \gamma=q.
\label{eq-a1gx-ig0-gq}
\end{gather}
The term \eqref{eq-llllllll-ffrggphh}
will appear as a summand of
\begin{equation}
(\overline{f}_{i_1}\tens\dots\tens\overline{f}_{i_{y-1}}\tens
\overline{r}_{i_y}\tens
\overline{g}_{i_{y+1}}\tens\dots\tens\overline{g}_{i_{z-1}}\tens
\overline{p}_{i_z}\tens
\overline{h}_{i_{z+1}}\tens\dots\tens\overline{h}_{i_x}) \mu^{(x)}
\label{eq-of-of-or-og-og-op-oh-oh-m}
\end{equation}
(which is a term of
$(\overline{r}\tens\overline{p})\theta_{nx}\mu^{(x)}$) if and only if
\eqref{eq-N-0iiiii} holds and
\begin{gather}
i_1+\dots+i_{y-1} \le l_1+\dots+l_\alpha, \qquad
l_1+\dots+l_{\alpha+1} \le i_1+\dots+i_y, \label{eq-iiyll-lliiy} \\
i_1+\dots+i_{z-1} \le l_1+\dots+l_{\alpha+\beta+1}, \qquad
l_1+\dots+l_{\alpha+\beta+2} \le i_1+\dots+i_z, \label{eq-iizll-lliiz} \\
y<z \text{ and } \forall\, 1\le\gamma\le x \qquad
i_\gamma=0 \implies \gamma\in\{y,z\}. \label{eq-yz-a1gx-ig0-gyz}
\end{gather}
A given non-decreasing sequence $N$ determines uniquely a sequence
$(i_1,\dots,i_x)$ of non-negative integers such that
\eqref{eq-N-0iiiii} holds.

If $l_{\alpha+1}>0$ or $l_1+\dots+l_\alpha$ does not belong to $N$,
then there exists exactly one element $y=y'$ such that
\eqref{eq-iiyll-lliiy} holds. If $l_{\alpha+1}=0$ and
$l_1+\dots+l_\alpha$ belongs to $N$, then there are at least 2 such
elements. Denote by $y'>0$ the least of them. Then
$l_1+\dots+l_\alpha=i_1+\dots+i_{y'-1}$ and $i_{y'}=0$. If $y$
satisfies both \eqref{eq-iiyll-lliiy} and \eqref{eq-yz-a1gx-ig0-gyz},
then $y\le y'$, hence, $y=y'$ is the only solution.

If $l_{\alpha+\beta+2}>0$ or $l_1+\dots+l_{\alpha+\beta+1}$ does not
belong to $N$, then there exists exactly one element $z=z'$ such that
\eqref{eq-iizll-lliiz} holds. If $l_{\alpha+\beta+2}=0$ and
$l_1+\dots+l_{\alpha+\beta+1}$ belongs to $N$, then there are at least
2 such elements. Denote by $z'$ the biggest of them. Then
$l_1+\dots+l_{\alpha+\beta+1}=i_1+\dots+i_{z'-1}$ and $i_{z'}=0$. If
$z$ satisfies both \eqref{eq-iizll-lliiz} and
\eqref{eq-yz-a1gx-ig0-gyz}, then $z'\le z$, hence, $z=z'$ is the only
solution.

Since $[i_1+\dots+i_{y'-1},i_1+\dots+i_{y'}]$ is the leftmost interval
with ends in $N$ containing
$[l_1+\dots+l_\alpha,l_1+\dots+l_{\alpha+1}]$, and the latter lies to
the left of
$[l_1+\dots+l_{\alpha+\beta+1},l_1+\dots+l_{\alpha+\beta+2}]$,
contained in the rightmost interval
$[i_1+\dots+i_{z'-1},i_1+\dots+i_{z'}]$, we deduce that $y'\le z'$.

If $y'=z'$, then $y'\le y<z\le z'$ can not be satisfied, hence,
\eqref{eq-iiyll-lliiy}--\eqref{eq-yz-a1gx-ig0-gyz} has no solutions
$(y,z)$. On the other hand, for $q=y'=z'$ the interval
$[l_1+\dots+l_\alpha,l_1+\dots+l_{\alpha+\beta+2}]$ is contained in
$[i_1+\dots+i_{q-1},i_1+\dots+i_q]$, that is, \eqref{eq-iill-llii}
holds. Only $l_{\alpha+1}$ and $l_{\alpha+\beta+2}$ might vanish,
both are contained in
$[l_1+\dots+l_\alpha,l_1+\dots+l_{\alpha+\beta+2}]$, hence, $i_\gamma$
might vanish only for $\gamma=q$, that is, \eqref{eq-a1gx-ig0-gq}
holds. Therefore, $q=y'$ satisfies conditions
\eqref{eq-iill-llii}--\eqref{eq-a1gx-ig0-gq}. This solution is unique,
since if \eqref{eq-iill-llii} is satisfied for $q=q'$, then
\eqref{eq-iiyll-lliiy} holds for $y=q'$.

If $y'<z'$, then $y=y'$, $z=z'$ is the only solution of system of
conditions \eqref{eq-iiyll-lliiy}--\eqref{eq-yz-a1gx-ig0-gyz}. This is
proved by examining the four cases which arise from the alternatives in the
two paragraphs that follow \eqref{eq-yz-a1gx-ig0-gyz}. Let us prove
that there are no solutions $q$ of the system of conditions
\eqref{eq-iill-llii}--\eqref{eq-a1gx-ig0-gq}. Suppose $q$ satisfies
these conditions, then $y=q$ satisfies \eqref{eq-iiyll-lliiy} and $z=q$
satisfies \eqref{eq-yz-a1gx-ig0-gyz}. Therefore, $y'\le q\le z'$,
$i_1+\dots+i_{y'-1}=i_1+\dots+i_{q-1}$ and
$i_1+\dots+i_q=i_1+\dots+i_{z'}$. Due to \eqref{eq-a1gx-ig0-gq} there
exists no more than one $\gamma$ such that $i_\gamma=0$. Thus, two
possibilities exist: either $y'=q<q+1=z'$, or $y'=q-1<q=z'$. In the
first case \eqref{eq-iizll-lliiz} and \eqref{eq-iill-llii} imply
\[ i_1+\dots+i_q \le l_1+\dots+l_{\alpha+\beta+1} \le
l_1+\dots+l_{\alpha+\beta+2} \le i_1+\dots+i_q, \]
hence, $i_1+\dots+i_q=l_1+\dots+l_{\alpha+\beta+1}$ and
$l_{\alpha+\beta+2}=0$. It follows that $i_{q+1}=0$, which contradicts
to \eqref{eq-a1gx-ig0-gq}. In the second case \eqref{eq-iill-llii} and
\eqref{eq-iiyll-lliiy} imply
\[ i_1+\dots+i_{q-1} \le l_1+\dots+l_\alpha \le
l_1+\dots+l_{\alpha+1} \le i_1+\dots+i_{q-1}, \]
hence, $i_1+\dots+i_{q-1}=l_1+\dots+l_\alpha$ and $l_{\alpha+1}=0$. It
follows that $i_q=0$. From \eqref{eq-iill-llii} we deduce that
$l_{\alpha+2}+\dots+l_{\alpha+\beta+2}=0$, which implies $i_{q+1}=0$
and this contradicts \eqref{eq-a1gx-ig0-gq}. We conclude that each
term \eqref{eq-llllllll-ffrggphh} either occurs in a unique term
\eqref{eq-of-ofRoh-oh-m} or in a unique term
\eqref{eq-of-of-or-og-og-op-oh-oh-m}. Therefore,
\eqref{eq-mrpH-Rm-mrpmM} is proven.

Since \eqref{eq-rpRn-ffrggphh} is proven, it implies the statement for
the transformation $(r\tens p)R$.
\end{proof}

Denote by $\ck A_\infty'$ the non-2-unital $\ck$\n-2-category, whose
objects are pairs $(\ca,\ca')$, consisting of an \ainf-category $\ca$
and a full \ainf-subcategory $\ca'\subset\ca$; 1\n-morphisms
$(\ca,\ca')\to(\cb,\cb')$ are \ainf-functors $f:\ca\to\cb$ such that
$(\Ob\ca')f\subset\Ob\cb'$;
\[ \ck A_\infty'\bigl((\ca,\ca'),(\cb,\cb')\bigr)(f,g) =
\bigl(A_\infty(\ca,\cb)(f,g),m_1\bigr), \]
and the operations are induced by those of $\ck A_\infty$.

\begin{corollary}\label{cor-D-K2-functor}
The following assignment defines a strict $\ck$\n-2-functor
\begin{align*}
\Dr: \ck A_\infty' & \longrightarrow \ck A_\infty, \\
(\ca,\ca') & \longmapsto \Dr(\ca|\ca'), \\
f: (\ca,\ca') \to (\cb,\cb') & \longmapsto
\overline{f}: \Dr(\ca|\ca') \to \Dr(\cb|\cb'), \\
\bigl(sA_\infty((\ca,\ca'),(\cb,\cb'))(f,g),B_1\bigr) & \longrightarrow
\bigl(sA_\infty(\Dr(\ca|\ca'),\Dr(\cb|\cb'))
(\overline{f},\overline{g}),\overline{B}_1\bigr),
\qquad r\mapsto\overline{r}.
\end{align*}
\end{corollary}

\begin{proof}
Since the coderivation
$(r\tens p)R:Ts\overline{\ca}\to Ts\overline{\cb}$ restricts to a
coderivation $(r\tens p)R:Ts\Dr(\ca|\ca')\to Ts\Dr(\cb|\cb')$ by
\propref{pro-homo-rpR-explicit}, $\Dr$ preserves the vertical
composition of 2\n-morphisms modulo homotopy by
\eqref{eq-ovrpB2-ovrovpovB2}.
\end{proof}

\begin{corollary}\label{cor-D-2-functor-Ap-A}
Let $A_\infty'$ be a non-2-unital 2\n-category, whose objects and
1\n-morphisms are the same as for $\ck A_\infty'$, and 2\n-morphisms
are equivalence classes of natural \ainf-transformations:
\[ A_\infty'\bigl((\ca,\ca'),(\cb,\cb')\bigr)(f,g)
= H^0\bigl(A_\infty(\ca,\cb)(f,g),m_1\bigr),
\]
and the operations are induced by those of $A_\infty$. Then the
following assignment defines a strict 2\n-functor
\begin{align*}
\Dr: A_\infty' & \longrightarrow A_\infty, \\
(\ca,\ca') & \longmapsto \Dr(\ca|\ca'), \\
f: (\ca,\ca') \to (\cb,\cb') & \longmapsto
\overline{f}: \Dr(\ca|\ca') \to \Dr(\cb|\cb'), \\
r:f \to g: (\ca,\ca') \to (\cb,\cb') & \longmapsto
\overline{r}:\overline{f}\to\overline{g}:\Dr(\ca|\ca')\to\Dr(\cb|\cb').
\end{align*}
\end{corollary}

The corollary follows from \corref{cor-D-K2-functor} by taking the
0\n-th cohomology.

\section{Unitality}
\begin{proposition}\label{pro-C-unital-D(C|B)-unital}
Let $\cb$ be a full subcategory of a unital \ainf-category $\cc$. Then
the \ainf-category $\Dr(\cc|\cb)$ is also unital. If $\uni^\cc$ is a
unit transformation of $\cc$, then $\Dr(\uni^\cc)$ is a unit
transformation of $\Dr(\cc|\cb)$.
\end{proposition}

\begin{proof}
The idempotent property $(\uni^\cc\tens\uni^\cc)B_2\equiv\uni^\cc$
implies by \corref{cor-D-K2-functor} that
\[ (\Dr(\uni^\cc)\tens\Dr(\uni^\cc))\overline{B}_2 \equiv
\Dr((\uni^\cc\tens\uni^\cc)B_2) \equiv \Dr(\uni^\cc), \]
so $\Dr(\uni^\cc)$ is an idempotent as well. Consider its 0\n-th
component
\[ \sS{_X}\Dr(\uni^\cc)_0 = \bigl[ \kk \rTTo^{\sS{_X}\uni^\cc_0}
s\cc(X,X) \rMono s\Dr(\cc|\cb)(X,X) \bigr]. \]
We have to prove that
\[ (\sS{_X}\Dr(\uni^\cc)_0\tens1)\overline{b}_2,\quad
(1\tens\sS{_Y}\Dr(\uni^\cc)_0)\overline{b}_2:
s\Dr(\cc|\cb)(X,Y) \to s\Dr(\cc|\cb)(X,Y) \]
are homotopy invertible.

Consider the following $\ZZ_{>0}$\n-grading of the \ainf-category
$\Dr(\cc|\cb)$
\begin{align*}
  G^k &= T^ks\cc\cap s\Dr(\cc|\cb), \qquad k\ge1, \\
  G^k(X,Y) &= \oplus_{C_1,\dots,C_{k-1}\in\Ob\cb} s\cc(X,C_1)\tens
s\cc(C_1,C_2)\tens\dots\tens s\cc(C_{k-1},Y).
\end{align*}
Denote also $C_0=X$, $C_k=Y$. The corresponding increasing filtration
\[ 0 = \Phi_0 \subset \Phi_1 \subset\dots\subset \Phi_n \subset
\Phi_{n+1} \subset\dots\subset s\Dr(\cc|\cb) \]
is made of $\Phi_n=\oplus_{k=1}^nG^k$. The $\kk$\n-linear maps
\[ \overline{b}_1=b,\quad (\sS{_X}\uni^\cc_0\tens1)\overline{b}_2,
\quad (1\tens\sS{_Y}\uni^\cc_0)\overline{b}_2:
s\Dr(\cc|\cb)(X,Y) \to s\Dr(\cc|\cb)(X,Y) \]
preserve the filtration. Consider the $\ZZ_{>0}\times\ZZ$-graded
quiver, associated with this filtration. The above maps induce on
graded components $G^k$ the following maps:
\begin{align}
d_k = \sum_{\alpha+1+\beta=k}
1^{\tens\alpha}\tens b_1\tens1^{\tens\beta} &: G^k(X,Y) \to G^k(X,Y),
\label{eq-dk-1b1-Gk-Gk} \\
(\sS{_X}\uni^\cc_0\tens1)b_2\tens1^{\tens k-1} &: G^k(X,Y) \to
G^k(X,Y), \label{eq-i1b1-Gk-Gk} \\
1^{\tens k-1}\tens(1\tens\sS{_Y}\uni^\cc_0)b_2 &: G^k(X,Y) \to
G^k(X,Y). \label{eq-11ib-Gk-Gk}
\end{align}
Let $h$, $h'$ be homotopies as in
\begin{gather*}
(1\tens\sS{_Y}\uni^\cc_0)b_2 = 1+h_{C_{k-1},Y}b_1 + b_1h_{C_{k-1},Y} :
s\cc(C_{k-1},Y) \to s\cc(C_{k-1},Y), \\
(\sS{_X}\uni^\cc_0\tens1)b_2 = -1 + h'_{X,C_1}b_1 + b_1h'_{X,C_1} :
s\cc(X,C_1) \to s\cc(X,C_1).
\end{gather*}
Using them we will present map \eqref{eq-11ib-Gk-Gk} restricted to
\( s\cc(X,C_1)\tens\dots\tens s\cc(C_{k-2},C_{k-1})\tens
s\cc(C_{k-1},Y) \)
as follows:
\begin{align*}
& 1^{\tens k-1}\tens(1\tens\sS{_Y}\uni^\cc_0)b_2
= 1^{\tens k-1}\tens(1+hb_1+b_1h) \\
&= 1 + (1^{\tens k-1}\tens h) \sum_{\alpha+1+\beta=k}
1^{\tens\alpha}\tens b_1\tens1^{\tens\beta} +
\Bigl( \sum_{\alpha+1+\beta=k}
1^{\tens\alpha}\tens b_1\tens1^{\tens\beta} \Bigr)
(1^{\tens k-1}\tens h) \\
&= 1 + (1^{\tens k-1}\tens h)d_k + d_k(1^{\tens k-1}\tens h).
\end{align*}
Let us define a $\kk$\n-linear map
$H:s\Dr(\cc|\cb)(X,Y)\to s\Dr(\cc|\cb)(X,Y)$ of degree $-1$ as a direct
sum of maps
\begin{multline*}
1^{\tens k-1}\tens h_{C_{k-1},Y}: s\cc(X,C_1)\tens\dots\tens
s\cc(C_{k-2},C_{k-1})\tens s\cc(C_{k-1},Y) \\
\to s\cc(X,C_1)\tens\dots\tens s\cc(C_{k-2},C_{k-1})\tens
s\cc(C_{k-1},Y).
\end{multline*}
Since $H$ preserves the subquivers $G^k$, it preserves also the
filtration $\Phi_n$. Therefore, the chain (with respect to
$\overline{b}_1$) map
\[ 1+N \overset{\text{def}}= (1\tens\sS{_Y}\uni^\cc_0)\overline{b}_2
- H\overline{b}_1 - \overline{b}_1H:
s\Dr(\cc|\cb)(X,Y) \to s\Dr(\cc|\cb)(X,Y)
\]
preserves the filtration and the associated map of graded complexes is
identity. Hence, $N$ has a strictly lower triangular matrix with
respect to the decomposition $s\Dr(\cc|\cb)=\oplus_{k\ge1}G^k$.
Therefore, the map $1+N$ is invertible with an inverse
$\sum_{i=0}^\infty(-N)^i$. Hence,
$(1\tens\sS{_Y}\uni^\cc_0)\overline{b}_2$ is homotopy invertible.

Similarly,
\begin{align*}
& (\sS{_X}\uni^\cc_0\tens1)b_2\tens1^{\tens k-1}
= (-1+h'b_1+b_1h')\tens1^{\tens k-1} \\
&= -1 + (h'\tens1^{\tens k-1}) \sum_{\alpha+1+\beta=k}
1^{\tens\alpha}\tens b_1\tens1^{\tens\beta} +
\Bigl( \sum_{\alpha+1+\beta=k}
1^{\tens\alpha}\tens b_1\tens1^{\tens\beta} \Bigr)
(h'\tens1^{\tens k-1}) \\
&= -1 + (h'\tens1^{\tens k-1})d_k + d_k(h'\tens1^{\tens k-1}).
\end{align*}
Define a map $H'$ as a direct sum of maps
\begin{multline*}
h'_{X,C_1}\tens1^{\tens k-1}:
s\cc(X,C_1)\tens s\cc(C_1,C_2)\tens\dots\tens s\cc(C_{k-1},Y) \\
\to s\cc(X,C_1)\tens s\cc(C_1,C_2)\tens\dots\tens s\cc(C_{k-1},Y).
\end{multline*}
Then the chain map
\[ -1+N' \overset{\text{def}}= (\sS{_X}\uni^\cc_0\tens1)\overline{b}_2
- H'\overline{b}_1 - \overline{b}_1H':
s\Dr(\cc|\cb)(X,Y) \to s\Dr(\cc|\cb)(X,Y) \]
preserves the filtration and gives $-1$ on the diagonal. Hence,
$N'$ is strictly lower triangular, and $-1+N'$ is invertible with
an inverse $-\sum_{i=0}^\infty(N')^i$. Therefore,
$(\sS{_X}\uni^\cc_0\tens1)\overline{b}_2$ is homotopy invertible.
\end{proof}

\begin{corollary}\label{cor-overline-unital}
If an \ainf-category $\cc$ is unital, then $\overline\cc$ is
unital with a unit transformation $\overline{\uni^\cc}$.
\end{corollary}

Indeed, $\overline\cc=\Dr(\cc|\cc)$.

\begin{corollary}\label{cor-underline-unital}
If an \ainf-category $\cc$ is unital, then $\underline\cc$ is
unital with a unit transformation $\underline{\uni^\cc}$.
\end{corollary}

\begin{proof}
The \ainf-functor $\bm^{-1}:\underline\cc\to\overline\cc$ is
invertible and $\overline\cc$ is unital. Hence, by
\cite[Section~8.12]{Lyu-AinfCat} $\underline\cc$ is unital and
$\bm^{-1}\uni^{\overline\cc}\bm
=\bm^{-1}\overline{\uni^\cc}\bm=\underline{\uni^\cc}$
is its unit transformation.
\end{proof}

\begin{remark}
Functors $\ju^\cc:\cc\to\Dr(\cc|\cb)$,
$\ju^\cc:\cc\to\overline\cc$,
$\underline{j}^\cc:\cc\to\underline\cc$ are unital. This follows
from \corref{cor-ju-cylinder} for $r=\uni^\cc$ and from
commutative diagram \eqref{dia-underj-bm-1}.
\end{remark}

\begin{corollary}\label{cor-overline-underline-functor-unital}
Let $i:\cc\to\ci$ be a unital \ainf-functor. Then the \ainf-functors
$\iu:\overline{\cc}\to\overline{\ci}$ and
$\underline{i}:\underline{\cc}\to\underline{\ci}$ are unital as well.
\end{corollary}

\begin{proof}
Since $i\uni^\ci\equiv\uni^\cc i$, we have
$\iu\overline{\uni^\ci}\equiv\overline{\uni^\cc}\iu$ by
\corref{cor-2fun-ovmi-A-A}. Therefore, $\iu$ is unital by
\corref{cor-overline-unital}. We have also
$\underline{i}\,\underline{\uni^\ci}
\equiv\underline{\uni^\cc}\,\underline{i}$
by \corref{cor-2fun-unmi-A-A}. Hence, $\underline{i}$ is unital by
\corref{cor-underline-unital}.
\end{proof}

\begin{corollary}\label{cor-D-functor-unital}
Let $i:\cc\to\ci$ be a unital \ainf-functor, which maps objects of a
full \ainf-subcategory $\cb\subset\cc$ to objects of a full
\ainf-subcategory $\cj\subset\ci$. Then the \ainf-functor
$\Dr(i):\Dr(\cc|\cb)\to\Dr(\ci|\cj)$ is unital as well.
\end{corollary}

\begin{proof}
Since $i\uni^\ci\equiv\uni^\cc i$, we have
$\Dr(i)\Dr(\uni^\ci)\equiv\Dr(\uni^\cc)\Dr(i)$ by
\corref{cor-D-2-functor-Ap-A}. Therefore, $\Dr(i)$ is unital by
\propref{pro-C-unital-D(C|B)-unital}.
\end{proof}

Summing up, when we restrict $\unmi$, $\ovmi$ or $\Dr$ to unital
\ainf-categories, we get strict 2\n-functors of 1\n-unital
left-2-unital ($\ck$\n-)2-categories (see Definition~A.3 and
Corollary~7.11
 \href{http://arXiv.org/abs/math.CT/0210047}{\tt arXiv:\linebreak[1]math.CT/\linebreak[1]0210047}).
When we restrict $\unmi$, $\ovmi$ or $\Dr$ further to unital
\ainf-categories and unital \ainf-functors, we get strict 2\n-functors
of ordinary 1-2-unital ($\ck$\n-)2-categories.

\section{Contractibility}
A chain complex $C$ is \emph{contractible} if $\id_C$ is
null-homotopic. We say that an \ainf-category is contractible if all
its complexes of morphisms are contractible. Such \ainf-categories
behave like categories with zero morphisms only, although
contractibility might  not be obvious. An example of this kind is
provided by $\underline\cc$ and $\overline\cc$, when $\cc$ is unital.
In this section we also collect various notions of contractibility for
\ainf-functors. For unital \ainf-functors all these definitions become
equivalent.

\begin{proposition}\label{pro-C1-C5-contract}
Let $\cb$ be a unital \ainf-category. Let $f:\ca\to\cb$ be an
\ainf-functor. Then the following conditions are equivalent:
\begin{enumerate}
\item[(C1)] For any $X\in\Ob\ca$ and any $V\in\Ob\cb$ the complex
$(s\cb(Xf,V),b_1)$ is contractible;

\item[(C2)] For any $U\in\Ob\cb$ and any $Y\in\ca$ the complex
$(s\cb(U,Yf),b_1)$ is contractible;

\item[(C3)] For any object $X$ of $\ca$ the complex
$(s\cb(Xf,Xf),b_1)$ is acyclic;

\item[(C4)] For any object $X$ of $\ca$ there is an element
$\sS{_X}v\in(s\cb)^{-2}(Xf,Xf)$ such that
$\sS{_{Xf}}\uni^\cb_0=\sS{_X}vb_1$;

\item[(C5)] $f\uni^\cb\equiv0:f\to f:\ca\to\cb$.
\end{enumerate}
\end{proposition}

\begin{proof}
Clearly, (C1)$\implies$(C3)$\implies$(C4), (C2)$\implies$(C3) and
(C5)$\implies$(C4).

(C4)$\implies$(C1):
Consider a $\kk$\n-linear map
$(\sS{_X}v\tens1)b_2:s\cb(Xf,V)\to s\cb(Xf,V)$ of degree $-1$. Its
commutator with $b_1$ is
\begin{equation*}
(\sS{_X}v\tens1)b_2b_1 + b_1(\sS{_X}v\tens1)b_2
= - (\sS{_X}vb_1\tens1)b_2 = - (\sS{_{Xf}}\uni^\cb_0\tens1)b_2
\sim 1: s\cb(Xf,V)\to s\cb(Xf,V)
\end{equation*}
by \cite[Lemma~7.4]{Lyu-AinfCat}. Therefore, $s\cb(Xf,V)$ is
contractible.

(C4)$\implies$(C2):
Consider a $\kk$\n-linear map
$(1\tens\sS{_X}v)b_2:s\cb(U,Xf)\to s\cb(U,Xf)$ of degree $-1$. Its
commutator with $b_1$ is
\begin{equation*}
(1\tens\sS{_X}v)b_2b_1 + b_1(1\tens\sS{_X}v)b_2 =
- (1\tens\sS{_X}vb_1)b_2 = - (1\tens\sS{_{Xf}}\uni^\cb_0)b_2
\sim -1: s\cb(U,Xf) \to s\cb(U,Xf)
\end{equation*}
by \cite[Lemma~7.4]{Lyu-AinfCat}. Therefore, $s\cb(U,Xf)$ is
contractible.

(C1)$\implies$(C5):
We look for an $(f,f)$\n-coderivation $v$ of degree $-2$ such that
$vb-bv=f\uni^\cc$. We choose its 0\n-th component as
$\sS{_X}v_0:\kk\to(s\cb)^{-2}(Xf,Xf)$, $1\mapsto\sS{_X}v$, where
$\sS{_X}v$ satisfies condition (C4).

Let $n$ be a positive integer. Assume that $(v_0,v_1,\dots,v_{n-1})$
are already found such that an $(f,f)$\n-coderivation
$\tilde{v}=(v_0,v_1,\dots,v_{n-1},0,0,\dots)$ of degree $-2$ satisfies
equations $\lambda_m=0$ for $m<n$, where the $(f,f)$\n-coderivation
$\lambda$ of degree $-1$ is $\lambda=f\uni^\cc-\tilde{v}b+b\tilde{v}$.
To make the induction step, we look for a map
\[ v_n: s\ca(X_0,X_1)\tens\dots\tens s\ca(X_{n-1},X_n)
\to s\cb(X_0f,X_nf), \]
such that
\begin{equation}
v_nb_1 - \sum_{q+t=n}(1^{\tens q}\tens b_1\tens1^{\tens t})v_n
= \lambda_n.
\label{eq-vb-bv-lambda-n}
\end{equation}
The identity $\lambda b+b\lambda=0$ implies
\[ \lambda_nd = \lambda_nb_1
+ \sum_{q+t=n}(1^{\tens q}\tens b_1\tens1^{\tens t})\lambda_n = 0, \]
where $d$ is the differential in the complex
\begin{equation}
\Hom\bigl(s\ca(X_0,X_1)\tens\dots\tens s\ca(X_{n-1},X_n),
s\cb(X_0f,X_nf)\bigr)
\label{eq-Hom-aaa-b}
\end{equation}
(all complexes are equipped with the differential $b_1$). Since
$s\cb(X_0f,X_nf)$ is contractible, so is complex~\eqref{eq-Hom-aaa-b}.
As it is acyclic, there exists $v_n$ such that $v_nd=\lambda_n$, that
is, \eqref{eq-vb-bv-lambda-n} holds. Induction finishes the
construction of $v$.
\end{proof}

\begin{proposition}\label{pro-C6-C9-contr}
Let $\ca$ be a unital \ainf-category. Let $f:\ca\to\cb$ be an
\ainf-functor. Then the following conditions are equivalent:
\begin{enumerate}
\item[(C6)] For all objects $X$, $Y$ of $\ca$ the chain map
$f_1:(s\ca(X,Y),b_1)\to(s\cb(Xf,Yf),b_1)$ is homotopic to 0;

\item[(C7)] For any object $X$ of $\ca$ the chain map
$f_1:(s\ca(X,X),b_1)\to(s\cb(Xf,Xf),b_1)$ is homotopic to 0;

\item[(C8)] For any object $X$ of $\ca$ we have
$H^\bull(f_1)=0:H^\bull(s\ca(X,X),b_1)\to H^\bull(s\cb(Xf,Xf),b_1)$;

\item[(C9)] For any object $X$ of $\ca$ there is an element
$\sS{_X}w\in(s\cb)^{-2}(Xf,Xf)$ such that
$\sS{_X}\uni^\ca_0f_1=\sS{_X}wb_1$.
\end{enumerate}
\end{proposition}

\begin{proof}
Clearly, (C6)$\implies$(C7)$\implies$(C8)$\implies$(C9).

(C9)$\implies$(C6):
Since $f$ and $b$ commute, we have
\begin{gather*}
(1\tens\uni^\ca_0)f_2b_1 + (1\tens\uni^\ca_0)(f_1\tens f_1)b_2 =
(1\tens\uni^\ca_0)b_2f_1+(1\tens\uni^\ca_0)(1\tens b_1+b_1\tens1)f_2, \\
b_1(1\tens\uni^\ca_0)f_2 + (1\tens\uni^\ca_0)f_2b_1
- (f_1\tens\sS{_Y}w)b_2b_1 - b_1(f_1\tens\sS{_Y}w)b_2
= (1\tens\uni^\ca_0)b_2f_1 \sim f_1
\end{gather*}
by \cite[Lemma~7.4]{Lyu-AinfCat}. Therefore, $f_1$ is homotopic to 0.
\end{proof}

\begin{proposition}\label{pro-C10-C11-unital-f}
Let $\ca$, $\cb$ be unital \ainf-categories. Let $f:\ca\to\cb$ be a
unital \ainf-functor. Then conditions (C1)--(C9) are equivalent to the
following conditions:
\begin{enumerate}
\item[(C10)] There is an isomorphism of \ainf-functors
$f\simeq\co^f:\ca\to\cb$, where $\co^f$ is defined as follows:
$X\co^f=Xf$, $\co^f_k=0$ for all $k\ge1$;

\item[(C11)] $\uni^\ca f\equiv0:f\to f:\ca\to\cb$.
\end{enumerate}
\end{proposition}

\begin{proof}
Unitality implies that (C5) and (C11) are equivalent, and that (C4) and
(C9) are equivalent.

(C5)$\implies$(C10):
Consider zero natural \ainf-transformations $0:f\to\co^f:\ca\to\cb$ and
$0:\co^f\to f:\ca\to\cb$. Their composition in one order
$0\cdot0=0:f\to f:\ca\to\cb$ is equivalent to $f\uni^\cb=\sS{_f}1s$ by
(C5). Their composition in the other order
$0\cdot0=0:\co^f\to\co^f:\ca\to\cb$ is equivalent to $\co^f\uni^\cb$.
Indeed, there exists an $(f,f)$-coderivation $w$ of degree $-2$ such
that $f\uni^\cb=wb-bw$. In particular,
$\sS{_X}(f\uni^\cb)_0=\sS{_{Xf}}\uni^\cb_0=\sS{_X}w_0b_1$. Consider the
$(\co^f,\co^f)$-coderivation $v$ of degree $-2$, given by its components
$v_0=w_0$ and $v_k=0$ for $k>0$. Then
$\sS{_X}(\co^f\uni^\cb)_0=\sS{_{Xf}}\uni^\cb_0=\sS{_X}v_0b_1$ and
\[ (\co^f\uni^\cb)_n = 0 = v_nb_1 -
\sum_{q+k+t=n} (1^{\tens q}\tens b_k\tens1^{\tens t}) v_{q+1+t}
= (vb-bv)_n \]
for $n>0$. Therefore, $\co^f\uni^\cb=vb-bv$.

(C10)$\implies$(C4):
Since $f$ is unital, isomorphic to it \ainf-functor $\co^f$ is unital.
Thus,
\[ \co^f\uni^\cb\equiv\uni^\ca\co^f = 0 : \co^f \to \co^f:\ca\to\cb. \]
Therefore, there exists an $(\co^f,\co^f)$-coderivation $v$ of degree
$-2$ such that $\co^f\uni^\cb=vb-bv$. In particular,
$\sS{_{Xf}}\uni^\cb_0=\sS{_X}(\co^f\uni^\cb)_0=\sS{_X}v_0b_1$, hence,
(C4) holds.
\end{proof}

\begin{definition}\label{def-contractible-functor-category}
Let $\ca$ be a unital \ainf-category. An \ainf-functor $f:\ca\to\cb$ is
\emph{contractible} if it satisfies equivalent conditions (C6)--(C9) of
\propref{pro-C6-C9-contr}. An \ainf-category $\ca$ is
\emph{contractible} if complexes $(s\ca(X,Y),b_1)$ are contractible for
all objects $X$, $Y$ of $\ca$.
\end{definition}

A contractible \ainf-category $\ca$ is unital. Indeed,
\(\sS{_X}\uni^\ca_0=0\) are unit elements of $\ca$. The identity
\ainf-functor \(\id:\ca\to\ca\) is contractible if and only if $\ca$ is
contractible. A unital \ainf-functor $f$ is contractible if and only if
equivalent conditions (C1)--(C11) hold.

\begin{example}\label{exam-undc-overc-contractible}
If $\cc$ is a unital \ainf-category, then $\underline\cc$,
$\overline\cc$ are contractible. Indeed, by Corollaries
\ref{cor-overline-unital} and \ref{cor-underline-unital} these
categories are unital. In particular, for all objects $X$, $Y$ of
$\cc$ the chain map
\[ 0 = (1\tens\sS{_Y}\uni^\cc_0)\underline{b}_2: s\underline\cc(X,Y)
\to s\underline\cc(X,Y) \]
is homotopy invertible. Hence,
$(s\underline\cc(X,Y),\underline{b}_1)
=(s\overline\cc(X,Y),\overline{b}_1)$
is contractible. By \propref{pro-C1-C5-contract} (C2) \ainf-categories
$\underline\cc$ and $\overline\cc$ are contractible.
\end{example}

\begin{example}
Let $\cb$ be a full subcategory of a unital \ainf-category $\cc$.
Then the \ainf-functor
$\ju'=\bigl(\cb \rMono \cc \rTTo^\ju \Dr(\cc|\cb)\bigr)$ is
contractible according to criterion (C4): for any object $X$ of $\cb$
\[ (\sS{_X}\uni^\cc_0\tens\sS{_X}\uni^\cc_0)\overline{b}_1
= (\sS{_X}\uni^\cc_0\tens\sS{_X}\uni^\cc_0)b_2
= \sS{_X}\uni^\cc_0 + \sS{_X}v_0b_1
= \sS{_X}\uni^{\Dr(\cc|\cb)}_0 + \sS{_X}v_0\overline{b}_1. \]
\end{example}

\begin{proposition}\label{pro-C0-C10'-contractible-category}
A unital \ainf-category $\ca$ is contractible if and only if the
following equivalent conditions hold:
\begin{enumerate}
\item[(C0)] $\ca$ is equivalent in $A_\infty^u$ to an \ainf-category
$\co$, such that $\co(U,V)=0$ for all objects $U$, $V$ of $\co$;

\item[(C$1'$)] For all objects $X$, $Y$ of $\ca$ the complex
$(s\ca(X,Y),b_1)$ is contractible;

\item[(C$2'$)] For any object $X$ of $\ca$ the complex
$(s\ca(X,X),b_1)$ is contractible;

\item[(C$3'$)] For any object $X$ of $\ca$ the complex
$(s\ca(X,X),b_1)$ is acyclic;

\item[(C$4'$)] For any object $X$ of $\ca$ there is an element
$\sS{_X}v\in(s\ca)^{-2}(X,X)$ such that
$\sS{_X}\uni^\ca_0=\sS{_X}vb_1$;

\item[(C$5'$)] $\uni^\ca\equiv0:\id_\ca\to\id_\ca:\ca\to\ca$;

\item[(C$10'$)] There is an isomorphism of \ainf-functors
$\id_\ca\simeq\co^{\id}:\ca\to\ca$, where $X\co^{\id}=X$,
$\co^{\id}_k=0$ for all $k\ge1$.
\end{enumerate}
\end{proposition}

\begin{proof}
Conditions (C$1'$)--(C$5'$), (C$10'$) are just conditions (C1)--(C10) for
$f=\id_\ca$, hence, they are equivalent to contractibility of $\ca$.
Notice that any $\co$ as in (C0) is strictly unital.

(C$10'$)$\implies$(C0):
Denote by $\co^\ca$ the \ainf-category, whose class of objects is
$\Ob\ca$, and $\co^\ca(X,Y)=0$ for all objects
$X,Y\in\Ob\ca=\Ob\co^\ca$. Let $\phi:\co^\ca\to\ca$,
$\psi:\ca\to\co^\ca$ be the unique \ainf-functors such that $X\phi=X$,
$X\psi=X$ for all $X\in\Ob\ca$. Then $\phi\psi=\id_\ca$ and
$\psi\phi=\co^{\id_\ca}\simeq\id_\ca$ by (C$10'$).

(C0)$\implies$(C$3'$):
There is a 2\n-functor $H^\bull:A_\infty^u\to\Cat$,
$\Ob H^\bull(\ca)=\Ob\ca$,
$H^\bull(\ca)(X,Y)$\linebreak[1]${}=H^\bull(\ca(X,Y),b_1)$. If $\ca$ is
equivalent to $\co$ in $A_\infty^u$, then $H^\bull(\ca)$ is
equivalent to $H^\bull(\co)$ in $\Cat$. Clearly, $H^\bull(\co)(U,V)=0$
for all objects $U$, $V$ of $\co$. Hence, $H^\bull(\ca)(X,Y)=0$ for all
objects $X$, $Y$ of $\ca$.
\end{proof}

\begin{corollary}\label{cor-contractible-r-equivalent-0}
If $\cb$ is contractible, then for any \ainf-category $\ca,$ any natural
\ainf-transformation $r:f\to g:\ca\to\cb$ is equivalent to 0.
\end{corollary}

\begin{remark}\label{rem-O-0-equivalent-1}
Any \ainf-category $\co$, such that $\co(U,V)=0$ for all objects $U$,
$V$ of $\co$ with non-empty $\Ob\co$ is equivalent to
1-object-1-morphism \ainf-category $\1$, such that $\Ob\1=\{*\}$ and
$\1(*,*)=0$. Indeed, choose an object $Z\in\Ob\co$. Consider
\ainf-functors $\phi:\co\to\1$, $U\mapsto*$ and $\psi:\1\to\co$,
$*\mapsto Z$. We have $\psi\phi=\id_{\1}$ and $\phi\psi$ is isomorphic
to $\id_\co$ via inverse  2\n-morphisms
$0:\phi\psi\to\id_\co:\co\to\co$ and $0:\id_\co\to\phi\psi:\co\to\co$.
\end{remark}

\begin{remark}
Let $\cc$ be strictly unital. Then $\Dr(\cc|\cb)$ has a strict
unit $\uni^{\Dr(\cc|\cb)}$ described in \secref{sec-str-unit}. On
the other hand, $\Dr(\uni^\cc)$ is a unit of $\Dr(\cc|\cb)$ as
well. Hence, $\uni^{\Dr(\cc|\cb)}\equiv\Dr(\uni^\cc)$ by
\cite[Corollary~7.10]{Lyu-AinfCat}.
\end{remark}

\section{The case of a contractible subcategory}
Taking the quotient $\Dr(\ce|\cf)$ can be interpreted as contracting
the full \ainf-subcategory $\cf\subset\ce$. If $\cf$ were already
contractible, one would expect that no further contracting is required.
And, in fact, if $\ce$ is unital, we shall prove below that
$\Dr(\ce|\cf)$ is equivalent to $\ce$. In the proof we shall construct
inductively a new \ainf-structure on $\ce$. So first of all, we
consider direct limits of \ainf-structures on a given graded
$\kk$\n-linear quiver.

\begin{lemma}\label{lem-lim-prod}
Let $\cb$ be a graded $\kk$\n-quiver. Let $\ca_k$, $k\ge1$, be a
sequence of \ainf-categories, whose underlying graded $\kk$\n-quiver is
$\cb$. Let $\ca_1\xrightarrow{\ f^1\ }\ca_2\xrightarrow{\ f^2\ }\dots$
be a sequence of \ainf-functors, such that $f^k_1=\id_{\ca_k}$ for all
$k$, and let $N_i$, $i\ge2$ be an increasing sequence of positive
integers, such that $f^k_i=0$ for $k\ge N_i$. Then there exists a
direct (2\n-)limit $\ca=\tcolimit_{f^i} \ca_i$ of this diagram, and
the structure \ainf-functors $p^k:\ca_k\to\ca$ are invertible.
\end{lemma}

\begin{proof}
If $g:\cd\to \cc$ is an \ainf-functor, such that $g_1=\id_\cd$ and
$g_i=0$ for $i=2,\dots,k$, then for any such $i$ there exists a
commutative diagram
\begin{diagram}
s\cd^{\tens i} & \rTTo^{\id_{s\cd^{\tens i}}} & s\cc^{\tens i} \\
\dTTo<{b_i} && \dTTo>{b_i} \\
s\cd & \rTTo^{\id_{s\cd}} & s\cc
\end{diagram}
which allows us to identify the \ainf-operations $b_i$, $i=1,\dots,k$ on
$\cd$ and $\cc$.

Due to this remark, we take $\cb$ as the underlying $\kk$\n-quiver of
$\ca$, we set $b_1:s\ca\to s\ca$ to be $b_1:s\ca_1\to s\ca_1$ and we
set $b_i:s\ca^{\tens i}\to s\ca$ equal
$b_i:s\ca_{N_i}^{\tens i}\to s\ca_{N_i}$. This equips $\ca$ with an
\ainf-structure. We define $p^k$ setting its $i$\n-th component equal
to $p^k_i=(f^kf^{k+1}\dots f^l)_i$ for $l=\max(k,N_i)$.

Given an \ainf-category $\cc$ and \ainf-functors $\pi^k:\ca_k\to\cc$,
$k=1,\dots$, such that $\pi^{k}=f^k\pi^{k+1}$, then there exists the
unique \ainf-functor $\pi:\ca\to\cc$, such that $\pi^k=p^k\pi$, defined
by $\pi_i=\pi^{N(i)}_i$, $i\ge1$. It shows that the constructed $\ca$
is a direct limit of the diagram $(\ca_i,f^i,i\geqslant1)$.
\end{proof}

\begin{lemma}\label{lem-iso-on-stupid}
Let $\ce$ be an \ainf-category and let $\cf$ be its full
\ainf-subcategory such that the complex of $\kk$\n-modules
$(s\ce(X,Y),b_1)$ is contractible provided at least one of $X$, $Y$
belongs to $\Ob\cf$. Denote by $D_n(\ce|\cf)$, $n=2,3,\dots$ the
$\kk$\n-submodule in $(s\ce)^{\tens n}$, which is a sum of
$s\ce(X_0,X_1)\tens s\ce(X_1,X_2)\tens\dots\tens
s\ce(X_{n-2},X_{n-1})\tens s\ce(X_{n-1},X_n)$,
such that at least for one $i=0,\dots,n$ object $X_i$ belongs to
$\Ob\cf$. Then there exists an invertible \ainf-functor
$g:\ce\to\ce_\cf$, such that \ainf-category $(\ce_\cf,b')$ coincides
with $\ce$ as a graded differential (with respect to $b_1$)
$\kk$\n-quiver (that is, $g_1=\id_\ce$) and $(D_n(\ce|\cf))b'_n=0$ for
any $n>1$.
\end{lemma}

\begin{proof}
We construct a chain of \ainf-isomorphisms $f^i:\ce_i\to\ce_{i+1}$,
$i=1,\dots$, $\ce_1=\ce$ as in \lemref{lem-lim-prod} and then we set
$\ce_\cf=\tcolimit_{f^i}\ce_i$ and $g=p^1$. For the constructed $\ce_j$
and $f^j$, $j=1,\dots$ the following will hold:
1) $b_k|_{D_k(\ce|\cf)}=0$ holds in $\ce_j$ for $k=2,\dots,j$;
2) $f_i^j=0$ for $i\ne1,j+1$.

Given an \ainf-category $\ce_{k}$ and a $\kk$\n-quiver morphism
$f^{k}=(\id_\ce,0,\dots,0,f^{k}_{k+1},0,\dots):Ts\ce_{k}\to s\ce$
of degree $0$, they define (following \remref{rem-move-diff}) a unique
\ainf-category structure $\ce_{k+1}$ on the graded $\kk$\n-quiver
$\ce$ such that $f^{k}$ turns into an \ainf-functor
$f^{k}:\ce_{k}\to\ce_{k+1}$. Assume $f^j:\ce_{j}\to\ce_{j+1}$,
$j=1,\dots,k-1$ are constructed (the reasoning is valid for $k=1$ too).

Let us fix for any sequence $X_0,X_1,\dots,X_n$ as in the hypothesis of
the lemma an index $l(X_0,\dots,X_n)$, such that $X_{l(X_0,\dots,X_n)}$
belongs to $\Ob\cf$.
Any choice of $f^{k}_{k+1}:T^{k+1}s\ce_k\to s\ce$ determines by
\remref{rem-move-diff} an \ainf-category $\ce_{k+1}$ with the
operations \(b_p=b^{\ce_{k+1}}_p\). Notice that the conditions
$b^{\ce_{k+1}}_2|_{D_2(\ce|\cf)}=0$, \dots,
$b^{\ce_{k+1}}_k|_{D_k(\ce|\cf)}=0$ hold
automatically: in view of $f_1^{k}=\id_\ce$ in $\ce_{k}$ and
$f_j^{k}=0$, $1<j\le i\le k$ the $i$\n-th condition
\eqref{eq-fff-b-b-f} shows, that $b_1,\dots,b_k$ in $\ce_{k}$ and
$\ce_{k+1}$ coincide. The $(k+1)$\n-th condition \eqref{eq-b-b-0} for
$\ce_{k}$ on $D_{k+1}(\ce|\cf)$ turns into
\[ \sum_{r+1+t=k+1}
(1^{\tens r}\tens b_1\tens1^{\tens t})b_{k+1}+b_{k+1}b_1=0 \]
(for any other summand in the sum
$\sum_{r+n+t=k+1}(1^{\tens r}\tens b_n\tens1^{\tens t})b_{r+1+t}$
either the first or the second factor vanishes by induction). On the
other hand, by the induction assumptions the $({k+1})$\n-th condition
\eqref{eq-fff-b-b-f} turns on $D_{k+1}(\ce|\cf)$ into
\begin{equation}\label{eq-con-on-f}
(f_{1}^{k})^{\tens (k+1)}b_{k+1}+f_{k+1}^{k}b_1=
\sum_{r+1+t={k+1}} (1^{\tens r}\tens b_1\tens1^{\tens t})f_{{k+1}}^{k}
+b_{k+1}f_{1}^{k}: sD_{k+1}(\ce|\cf)\to s\ce_{k+1}.
\end{equation}
We consider the condition ``$b_{k+1}|_{D_{k+1}(\ce|\cf)}=0$ holds in
the \ainf-category $\ce_{k+1}$'' as an equation with respect to
$f^{k}_{k+1}$. Denote $l=\min\{l(X_0,\dots,X_{k+1}),k\}$. Choose a
contracting homotopy, i.e. $\kk$\n-module morphism
$h:s\ce(X_l,X_{l+1})\to s\ce(X_l,X_{l+1})$ of degree $-1$ such that
$hb_1+b_1h=\id_{s\ce(X_l,X_{l+1})}$. We define $f_{k+1}^{k}$ on
$D_{k+1}(\ce|\cf)(X_0,\dots,X_{k+1})$ as
\[ f_{k+1}^{k}=-(1^{\tens l}\tens h\tens 1^{\tens (k-l)})b_{k+1}. \]
On $s\ce(X_0,X_1)\tens\dots\tens s\ce(X_{n-1},X_n)$ such that all
$X_i\notin\Ob\cf,$ we set $f_{k+1}^{k}=0$.

Compare equations~\eqref{eq-con-on-f} restricted to
$D_{k+1}(\ce|\cf)(X_0,\dots,X_{k+1})$ with the following computation
\begin{multline*}
f^k_{k+1}b_1=-(1^{\tens l}\tens h\tens 1^{\tens(k-l)})b_{k+1}b_1
=(1^{\tens l}\tens h\tens 1^{\tens (k-l)})
\sum_{r+1+t={k+1}} (1^{\tens r}\tens b_1\tens1^{\tens t})b_{{k+1}} \\
=\sum_{r+1+t={k+1}} (1^{\tens r}\tens b_1\tens1^{\tens t})
(-(1^{\tens l}\tens h\tens 1^{\tens (k-l)}))b_{k+1} \\
+(1^{\tens l}\tens b_1h\tens 1^{\tens (k-l)})b_{k+1}
+(1^{\tens l}\tens hb_1\tens 1^{\tens (k-l)})b_{k+1}
=\sum_{r+1+t={k+1}}
(1^{\tens r}\tens b_1\tens1^{\tens t})f^k_{k+1}+b_{{k+1}}.
\end{multline*}
We deduce that in $\ce_{k+1}$ the restriction of $b_{k+1}$ to
$D_{k+1}(\ce|\cf)$ vanishes. Now the lemma follows from the definition
of the limit morphism $g:\ce\to\tcolimit_{f^i}\ce_i$.
\end{proof}

\begin{lemma}\label{lem-drin-splits}
Let $\cf\subset\ce_\cf$ be a full \ainf-subcategory such that
$b_n\big|_{D_n(\ce_\cf|\cf)}$ vanishes for all $n\ge2$. Then the
canonical strict embedding $\ju:\ce_\cf\to\Dr(\ce_\cf|\cf)$ admits a
splitting strict \ainf-functor $\pi:\Dr(\ce_\cf|\cf)\to\ce_\cf,$ that
is, $\ju\pi=\id_{\ce_\cf}$. Its first component is the projection
\[ \pi_1 = \bigl( s\Dr(\ce_\cf|\cf) \hookrightarrow T^+s\ce_\cf
\rTTo^{\pr_1} s\ce_\cf \bigr). \]
\end{lemma}

\begin{proof}
Denote $\ca=\Dr(\ce_\cf|\cf)\cap(s^{-1}T^{\ge2}s\ce_\cf)=\Ker\pi_1$.
Then $\Dr(\ce_\cf|\cf)=\ce_\cf\oplus\ca$ as a graded $\kk$\n-quiver and
$s\ca\subset\oplus_{n\ge2}D_n(\ce|\cf)$. Let us check that $\pi$ is
an \ainf-functor, that is, $\pi_1^{\tens k}b_k=\overline{b}_k\pi_1$ for
all $k\ge1$.

Notice that
$\ca(X,Y)=\oplus^{n\ge2}_{X_1,\dots,X_{n-1}\in\Ob\cf}
s\ce(X,X_1)\tens\dots\tens s\ce(X_{n-1},Y)$,
and each substring of such a tensor product of length $k\ge2$ is in
$D_k(\ce|\cf)$. The restriction to $s\ce_\cf^{\tens k}$ of the equation
$\pi_1^{\tens k}b_k=\overline{b}_k\pi_1$ follows from
\eqref{eq-Bn-sum-1bm1} or \corref{cor-mu-ainf-functor}. Both sides,
restricted to $s\cx_1\tens\dots\tens s\cx_k$ where $\cx_j$ is
$\ce_\cf$ or $\ca$, vanish if at least one $\ca$ is present. Indeed,
$\overline{b}_k$ vanishes in that case if $k>1$ due to
\eqref{eq-Bn-sum-1bm1}. For $k=1,$
$(s\ca)\overline{b}_1=(s\ca)b\subset s\ca$ as equation
$\overline{b}_1\big|_\ca=b\big|_\ca=
\sum_{q+1+t=n}1^{\tens q}\tens b_1\tens1^{\tens t}:s\ca\to s\ca$
shows. Other claims are clear.
\end{proof}

\begin{proposition}\label{pro-D(E|contractible-F)-E}
Let $\cf$ be a contractible full subcategory of a unital \ainf-category
$\ce$. Then there exists a quasi-inverse to the canonical strict
embedding $\ju^\ce:\ce\to\Dr(\ce|\cf)$ unital \ainf-functor
$\pi^\ce:\Dr(\ce|\cf)\to\ce$ such that $\ju^\ce\pi^\ce=\id_\ce$. In
particular, $\Dr(\ce|\cf)$ is equivalent to $\ce$.
\end{proposition}

\begin{proof}
First of all, we prove the statements for the full embedding
$\cf\subset\ce_\cf$ constructed in \lemref{lem-iso-on-stupid}.
Since $\ce$ is unital and $g:\ce\to\ce_\cf$ is invertible, the
\ainf-category $\ce_\cf$ is unital by
\cite[Section~8.12]{Lyu-AinfCat}. Let us prove that the
\ainf-functor $\pi^{\ce_\cf}=\pi:\Dr(\ce_\cf|\cf)\to\ce_\cf$ from
\lemref{lem-drin-splits} is unital, that is,
$\pi\uni^{\ce_\cf}\equiv\Dr(\uni^{\ce_\cf})\pi$.

We look for a 3\n-morphism
\[ v: \pi\uni^{\ce_\cf} \to \Dr(\uni^{\ce_\cf})\pi:
\pi \to \pi:  \Dr(\ce_\cf|\cf) \to \ce_\cf. \]
Since
$(\pi\uni^{\ce_\cf})_0=\uni^{\ce_\cf}_0=(\Dr(\uni^{\ce_\cf})\pi)_0$, we
may take $v_0=0$. Let us proceed by induction. Assume that we have
already found components $(v_0,v_1,\dots,v_{n-1})$ of $v$ such that
$v_m$ vanishes on $(s\ce_\cf)^{\tens m}$ for all $m<n$. Define a
$(\pi,\pi)$\n-transformation
$\tilde{v}=(v_0,v_1,\dots,v_{n-1},0,0,\dots)$ by these components.
Denote by $\lambda$ the $(\pi,\pi)$\n-transformation
$\pi\uni^{\ce_\cf}-\Dr(\uni^{\ce_\cf})\pi
-\tilde{v}b+\overline{b}\tilde{v}$.
Our assumption is that $\lambda_m=0$ for $m<n$. Clearly,
$\lambda b+\overline{b}\lambda=0$. This implies
\[ 0 = (\lambda b+\overline{b}\lambda)_n = \lambda_nb_1 +\sum_{q+1+t=n}
(1^{\tens q}\tens\overline{b}_1\tens1^{\tens t})\lambda_n, \]
that is,
$\lambda_n\in\Hom^{-1}\bigl((s\Dr(\ce_\cf|\cf))^{\tens n},s\ce_\cf\bigr)$
is a cocycle. We wish to prove that it is a coboundary of an element
$v_n\in\Hom^{-2}\bigl((s\Dr(\ce_\cf|\cf))^{\tens n},s\ce_\cf\bigr)$,
that is,
\[ \lambda_n = v_nb_1 - \sum_{q+1+t=n}
(1^{\tens q}\tens\overline{b}_1\tens1^{\tens t})v_n = v_nd. \]
We have
$(s\Dr(\ce_\cf|\cf))^{\tens n}
=\oplus_{\cx_i\in\{\ce_\cf,\ca\}}s\cx_1\tens\dots\tens s\cx_n$.
If at least one $\ca$ is present in $s\cx_1\tens\dots\tens s\cx_n$,
then this summand is contractible, hence,
$\Hom(s\cx_1\tens\dots\tens s\cx_n,s\ce_\cf)$ is contractible. Therefore,
there are such $v_n'\in\Hom^{-2}(s\cx_1\tens\dots\tens s\cx_n,s\ce_\cf)$
that $v_n'd=\lambda_n$ on $s\cx_1\tens\dots\tens s\cx_n$. It remains to
look at the case $s\cx_1\tens\dots\tens s\cx_n=(s\ce_\cf)^{\tens n}$.
By restriction to this submodule, we have
\begin{gather*}
(\pi\uni^{\ce_\cf})_n = \uni^{\ce_\cf}_n =
(\Dr(\uni^{\ce_\cf})\pi)_n:(s\ce_\cf)^{\tens n}\to s\ce_\cf, \\
(\tilde{v}b-\overline{b}\tilde{v})_n=0:(s\ce_\cf)^{\tens n}\to s\ce_\cf.
\end{gather*}
Therefore, $\lambda_n$ vanishes on $(s\ce_\cf)^{\tens n}$, and
$v_n\big|_{(s\ce_\cf)^{\tens n}}=0$ satisfies the equation and the
induction assumptions.

Define the following natural \ainf-transformations
\begin{align*}
r &: \id \to \pi\ju: \Dr(\ce_\cf|\cf) \to \Dr(\ce_\cf|\cf), \\
p &: \pi\ju \to \id: \Dr(\ce_\cf|\cf) \to \Dr(\ce_\cf|\cf)
\end{align*}
via its components, restricted to
$s\cx_1\tens\dots\tens s\cx_n\subset(s\Dr(\ce_\cf|\cf))^{\tens n}$,
namely, $r_k=\uni^{\ce_\cf}_k$ and $p_k=\uni^{\ce_\cf}_k$ if
$s\cx_1\tens\dots\tens s\cx_n=(s\ce_\cf)^{\tens n}$, and $r_k=0$,
$p_k=0$ otherwise. We have to check the equation
$r\overline{b}+\overline{b}r=0$. Its restriction to
$(s\ce_\cf)^{\tens n}$ holds because on this submodule $\overline{b}$
can be replaced with $b$ and $r$ with $\uni$. If
$\{\cx_1,\dots,\cx_n\}$ contains $\ca,$ then all terms in the following
sums vanish on $s\cx_1\tens\dots\tens s\cx_n$, hence,
\[ \sum_{q+k+t=n}
(1^{\tens q}\tens r_k\tens(\pi\ju)_1^{\tens t})\overline{b}_{q+1+t} +
\sum_{q+k+t=n}(1^{\tens q}\tens\overline{b}_k\tens1^{\tens t})r_{q+1+t}
= 0.
\]
In the same way we prove that $p\overline{b}+\overline{b}p=0$. Thus,
$r$ and $p$ are, indeed, natural \ainf-transformations.

Let us prove now that $r$ and $p$ are inverse to each other
2\n-morphisms, that is,
\begin{align*}
(r\tens p)B_2 &\equiv \Dr(\uni^{\ce_\cf}): \id\to\id:
\Dr(\ce_\cf|\cf) \to \Dr(\ce_\cf|\cf), \\
(p\tens r)B_2 &\equiv \pi\ju\Dr(\uni^{\ce_\cf}):
\pi\ju \to \pi\ju: \Dr(\ce_\cf|\cf) \to \Dr(\ce_\cf|\cf).
\end{align*}
We look for a 3\n-morphism
\[ v: (r\tens p)B_2 \to \Dr(\uni^{\ce_\cf}): \id\to\id:
\Dr(\ce_\cf|\cf) \to \Dr(\ce_\cf|\cf). \]
Let us look at the restriction of the equation
\begin{equation}
(r\tens p)B_2-\Dr(\uni^{\ce_\cf})=v\overline{b}-\overline{b}v
\label{eq-rpB2-Dr11-vb-bv}
\end{equation}
to $(s\ce_\cf)^{\tens n}$. First of all, $(\pi\ju)_1=1$ on $s\ce_\cf$.
Summands $\overline{b}_k$, contained in $B_2$, are applied to elements
of $(s\ce_\cf)^{\tens k}$ only. Hence, they can be replaced with $b_k$.
Therefore, $B^{\Dr(\ce_\cf|\cf)}_2$ is replaced with $B^{\ce_\cf}_2$ on
$(s\ce_\cf)^{\tens n}$. The problem of finding a 3\n-morphism
\[ w: (\uni^{\ce_\cf}\tens\uni^{\ce_\cf})B_2 \to \uni^{\ce_\cf}:
\id \to \id: \ce_\cf \to \ce_\cf
\]
is solvable. We set the restriction of $v_n$ to $(s\ce_\cf)^{\tens n}$
equal $v_n=w_n:(s\ce_\cf)^{\tens n}\to s\ce_\cf$ and it solves
\eqref{eq-rpB2-Dr11-vb-bv} on this submodule. The restriction of
equation~\eqref{eq-rpB2-Dr11-vb-bv} to $s\cx_1\tens\dots\tens s\cx_n$
that contains factor $\ca$, can be solved by induction due to
contractibility of $s\cx_1\tens\dots\tens s\cx_n$ as above. Thus $v$ is
constructed.

Similarly a 3\n-morphism
\[ u: (p\tens r)B_2 \to \pi\ju\Dr(\uni^{\ce_\cf}):
\pi\ju \to \pi\ju: \Dr(\ce_\cf|\cf) \to \Dr(\ce_\cf|\cf) \]
is constructed.

The property $\ju\pi=\id_{\ce_\cf}$ is proved in
\lemref{lem-drin-splits}.

Now we turn to the general case. The invertible \ainf-functor
$g:\ce\to\ce_\cf$, constructed in \lemref{lem-iso-on-stupid} is
the identity on objects. Denoting $\pi^{\ce_\cf}=\pi$ as above,
$g'=g\big|_\cf:\cf\to\cf$, and
$\pi^\ce=\overline{g}\pi^{\ce_\cf}g^{-1}$, we get a diagram
\begin{diagram}
\cf & \rMono & \ce & \pile{\rTTo^{\ju^\ce}\\ \lTTo_{\pi^\ce}}
& \Dr(\ce|\cf) \\
\dTTo<{g'} && \dTTo>g && \dTTo>{\overline{g}} \\
\cf & \rMono & \ce_\cf &
\pile{\rTTo^{\ju^{\ce_\cf}}\\ \lTTo_{\pi^{\ce_\cf}}} & \Dr(\ce_\cf|\cf)
\end{diagram}
All the required properties of $\pi^\ce$ follow immediately from those
of $\pi^{\ce_\cf}$.
\end{proof}

\subsection{Reducing a full contractible subcategory to 0.}
Let $\cf$ be a full contractible subcategory of a unital \ainf-category
$\ce$. Let us consider another \ainf-category $\ce\plquot\cf$, whose
class of objects is $\Ob\ce$. Here $\plquot$ stands for the plain
quotient. The morphisms are $\ce\plquot\cf(X,Y)=\ce(X,Y)$, if
$X,Y\in\Ob\ce-\Ob\cf$ and $\ce\plquot\cf(X,Y)=0$ otherwise. The
component of the differential for $\ce\plquot\cf$
\[ b_n:
s\ce\plquot\cf(X_0,X_1)\tens\dots\tens s\ce\plquot\cf(X_{n-1},X_n)
\to s\ce\plquot\cf(X_0,X_n) \]
equals $b_n$ for $\ce$ if $X_0,\dots,X_n\notin\Ob\cf$, and vanishes
otherwise.

There is a strict embedding $e:\ce\plquot\cf\to\ce$, which is the identity
on objects, $e_1=\id:s\ce\plquot\cf(X,Y)\to s\ce(X,Y)$ if
$X,Y\notin\Ob\cf$ and vanishes otherwise.The identity
$e_1^{\tens n}b^\ce_n=b^{\ce\plquot\cf}_ne_1$ is obvious.

If $\ce$ is strictly unital with the strict unit $\uni^\ce$, then
$\ce\plquot\cf$ is strictly unital with the strict unit
$\uni^{\ce\plquot\cf}$, defined as follows. Its 0\n-th component is
$\sS{_X}\uni^{\ce\plquot\cf}_0=\sS{_X}\uni^\ce_0$ if $X\notin\Ob\cf$,
and vanishes otherwise.

Let us consider the general case of a unital $\ce$. Each complex
$(s\ce(X,Y),b_1)$ is contractible if $X$ or $Y$ is an object of $\cf$
due to \propref{pro-C1-C5-contract} (C1), (C2). Therefore,
$e_1:s\ce\plquot\cf(X,Y)\to s\ce(X,Y)$ is homotopy invertible for all
pairs $X$, $Y$ of objects of $\ce$. Consider the following data:
identity map $h=\id:\Ob\ce\to\Ob\ce\plquot\cf$ and $\kk$\n-linear maps
$\sS{_X}r_0=\sS{_X}p_0=\sS{_X}\uni^\ce_0:\kk\to(s\ce)^{-1}(X,X)$.
Clearly, $\sS{_X}\uni^\ce_0b_1=0$ and
$(\sS{_X}\uni^\ce_0\tens\sS{_X}\uni^\ce_0)b_2-\sS{_X}\uni^\ce_0 \in\im
b_1$.
Therefore, the hypotheses of Theorem~8.8 of \cite{Lyu-AinfCat} are
satisfied. By this theorem we conclude that $\ce\plquot\cf$ is unital
and $e:\ce\plquot\cf\to\ce$ is a unital \ainf-equivalence.

\section{\texorpdfstring
{An $A_\infty$-functor related by an $A_\infty$-transformation to a given $A_\infty$-functor}
{An A8-functor related by an A8-transformation to a given A8-functor}}
 \label{sec-constr-inj-res-ainf-fun}
Given an \ainf-functor $f$ and the 0\n-th component $r_0$ of a natural
\ainf-transformation $r:f\to g$, we construct the \ainf-functor $g$ and
extend $r_0$ to the whole \ainf-transformation $r$. We do it under
additional assumptions on $r_0$ which are satisfied, for instance,
when $r_0$ is invertible. In the next section we apply this
construction to the case of the quasi-isomorphisms $r_0$.

\subsection{Assumptions.}\label{sec-inj-res-ainf-fun}
Let $\cb$, $\cc$ be \ainf-categories, let $f:\cb\to\cc$ be an
\ainf-functor and let $g:\Ob\cb\to\Ob\cc$ be a map. Assume that for
each object $X\in\Ob\cb,$ there is an element $r_X\in\cc^0(Xf,Xg)$ such
that $r_Xsb_1=0$. For any object $Y\in\Ob\cb,$ this element determines a
chain map
\[ (r_Xs\tens1)b_2: s\cc(Xg,Yg) \to s\cc(Xf,Yg), \qquad
p\mapsto (-)^p(r_Xs\tens p)b_2. \]
Finally, we assume that for any chain complex of $\kk$\n-modules of the
form
$N=s\cb(X_0,X_1)\tens_\kk s\cb(X_1,X_2)\tens_\kk\dots\tens_\kk
s\cb(X_{n-1},X_n)$, $n\ge0$,
the following chain map
\begin{equation}
u = \Hom(N,(r_Xs\tens1)b_2): \Hom_\kk^\bull(N,s\cc(Xg,Yg)) \to
\Hom_\kk^\bull(N,s\cc(Xf,Yg))
\label{eq-Hom-qis-u}
\end{equation}
is a quasi-isomorphism. For $n=0,$ we have $N=\kk$ and the 0\n-th
condition means that $(r_Xs\tens1)b_2$ is a quasi-isomorphism.

\begin{proposition}\label{prop-ainf-fr-g-tr-r}
Under the above assumptions, the map $g:\Ob\cb\to\Ob\cc$ extends to an
\ainf-functor $g:\cb\to\cc$. There exists a natural
\ainf-transformation $r:f\to g:\cb\to\cc$ such that its 0\n-th
component is $r_0:\kk\to s\cc(Xf,Xg)$, $1\mapsto r_Xs$.
\end{proposition}

All statements in this section (existence of the \ainf-functor $g$ and
the natural \ainf-transformation $r$, their uniqueness in a certain
sense, unitality of $g$ and invertibility of $r$) are proved in a
similar fashion, using acyclicity of the cone of the quasi-isomorphism
$u$.

\begin{proof}
The components $g_0=0$ and $r_0$ are already known. Let us build the
remaining components by induction. Assume that we have already found
components $g_m$, $r_m$ of the sought for $g$, $r$ for $m<n$, such that the
equations
\begin{gather*}
gb\pr_1=bg\pr_1:s\cb(X_0,X_1)\tens_\kk\dots\tens_\kk s\cb(X_{m-1},X_m)
\to s\cc(X_0g,X_mg), \\
(rb+br)\pr_1=0:s\cb(X_0,X_1)\tens_\kk\dots\tens_\kk s\cb(X_{m-1},X_m)
\to s\cc(X_0f,X_mg)
\end{gather*}
are satisfied for all $m<n$. Under these assumptions, we will find such
$g_n$, $r_n$ that the above equations are satisfied for $m=n$. Let us
write down these equations explicitly. The terms which contain unknown
maps $g_n$, $r_n$ are singled out on the left hand side. The right hand
side consists of already known terms:
\begin{multline}
-g_nb_1 + \sum_{q+1+t=n}(1^{\tens q}\tens b_1\tens1^{\tens t})g_n \\
= \sum_{l>1;i_1+\dots+i_l=n}
(g_{i_1} \tens g_{i_2} \tens\dots\tens g_{i_l}) b_l
- \sum_{k>1;q+k+t=n} (1^{\tens q}\tens b_k\tens1^{\tens t}) g_{q+1+t},
\label{eq-gb-bg-Tnb}
\end{multline}
\begin{multline}
r_nb_1 + \sum_{q+1+t=n}(1^{\tens q}\tens b_1\tens1^{\tens t})r_n
+ (r_0\tens g_n)b_2 \\
= - \hspace*{-3mm}
\sum_{\substack{k<n;(q,k,t)\ne(0,0,1)\\i_1+\dots+i_q+k+j_1+\dots+j_t=n}}
\hspace*{-3mm} (f_{i_1}\tens\dots\tens f_{i_q}\tens r_k\tens
 g_{j_1}\tens\dots\tens g_{j_t})b_{q+1+t}
- \sum_{\substack{k>1\\q+k+t=n}}
(1^{\tens q}\tens b_k\tens1^{\tens t})r_{q+1+t}.
\label{eq-rb+br-TnB}
\end{multline}
Let us prove that there exist $\kk$\n-linear maps
$g_n:s\cb(X_0,X_1)\tens_\kk\dots\tens_\kk s\cb(X_{n-1},X_n)\to
s\cb(X_0g,X_ng)$,
$r_n:s\cb(X_0,X_1)\tens_\kk\dots\tens_\kk s\cb(X_{n-1},X_n)\to
s\cb(X_0f,X_ng)$
which solve the above equations.

Since the map $u=\Hom(N,(r_{X_0}s\tens1)b_2)$ from \eqref{eq-Hom-qis-u}
is a quasi-isomorphism, $\Cone(u)$ is acyclic. As a differential graded
$\kk$\n-module
\begin{align*}
\Cone(u) &= \Hom_\kk^\bull(N,s\cc(X_0f,X_ng)) \oplus
\Hom_\kk^\bull(N,s\cc(X_0g,X_ng))[1], \\
(v,p)d &= (vd+pu,-pd),
\end{align*}
where $N=s\cb(X_0,X_1)\tens_\kk\dots\tens_\kk s\cb(X_{n-1},X_n)$.
Denote by $\lambda_n\in\Hom_\kk^1(N,s\cc(X_0g,X_ng))$ the right hand
side of \eqref{eq-gb-bg-Tnb} and by
$\nu_n\in\Hom_\kk^0(N,s\cc(X_0f,X_ng))$ the right hand side of
\eqref{eq-rb+br-TnB}. Equations \eqref{eq-gb-bg-Tnb} and
\eqref{eq-rb+br-TnB} mean that
$(r_n,g_n)d=(\nu_n,\lambda_n)\in\Cone^0(u)$. Since $\Cone(u)$ is
acyclic, such a pair $(r_n,g_n)\in\Cone^{-1}(u)$ exists if and only if
$(\nu_n,\lambda_n)\in\Cone^0(u)$ is a cycle, that is, equations
$-\lambda_nd=0$, $\nu_nd+\lambda_nu=0$ are satisfied. Let us verify
them now.

Introduce a cocategory homomorphism $\tilde{g}:Ts\cb\to Ts\cc$ by its
components
$(g_1,\dots,g_{n-1}$, $0$, $0,\dots)$ (these are already known).
The map $\lambda=\tilde{g}b-b\tilde{g}$ is a
$(\tilde{g},\tilde{g})$-coderivation. Its components
$(\tilde{g}b-b\tilde{g})_k$ vanish for $0\le k\le n-1$, and
\[ (\tilde{g}b-b\tilde{g})_n = \sum_{l>1;i_1+\dots+i_l=n}
(g_{i_1} \tens g_{i_2} \tens\dots\tens g_{i_l}) b_l
- \sum_{k>1;q+k+t=n} (1^{\tens q}\tens b_k\tens1^{\tens t}) g_{q+1+t}
= \lambda_n \]
is the right hand side of \eqref{eq-gb-bg-Tnb}. This coderivation
commutes with $b$ since
$(\tilde{g}b-b\tilde{g})b+b(\tilde{g}b-b\tilde{g})=0$. Applying this
identity to $T^ns\cb$ and composing it with $\pr_1:Ts\cc\to s\cc,$ we
get an identity
\[ (\tilde{g}b-b\tilde{g})_nb_1 + \sum_{q+1+t=n}
(1^{\tens q}\tens b_1\tens1^{\tens t})(\tilde{g}b-b\tilde{g})_n = 0, \]
which means precisely that $\lambda_nd=0$.

Introduce a $(f,\tilde{g})$-coderivation $\tilde{r}:Ts\cb\to Ts\cc$ by
its components
$(r_0,r_1,\dots,r_{n-1},0$, $0$, $\dots)$ (these are already
known). The commutator $\tilde{r}b+b\tilde{r}$ has the following
property:
\[ (\tilde{r}b+b\tilde{r})\Delta = \Delta\bigl[
f\tens(\tilde{r}b+b\tilde{r}) + (\tilde{r}b+b\tilde{r})\tens\tilde{g}
+ \tilde{r}\tens(\tilde{g}b-b\tilde{g})\bigr]. \]
Let us construct a map
$\theta=[\tilde{r}\tens(\tilde{g}b-b\tilde{g})]\theta:Ts\cb\to Ts\cc$
for the data $f \rTTo^{\tilde{r}} \tilde{g}
\rTTo^{\tilde{g}b-b\tilde{g}} \tilde{g}:Ts\cb\to Ts\cc$
as in Section~3 of \cite{Lyu-AinfCat} (see also
\secref{sec-convent-nota}). Its components
$\theta_{kl}=\theta\big|_{T^ks\cb}\pr_l:T^ks\cb\to T^ls\cc$ are given
by formula~\eqref{eq-theta-kl-short}
\begin{equation*}
\theta_{kl} =
\sum_{\substack{\alpha+\beta+\gamma+2=l\\a+j+c+t+e=k}}
f_{a\alpha}\tens\tilde{r}_j\tens\tilde{g}_{c\beta}\tens
(\tilde{g}b-b\tilde{g})_t\tens\tilde{g}_{e\gamma}.
\end{equation*}
By Proposition~3.1 of \cite{Lyu-AinfCat} the map $\theta$ satisfies the
equation
\[ \theta\Delta = \Delta\bigl[f\tens\theta + \theta\tens\tilde{g} +
\tilde{r}\tens(\tilde{g}b-b\tilde{g})\bigr]. \]
Therefore,
$\nu=-\tilde{r}b-b\tilde{r}+
[\tilde{r}\tens(\tilde{g}b-b\tilde{g})]\theta:Ts\cb\to Ts\cc$
is a $(f,\tilde{g})$-coderivation. Since $\theta_{k1}=0$ for all $k$,
the components $\nu_k$ vanish for $k<n$, and
\begin{multline*}
\nu_n = -
\sum_{\substack{k<n;(q,k,t)\ne(0,0,1)\\i_1+\dots+i_q+k+j_1+\dots+j_t=n}}
(f_{i_1}\tens\dots\tens f_{i_q}\tens r_k\tens
 g_{j_1}\tens\dots\tens g_{j_t})b_{q+1+t} \\
- \sum_{k>1;q+k+t=n}(1^{\tens q}\tens b_k\tens1^{\tens t})r_{q+1+t}
\end{multline*}
which is the right hand side of \eqref{eq-rb+br-TnB}. We have an
obvious identity
\[ \nu b - b\nu = (-\tilde{r}b-b\tilde{r}+\theta)b
-b(-\tilde{r}b-b\tilde{r}+\theta) = \theta b - b\theta. \]
Applying this identity to $T^ns\cb$ and composing it with
$\pr_1:Ts\cc\to s\cc,$ we get an identity
\[ \nu_nb_1 - \sum_{q+1+t=n}(1^{\tens q}\tens b_1\tens1^{\tens t})\nu_n
= \bigl[r_0\tens(\tilde{g}b-b\tilde{g})_n\bigr]b_2, \]
since $[\tilde{r}\tens(\tilde{g}b-b\tilde{g})]\theta_{nl}$ vanishes for
$l\ne2$, and equals $r_0\tens(\tilde{g}b-b\tilde{g})_n$ for $l=2$. The
above equation means precisely that $\nu_nd=-\lambda_nu$. Indeed,
$(\tilde{g}b-b\tilde{g})_n(r_{X_0}\tens1)
=-r_0\tens(\tilde{g}b-b\tilde{g})_n$.
Thus, proposition is proved by induction.
\end{proof}

\subsection{\texorpdfstring
{Transformations between the constructed $A_\infty$-functors.}
{Transformations between the constructed A8-functors.}}
    \label{sec-Transformations}
Let $\cb$, $\cc$ be \ainf-categories, let $f:\cb\to\cc$ be an
\ainf-functor, let $g:\Ob\cb\to\Ob\cc$ be a map, and assume that for
each object $X\in\Ob\cb$ there is a map
$r_0:\kk\to(s\cc)^{-1}(Xf,Xg)$ such that $r_0b_1=0$. Let the
assumptions of \secref{sec-inj-res-ainf-fun} hold. Let $g,g':\cb\to\cc$
be two \ainf-functors, whose underlying map is the given
$g:\Ob\cb\to\Ob\cc$. Let $r:f\to g:\cb\to\cc$, $r':f\to g':\cb\to\cc$
be natural \ainf-transformations, whose 0\n-th component $r_0=r'_0$ is
the given map $r_0:\kk\to(s\cc)^{-1}(Xf,Xg)$.

\begin{proposition}\label{pro-p-g-g-exists}
Under the above assumptions, there exists a natural \ainf-transformation
$p:g\to g':\cb\to\cc$ such that $r'=(f \rTTo^r g \rTTo^p g')$ in the
2\n-category $A_\infty$.
\end{proposition}

\begin{proof}
Let us construct a $(g,g')$\n-coderivation $p:Ts\cb\to Ts\cc$ of degree
$-1$ and a $(f,g')$\n-coderivation $v:Ts\cb\to Ts\cc$ of degree $-2$
such that
\begin{align*}
pb+bp &= 0, \\
(r\tens p)B_2 - r' &= [v,b],
\end{align*}
that is, $p:g\to g':\cb\to\cc$ is a 2\n-morphism and
$v:(r\tens p)B_2\to r':f\to g':\cb\to\cc$ is a 3\n-morphism. Let us
build the components of $p$ and $v$ by induction. We have $p_k=0$ and
$v_k=0$ for $k<0$. Given non-negative $n$, assume that we have already
found components $p_m$, $v_m$ of the sought for $p$, $v$ for $m<n$, such
that equations
\begin{gather*}
(pb+bp)\pr_1=0:s\cb(X_0,X_1)\tens\dots\tens s\cb(X_{m-1},X_m)
\to s\cc(X_0g,X_mg), \\
\{(r\tens p)B_2-r'-[v,b]\}\pr_1=0:
s\cb(X_0,X_1)\tens\dots\tens s\cb(X_{m-1},X_m)\to s\cc(X_0f,X_mg)
\end{gather*}
are satisfied for all $m<n$. Under these assumptions we will find such
$p_n$, $v_n$ that the above equations are satisfied for $m=n$. Notice
that for $m=n=0$ the source complexes reduce to $\kk$. Let us write
down these equations explicitly. The terms which contain unknown maps
$p_n$, $v_n$ are singled out on the left hand side. The right hand side
consists of already known terms:
\begin{multline}
- p_nb_1 - \sum_{q+1+t=n}(1^{\tens q}\tens b_1\tens1^{\tens t})p_n \\
= \sum_{q+k+t=n}^{k>1}(1^{\tens q}\tens b_k\tens1^{\tens t})p_{q+1+t}
+ \sum_{i_1+\dots+i_q+k+j_1+\dots+j_t=n}^{k<n}
(g_{i_1}\tens\dots\tens g_{i_q}\tens p_k\tens
g'_{j_1}\tens\dots\tens g'_{j_t})b_{q+1+t},
\label{eq-pb-bp-Tnb}
\end{multline}
\begin{multline}
v_nb_1 - \sum_{q+1+t=n}(1^{\tens q}\tens b_1\tens1^{\tens t})v_n
- (r_0\tens p_n)b_2 \\
=\sum^{k<n}_{a_1+\dots+a_\alpha+j+c_1+\dots+c_\beta+k+e_1+\dots+e_\gamma=n}
\hspace*{-13mm} (f_{a_1}\tens\dots\tens f_{a_\alpha}\tens
r_j\tens g_{c_1}\tens\dots\tens g_{c_\beta}\tens p_k\tens
g'_{e_1}\tens\dots\tens g'_{e_\gamma})b_{\alpha+\beta+\gamma+2}
- r'_n \\
- \sum^{k<n}_{i_1+\dots+i_q+k+j_1+\dots+j_t=n}
(f_{i_1}\tens\dots\tens f_{i_q}\tens v_k\tens
g'_{j_1}\tens\dots\tens g'_{j_t})b_{q+1+t}
+ \sum^{k>1}_{q+k+t=n}(1^{\tens q}\tens b_k\tens1^{\tens t})v_{q+1+t}.
\label{eq-vb+bv-TnB}
\end{multline}
The components of $(r\tens p)B_2$ are computed by formula~(5.1.3) of
\cite{Lyu-AinfCat}:
\[ [(r\tens p)B_2]_n = \sum_l (r\tens p)\theta_{nl}b_l.
\]
Denote by $\lambda_n\in\Hom_\kk^0(N,s\cc(X_0g,X_ng))$ the right hand
side of \eqref{eq-pb-bp-Tnb} and by
$\nu_n\in\Hom_\kk^{-1}(N,s\cc(X_0f,X_ng))$ the right hand side of
\eqref{eq-vb+bv-TnB}, where
$N=s\cb(X_0,X_1)\tens_\kk\dots\tens_\kk s\cb(X_{n-1},X_n)$. In
particular, $N=\kk$ for $n=0$. Equations \eqref{eq-pb-bp-Tnb} and
\eqref{eq-vb+bv-TnB} mean that
$(v_n,p_n)d=(\nu_n,\lambda_n)\in\Cone^{-1}(u)$. Since $\Cone(u)$ is
acyclic, such a pair $(v_n,p_n)\in\Cone^{-2}(u)$ exists if and only if
$(\nu_n,\lambda_n)\in\Cone^{-1}(u)$ is a cycle, that is, equations
$-\lambda_nd=0$, $\nu_nd+\lambda_nu=0$ are satisfied. Let us verify
them now.

Introduce a $(g,g')$-coderivation $\tilde{p}:Ts\cb\to Ts\cc$ of degree
$-1$ by its components $(p_0,p_1,\dots,p_{n-1},0,0$, $\dots)$. The
commutator $\lambda=\tilde{p}b+b\tilde{p}$ is also a
$(g,g')$-coderivation (of degree 0). Its components $\lambda_m$ vanish
for $m<n$. The component $\lambda_n=(\tilde{p}b+b\tilde{p})_n$ is the
right hand side of \eqref{eq-pb-bp-Tnb}. Consider the identity
\[ (\tilde{p}b+b\tilde{p})b-b(\tilde{p}b+b\tilde{p})=0. \]
Applying this identity to $T^ns\cb$ and composing it with
$\pr_1:Ts\cc\to s\cc,$ we get an identity
\[ \lambda_nb_1 -
\sum_{q+1+t=n}(1^{\tens q}\tens b_1\tens1^{\tens t})\lambda_n = 0, \]
that is, $\lambda_nd=0$.

Introduce a $(f,g')$-coderivation $\tilde{v}:Ts\cb\to Ts\cc$ of degree
$-2$ by its components $(v_0,v_1,\dots,v_{n-1},0,0$, $\dots)$. All
summands of the map $\nu=(r\tens\tilde{p})B_2-r'-[\tilde{v},b]$ are
$(f,g')$-coderivations of degree $-1$. Hence, the same holds for $\nu$.
The components $\nu_m$ vanish for $m<n$. The component $\nu_n$ is the
right hand side of \eqref{eq-vb+bv-TnB}. Consider the commutator
\begin{multline*}
[\nu,b] = \nu B_1
= (r\tens\tilde{p})B_2B_1 - r'B_1 - \tilde{v}B_1B_1
= - (r\tens\tilde{p})(1\tens B_1+B_1\tens1)B_2 \\
= - (r\tens\tilde{p}B_1)B_2 = - (r\tens\lambda)B_2.
\end{multline*}
Applying this identity to $T^ns\cb$ and composing it with
$\pr_1:Ts\cc\to s\cc$ we get an identity
\[ \nu_nb_1 + \sum_{q+1+t=n}(1^{\tens q}\tens b_1\tens1^{\tens t})\nu_n
= - (r_0\tens\lambda_n)b_2, \]
that is, $\nu_nd=-\lambda_nu$. Thus, the proposition is proved by
induction.
\end{proof}

\begin{proposition}[Uniqueness of the transformations]
\label{prop-p-unique}
Let assumptions of Sections \ref{sec-inj-res-ainf-fun} and
\ref{sec-Transformations} hold. The natural \ainf-transformation
$p:g\to g':\cb\to\cc$, such that $r'=(f \rTTo^r g \rTTo^p g')$ in the
2\n-category $A_\infty$, is unique up to an equivalence.
\end{proposition}

\begin{proof}
Assume that we have two such 2\n-morphisms $p,q:g\to g':\cb\to\cc$ and
two 3\n-morphisms $v:(r\tens p)B_2\to r':f\to g':\cb\to\cc$ and
$w:(r\tens q)B_2\to r':f\to g':\cb\to\cc$. We are looking for a
3\n-morphism $x:p\to q:g\to g':\cb\to\cc$ and the following
4\n-morphism, whose source depends on $x$. Assuming that $p-q=xB_1$ we
deduce that
\[ -(r\tens x)B_2B_1 = (r\tens xB_1)B_2
= (r\tens p)B_2 - (r\tens q)B_2. \]
Since $(r\tens q)B_2-r'=wB_1$, we find out that
$(r\tens p)B_2-r'=[w-(r\tens x)B_2]B_1$. Thus, we have two
3\n-morphisms with the common source and target
$v,w-(r\tens x)B_2:(r\tens p)B_2\to r':f\to g':\cb\to\cc$. We are
looking for a 4\n-morphism
\[ z:w-(r\tens x)B_2\to v:(r\tens p)B_2\to r':f\to g':\cb\to\cc, \]
as well as for $x$. In other terms, we have to find coderivations $x$
of degree $-2$ and $z$ of degree $-3$ such that the following equations
hold:
\begin{align*}
p-q &= xb-bx, \\
w - (r\tens x)B_2 - v &= zb + bz.
\end{align*}
Let us build the components of $x$ and $z$ by induction. We have
$x_k=0$ and $z_k=0$ for $k<0$. Given non-negative $n$, assume that we
have already found components $x_m$, $z_m$ of the sought $x$, $z$ for
$m<n$, such that equations
\[ (p-q)\pr_1=(xb-bx)\pr_1:s\cb(X_0,X_1)\tens\dots\tens s\cb(X_{m-1},X_m)
\to s\cc(X_0g,X_mg), \]
\begin{multline*}
[w-(r\tens x)B_2-v]\pr_1=(zb + bz)\pr_1: \\
s\cb(X_0,X_1)\tens\dots\tens s\cb(X_{m-1},X_m)\to s\cc(X_0f,X_mg)
\end{multline*}
are satisfied for all $m<n$. Under these assumptions, we will find such
$x_n$, $z_n$ that the above equations are satisfied for $m=n$. Notice
that for $m=n=0$ the source complexes reduce to $\kk$. Let us write
down these equations explicitly. The terms which contain unknown maps
$x_n$, $z_n$ are singled out on the left hand side. The right hand side
consists of already known terms:
\begin{multline}
- x_nb_1 + \sum_{\alpha+1+\beta=n}
(1^{\tens\alpha}\tens b_1\tens1^{\tens\beta})x_n = q_n - p_n \\
+ \sum_{i_1+\dots+i_\alpha+k+j_1+\dots+j_\beta=n}^{k<n}
\hspace*{-6mm} (g_{i_1}\tens\dots\tens g_{i_\alpha}\tens x_k\tens
g'_{j_1}\tens\dots\tens g'_{j_\beta})b_{\alpha+1+\beta}
- \sum_{\alpha+k+\beta=n}^{k>1}
(1^{\tens\alpha}\tens b_k\tens1^{\tens\beta})x_{\alpha+1+\beta},
\label{eq-xb-bx-Tnb}
\end{multline}
\begin{multline}
z_nb_1 + \sum_{\alpha+1+\beta=n}
(1^{\tens\alpha}\tens b_1\tens1^{\tens\beta})z_n + (r_0\tens x_n)b_2 \\
=-\sum^{k<n}_{a_1+\dots+a_\alpha+j+c_1+\dots+c_\beta+k+e_1+\dots+e_\gamma=n}
\hspace*{-9mm} (f_{a_1}\tens\dots\tens f_{a_\alpha}\tens
r_j\tens g_{c_1}\tens\dots\tens g_{c_\beta}\tens x_k\tens
g'_{e_1}\tens\dots\tens g'_{e_\gamma})b_{\alpha+\beta+\gamma+2} \\
+ w_n - v_n - \sum^{k<n}_{i_1+\dots+i_\alpha+k+j_1+\dots+j_\beta=n}
(f_{i_1}\tens\dots\tens f_{i_\alpha}\tens z_k\tens
g'_{j_1}\tens\dots\tens g'_{j_\beta})b_{\alpha+1+\beta} \\
- \sum^{k>1}_{\alpha+k+\beta=n}
(1^{\tens\alpha}\tens b_k\tens1^{\tens\beta})z_{\alpha+1+\beta}.
\label{eq-zb+bz-TnB}
\end{multline}
Denote by $\lambda_n\in\Hom_\kk^{-1}(N,s\cc(X_0g,X_ng))$ the right hand
side of \eqref{eq-xb-bx-Tnb} and by
$\nu_n\in\Hom_\kk^{-2}(N,s\cc(X_0f,X_ng))$ the right hand side of
\eqref{eq-zb+bz-TnB}, where
$N=s\cb(X_0,X_1)\tens_\kk\dots\tens_\kk s\cb(X_{n-1},X_n)$. In
particular, $N=\kk$ for $n=0$. Equations \eqref{eq-xb-bx-Tnb} and
\eqref{eq-zb+bz-TnB} mean that
$(z_n,x_n)d=(\nu_n,\lambda_n)\in\Cone^{-2}(u)$. Since $\Cone(u)$ is
acyclic, such a pair $(z_n,x_n)\in\Cone^{-3}(u)$ exists if and only if
$(\nu_n,\lambda_n)\in\Cone^{-2}(u)$ is a cycle, that is, equations
$-\lambda_nd=0$, $\nu_nd+\lambda_nu=0$ are satisfied. Let us verify
them now.

Introduce a $(g,g')$-coderivation $\tilde{x}:Ts\cb\to Ts\cc$ of degree
$-2$ by its components $(x_0,x_1,\dots,x_{n-1},0,0$, $\dots)$. The
commutator $\tilde{x}b-b\tilde{x}$ is also a $(g,g')$-coderivation (of
degree $-1$). Hence, the map $\lambda=-p+q+\tilde{x}b-b\tilde{x}$ is
also a $(g,g')$-coderivation of degree $-1$. Its components $\lambda_m$
vanish for $m<n$. The component $\lambda_n$ is the right hand side of
\eqref{eq-xb-bx-Tnb}. Consider the identity
\[ \lambda B_1 = -pB_1 + qB_1 + \tilde{x}B_1B_1 = 0. \]
Applying this identity to $T^ns\cb$ and composing it with
$\pr_1:Ts\cc\to s\cc$ we get an identity
\[ \lambda_nb_1 + \sum_{\alpha+1+\beta=n}
(1^{\tens\alpha}\tens b_1\tens1^{\tens\beta})\lambda_n = 0, \]
that is, $\lambda_nd=0$.

Introduce a $(f,g')$-coderivation $\tilde{z}:Ts\cb\to Ts\cc$ of degree
$-3$ by its components $(z_0,z_1,\dots,z_{n-1},0,0$, $\dots)$. All
summands of the map $\nu=w-(r\tens\tilde{x})B_2-v-[\tilde{z},b]$ are
$(f,g')$-coderivations of degree $-2$. Hence, the same holds for $\nu$.
The components $\nu_m$ vanish for $m<n$. The component $\nu_n$ is the
right hand side of \eqref{eq-zb+bz-TnB}. Consider the commutator
\begin{multline*}
[\nu,b] = \nu B_1
= wB_1 - (r\tens\tilde{x})B_2B_1 - vB_1 - \tilde{z}B_1B_1 \\
= (r\tens q)B_2 - r' + (r\tens\tilde{x}B_1)B_2 - (r\tens p)B_2 + r'
= [r\tens(q-p+\tilde{x}B_1)]B_2 = (r\tens\lambda)B_2.
\end{multline*}
Applying this identity to $T^ns\cb$ and composing it with
$\pr_1:Ts\cc\to s\cc$ we get an identity
\[ \nu_nb_1 - \sum_{\alpha+1+\beta=n}
(1^{\tens\alpha}\tens b_1\tens1^{\tens\beta})\nu_n
= (r_0\tens\lambda_n)b_2, \]
that is, $\nu_nd=-\lambda_nu$. Thus, the proposition is proved by
induction.
\end{proof}

\begin{corollary}\label{cor-p-invertible}
In the assumptions of \propref{prop-p-unique}, let  $\cc$ be unital.
Then the constructed 2\n-morphism $p:g\to g':\cb\to\cc$ is
invertible in $A_\infty$.
\end{corollary}

\begin{proof}
Exchanging the pairs $(g,r)$ and $(g',r')$, we see that there is a
2\n-morphism $t:g'\to g:\cb\to\cc$, such that
$r=(f \rTTo^{r'} g' \rTTo^t g)$. Therefore,
$r=(f \rTTo^r g \rTTo^{p\cdot t} g)$. Since $\cc$ is unital, there
is a unit 2\n-endomorphism $1_gs=g\uni^\cc:g\to g:\cb\to\cc$. It
satisfies the equation $r=(f \rTTo^r g \rTTo^{1_gs} g)$. The
uniqueness proved in \propref{prop-p-unique} implies that
$p\cdot t=1_gs$. Similarly, $t\cdot p=1_{g'}s$.
\end{proof}

\begin{proposition}[Unitality of $A_\infty$-functors]
\label{prop-g-unital}
Let the assumptions of \secref{sec-inj-res-ainf-fun} hold. If
\ainf-categories $\cb$, $\cc$ are unital and \ainf-functor
$f:\cb\to\cc$ is unital, then the \ainf-functor $g:\cb\to\cc$
constructed in \propref{prop-ainf-fr-g-tr-r} is unital as well.
\end{proposition}

\begin{proof}
We are given a 2\n-morphism $r:f\to g:\cb\to\cc$ and a 3\n-morphism
$v:f\uni^\cc\to \uni^\cb f:f\to f:\cb\to\cc$. We are looking for a
3\n-morphism $w:g\uni^\cc\to \uni^\cb g:g\to g:\cb\to\cc$ and a
4\n-morphism $x$, whose target depends on $w$. Let us describe $x$ now
for the above $w$.

We have the following 3\n-morphisms
\begin{alignat*}3
(r\tens \uni^\cc)M_{11} &:& (f\uni^\cc\tens r)B_2 &\to&
(r\tens g\uni^\cc)B_2 &:f\to g:\cb\to\cc, \\
(v\tens r)B_2 &:& (f\uni^\cc\tens r)B_2 &\to& (\uni^\cb f\tens r)B_2
&:f\to g:\cb\to\cc, \\
(r\tens w)B_2 &:& (r\tens \uni^\cb g)B_2 &\to& (r\tens g\uni^\cc)B_2
&:f\to g:\cb\to\cc, \\
(\uni^\cb\tens r)M_{11} &:& (r\tens \uni^\cb g)B_2 &\to&
(\uni^\cb f\tens r)B_2 &:f\to g:\cb\to\cc.
\end{alignat*}
Indeed, for the cocategory homomorphism
$M:TsA_\infty\boxtimes TsA_\infty\to TsA_\infty$ we have an equation
$MB=(1\boxtimes B+B\boxtimes1)M$, see Section~6 of \cite{Lyu-AinfCat}.
It implies, in particular, that
\begin{align*}
(r\tens \uni^\cc)M_{11}B_1 - (f\uni^\cc\tens r)B_2 + (r\tens g\uni^\cc)B_2
&= (r\tens \uni^\cc)(1\tens B_1+B_1\tens1)M_{11} = 0, \\
(v\tens r)B_2B_1 = (vB_1\tens r)B_2
&= (f\uni^\cc\tens r)B_2 - (\uni^\cb f\tens r)B_2, \\
(r\tens w)B_2B_1 = -(r\tens wB_1)B_2
&= (r\tens \uni^\cb g)B_2 - (r\tens g\uni^\cc)B_2, \\
(\uni^\cb\tens r)M_{11}B_1 -(r\tens \uni^\cb g)B_2 +(\uni^\cb f\tens r)B_2
&= (\uni^\cb\tens r)(1\tens B_1+B_1\tens1)M_{11} = 0.
\end{align*}

Linear combinations of the above maps form 3\n-morphisms with the same
source and target
\begin{align*}
(r\tens \uni^\cc)M_{11}-(v\tens r)B_2 &:
(\uni^\cb f\tens r)B_2\to(r\tens g\uni^\cc)B_2:f\to g:\cb\to\cc, \\
(r\tens w)B_2-(\uni^\cb\tens r)M_{11} &:
(\uni^\cb f\tens r)B_2\to(r\tens g\uni^\cc)B_2:f\to g:\cb\to\cc.
\end{align*}
We are looking for a 4\n-morphism between the above 3\n-morphisms
\begin{multline*}
x:(r\tens \uni^\cc)M_{11}-(v\tens r)B_2
\to (r\tens w)B_2-(\uni^\cb\tens r)M_{11} \\
:(\uni^\cb f\tens r)B_2\to(r\tens g\uni^\cc)B_2:f\to g:\cb\to\cc,
\end{multline*}
as well as for $w$.

In other words, we have to find a $(g,g)$-coderivation $w$ of degree
$-2$ and an $(f,g)$-coderivation $x$ of degree $-3$ such that the
following equations hold:
\begin{align*}
-wb+bw &= \uni^\cb g - g\uni^\cc, \\
xb+bx &=
(r\tens \uni^\cc)M_{11}-(v\tens r)B_2-(r\tens w)B_2+(\uni^\cb\tens r)M_{11}.
\end{align*}
Let us construct the components of $w$ and $x$ by induction. We have
$w_k=0$ and $x_k=0$ for $k<0$. Given non-negative $n$, assume that we
have already found components $w_m$, $x_m$ of the sought for
$x$, $z$ for
$m<n$, such that the above equations restricted to $T^ms\cb$ are
satisfied for all $m<n$. Under these assumptions, we will find such $w_n$,
$x_n$  for $m=n$. Let us write down these
equations explicitly. The terms which contain unknown maps $w_n$, $x_n$
are singled out on the left hand side. The right hand side consists of
already known terms:
\begin{multline}
- w_nb_1 + \sum_{q+1+t=n}(1^{\tens q}\tens b_1\tens1^{\tens t})w_n
= \sum_{i_1+\dots+i_q+k+j_1+\dots+j_t=n}^{k<n} \hspace*{-9mm}
(g_{i_1}\tens\dots\tens g_{i_q}\tens w_k\tens
g_{j_1}\tens\dots\tens g_{j_t})b_{q+1+t} \\
- \sum_{q+k+t=n}^{k>1}(1^{\tens q}\tens b_k\tens1^{\tens t})w_{q+1+t}
+ \sum_{q+k+t=n}(1^{\tens q}\tens\uni^\cb_k\tens1^{\tens t})g_{q+1}
- \uni^\cc_n,
\label{eq-wb+bw-Tnb-unit}
\end{multline}
\begin{multline}
x_nb_1 + \sum_{q+1+t=n}(1^{\tens q}\tens b_1\tens1^{\tens t})x_n
+ (r_0\tens w_n)b_2 \\
= - \sum^{k<n}_{i_1+\dots+i_q+k+j_1+\dots+j_t=n} \hspace*{-1mm}
(f_{i_1}\tens\dots\tens f_{i_q}\tens x_k\tens
g_{j_1}\tens\dots\tens g_{j_t})b_{q+1+t}
- \sum^{k>1}_{q+k+t=n}(1^{\tens q}\tens b_k\tens1^{\tens t})x_{q+1+t} \\
-\sum^{k<n}_{a_1+\dots+a_\alpha+j+c_1+\dots+c_\beta+k+e_1+\dots+e_\gamma=n}
\hspace*{-2mm} (f_{a_1}\tens\dots\tens f_{a_\alpha}\tens
v_j\tens f_{c_1}\tens\dots\tens f_{c_\beta}\tens r_k\tens
g_{e_1}\tens\dots\tens g_{e_\gamma})b_{\alpha+\beta+\gamma+2} \\
-\sum^{k<n}_{a_1+\dots+a_\alpha+j+c_1+\dots+c_\beta+k+e_1+\dots+e_\gamma=n}
\hspace*{-3mm} (f_{a_1}\tens\dots\tens f_{a_\alpha}\tens
r_j\tens g_{c_1}\tens\dots\tens g_{c_\beta}\tens w_k\tens
g_{e_1}\tens\dots\tens g_{e_\gamma})b_{\alpha+\beta+\gamma+2} \\
+ \sum_{i_1+\dots+i_q+k+j_1+\dots+j_t=n} \hspace*{-3mm}
(f_{i_1}\tens\dots\tens f_{i_q}\tens r_k\tens
g_{j_1}\tens\dots\tens g_{j_t})\uni^\cc_{q+1+t}
+ \sum_{q+k+t=n}(1^{\tens q}\tens\uni^\cb_k\tens1^{\tens t})r_{q+1+t}.
\label{eq-xb+bx-TnB-unit}
\end{multline}
Denote by $\lambda_n\in\Hom_\kk^{-1}(N,s\cc(X_0g,X_ng))$ the right hand
side of \eqref{eq-wb+bw-Tnb-unit} and by
$\nu_n\in\Hom_\kk^{-2}(N,s\cc(X_0f,X_ng))$ the right hand side of
\eqref{eq-xb+bx-TnB-unit}, where
$N=s\cb(X_0,X_1)\tens_\kk\dots\tens_\kk s\cb(X_{n-1},X_n)$. Equations
\eqref{eq-wb+bw-Tnb-unit} and \eqref{eq-xb+bx-TnB-unit} mean that
$(x_n,w_n)d=(\nu_n,\lambda_n)\in\Cone^{-2}(u)$. Since $\Cone(u)$ is
acyclic, such a pair $(x_n,w_n)\in\Cone^{-3}(u)$ exists if and only if
$(\nu_n,\lambda_n)\in\Cone^{-2}(u)$ is a cycle, that is, equations
$-\lambda_nd=0$, $\nu_nd+\lambda_nu=0$ are satisfied. Let us verify
them now.

Introduce a $(g,g)$-coderivation $\tilde{w}:Ts\cb\to Ts\cc$ of degree
$-2$ by its components $(w_0,w_1,\dots,w_{n-1},0,0$, $\dots)$. Hence,
the map $\lambda=\tilde{w}b-b\tilde{w}+\uni^\cb g-g\uni^\cc$ is also a
$(g,g)$-coderivation of degree $-1$. Its components $\lambda_m$ vanish
for $m<n$. The component $\lambda_n$ is the right hand side of
\eqref{eq-wb+bw-Tnb-unit}. Consider the identity
\[ \lambda B_1 = \tilde{w}B_1B_1+(\uni^\cb g)B_1-(g\uni^\cc)B_1 =0. \]
Applying this identity to $T^ns\cb$ and composing it with
$\pr_1:Ts\cc\to s\cc,$ we get an identity
\[ \lambda_nb_1 +
\sum_{q+1+t=n}(1^{\tens q}\tens b_1\tens1^{\tens t})\lambda_n = 0, \]
that is, $\lambda_nd=0$.

Introduce a $(f,g)$-coderivation $\tilde{x}:Ts\cb\to Ts\cc$ of degree
$-3$ by its components $(x_0,x_1,\dots,x_{n-1},0,0$, $\dots)$. All
summands of the map
\[ \nu = - \tilde{x}B_1 + (r\tens \uni^\cc)M_{11} - (v\tens r)B_2
- (r\tens\tilde{w})B_2+(\uni^\cb\tens r)M_{11}  \]
are $(f,g)$-coderivations of degree $-2$. Hence, the same holds for
$\nu$. The components $\nu_m$ vanish for $m<n$. The component $\nu_n$
is the right hand side of \eqref{eq-xb+bx-TnB-unit}. Consider its
differential
\begin{multline*}
\nu B_1
= - \tilde{x}B_1B_1 + (r\tens \uni^\cc)M_{11}B_1 - (v\tens r)B_2B_1
- (r\tens\tilde{w})B_2B_1 + (\uni^\cb\tens r)M_{11}B_1  \\
= (f\uni^\cc\tens r)B_2 - (r\tens g\uni^\cc)B_2 - (vB_1\tens r)B_2
+(r\tens\tilde{w}B_1)B_2+(r\tens \uni^\cb g)B_2-(\uni^\cb f\tens r)B_2 \\
= [r\tens(\tilde{w}B_1 -g\uni^\cc +\uni^\cb g)]B_2 = (r\tens\lambda)B_2.
\end{multline*}
Applying identity $\nu B_1=(r\tens\lambda)B_2$ to $T^ns\cb$ and
composing it with $\pr_1:Ts\cc\to s\cc$ we get an identity
\[ \nu_nb_1 - \sum_{q+1+t=n}(1^{\tens q}\tens b_1\tens1^{\tens t})\nu_n
= (r_0\tens\lambda_n)b_2, \]
that is, $\nu_nd=-\lambda_nu$. Thus, the proposition is proved by
induction.
\end{proof}

\subsection{Invertible transformations.}\label{sec-Invertible-trans}
Let $\cb$, $\cc$ be unital \ainf-categories, and let
$f,g:\Ob\cb\to\Ob\cc$ be maps. Assume that for each object $X$ of $\cb$
there are $\kk$\n-linear maps
\begin{alignat*}3
\sS{_X}r_0 &: \kk \to (s\cc)^{-1}(Xf,Xg), &\qquad
\sS{_X}p_0 &: \kk \to (s\cc)^{-1}(Xg,Xf), \\
\sS{_X}w_0 &: \kk \to (s\cc)^{-2}(Xf,Xf), &\qquad
\sS{_X}v_0 &: \kk \to (s\cc)^{-2}(Xg,Xg),
\end{alignat*}
such that
\begin{gather}
\sS{_X}r_0b_1=0, \qquad \sS{_X}p_0b_1=0, \notag \\
(\sS{_X}r_0\tens\sS{_X}p_0)b_2-\sS{_{Xf}}\uni^\cc_0 = \sS{_X}w_0b_1,
\label{eq-rp-i-pr-i} \\
(\sS{_X}p_0\tens\sS{_X}r_0)b_2 - \sS{_{Xg}}\uni^\cc_0 = \sS{_X}v_0b_1.
\notag
\end{gather}

\begin{proposition}\label{pro-Invertible-rpfgBC}
Let the assumptions of \secref{sec-Invertible-trans} hold and, moreover,
let $f:\cb\to\cc$ be a unital \ainf-functor. Then the map $g$ extends
to a unital \ainf-functor $g:\cb\to\cc$ and the given $r_0$, $p_0$ extend
to natural \ainf-transformations $r:f\to g:\cb\to\cc$,
$p:g\to f:\cb\to\cc$, inverse to each other.
\end{proposition}

\begin{proof}
Propositions \ref{prop-ainf-fr-g-tr-r} and \ref{prop-g-unital} imply
the existence and unitality of $g$. Indeed, since $(r_0\tens1)b_2$ is a
homotopy invertible chain map, the map $u=\Hom(N,(r_0\tens1)b_2)$ is
also homotopy invertible, hence a quasi-isomorphism. Existence of
$r:f\to g:\cb\to\cc$ is shown in \propref{prop-ainf-fr-g-tr-r}.
Existence of $p:g\to f:\cb\to\cc$, inverse to $r$ is proven in
\cite[Proposition~7.15]{Lyu-AinfCat}.
\end{proof}

\section{Derived categories}
Let $\ca$ be a $\fu$\n-small Abelian $\kk$\n-linear category, and let
$\cc=\Com(\ca)$ or $\cc=\Com^+(\ca)$ be the differential graded category of
complexes (resp. bounded below complexes) of objects of $\ca$. Denote
by $\cb=\Acyc(\ca)$ its full subcategory of acyclic complexes. Let
$\cd=\Dr(\cc|\cb)$ be the constructed $\fu$\n-small differential graded
category. We observe first that quasi-isomorphisms in $\cc$ become
(homotopy) invertible elements in $\cd$.

Assuming that the ground ring $\kk$ is a field, we turn to the procedure of
finding a K\n-injective resolution (if they exist) into a unital
\ainf-functor. Under these assumptions, we also show that
$\cd=\Dr(\cc|\cb)$ is \ainf-equivalent to $\ci\subset\cc$ -- the full
subcategory of K\n-injective complexes. Hence, $H^0(\cd)$ is equivalent
to the derived category $\Dr(\ca)$.

\subsection{Invertibility of quasi-isomorphisms.}\label{sec-invert-qis}
Assume that $X$, $Y$ are objects of $\cc$ and $q:X\to Y$ is a
quasi-isomorphism. In particular, $q\in\cc^0(X,Y)$, $qm_1=0$. Let us
prove that $r=qs\ju_1\in(s\cd)^{-1}(X,Y)$ is invertible in the sense
of \secref{sec-Invertible-trans}, that is, there are elements
$p\in(s\cd)^{-1}(Y,X)$, $w\in(s\cd)^{-2}(X,X)$, $v\in(s\cd)^{-2}(Y,Y)$,
such that
\begin{gather}
r\overline{b}_1=0, \qquad p\overline{b}_1=0, \notag \\
(r\tens p)\overline{b}_2 - 1_Xs = w\overline{b}_1,
\label{eq-rp-1-pr-1} \\
(p\tens r)\overline{b}_2 - 1_Ys = v\overline{b}_1. \notag
\end{gather}
Indeed, denote $C=\Cone(q)=(Y\oplus X[1],d^C)$, where
$(y,x)d^C=(yd^Y+xq,-xd^X)$ for $y\in Y^l$, $x\in X^{l+1}$. Since $q$ is
a quasi-isomorphism, $C$ is acyclic. There is a standard exact sequence
of complexes $0\to Y \rTTo^n C \rTTo^k X[1]\to0$ with the chain maps
$n$, $k$, $yn=(y,0)$, $(0,x)k=x$. From now on we denote by $n$, $k$
also the corresponding elements $n\in\cc^0(Y,C)$, $k\in\cc^1(C,X)$.
Define $p$ as
\( p = ns\tens ks \in (s\cc)^{-1}(Y,C)\tens(s\cc)^0(C,X) \subset
(s\cd)^{-1}(Y,X) \).
Then
\[ p\overline{b}_1 = pb = (ns\tens ks)(1\tens b_1+b_1\tens1+b_2)
= - (n\tens k)m_2s = -(nk)s = 0. \]
Denote by $h\in\cc^{-1}(X,C)$ the following $\kk$\n-linear embedding
$X\to C$, $X^l\to C^{l-1}=Y^{l-1}\oplus X^l$, $x\mapsto(0,x)$. Define
$w$ as
\( w = hs\tens ks \in (s\cc)^{-2}(X,C)\tens(s\cc)^0(C,X) \subset
(s\cd)^{-2}(X,X) \).
Then
\begin{align*}
w\overline{b}_1 &= wb = (hs\tens ks)(1\tens b_1+b_1\tens1+b_2)
= hm_1s\tens ks - (hk)s = (qn)s\tens ks - 1_Xs \\
&= (qs\tens ns)b_2\tens ks - 1_Xs
= [qs\tensor2(ns\tensor1ks)]\overline{b}_2 - 1_Xs
= (r\tens p)\overline{b}_2 - 1_Xs.
\end{align*}
Indeed, $qn=hm_1=hd+dh:X\to C$ as explicit computation shows:
\[ xqn = (xq,0) = (xq,-xd^X)+(0,xd^X) = (0,x)d^C+(0,xd) = x(hd+dh). \]
Denote by $z\in\cc^0(C,Y)$ the following $\kk$\n-linear projection
$z:C\to Y$, $(y,x)\mapsto y$. Define $v$ as
\( v = -ns\tens zs \in (s\cc)^{-1}(Y,C)\tens(s\cc)^{-1}(C,Y) \subset
(s\cd)^{-2}(Y,Y) \).
Then
\begin{align*}
v\overline{b}_1 &= vb = -(ns\tens zs)(1\tens b_1+b_1\tens1+b_2)
= -ns\tens zm_1s - (nz)s = ns\tens(kq)s - 1_Ys \\
&= ns\tens(ks\tens qs)b_2 - 1_Ys
= [(ns\tensor1ks)\tensor2qs]\overline{b}_2 - 1_Ys
= (p\tens r)\overline{b}_2 - 1_Ys.
\end{align*}
Indeed, $-kq=zm_1=zd-dz:C\to Y$ as explicit computation shows:
\[ -(y,x)kq = -xq = yd-yd-xq = (y,x)zd-(y,x)d^Cz = (y,x)(zd-dz). \]
Thus, equations~\eqref{eq-rp-1-pr-1} hold true.

\subsection{K-injective complexes.}\label{sec-K-injective-complexes}
A complex $A\in\Ob\cc$ is K\n-injective if and only if for every
quasi-isomorphism $t:X\to Y\in\cc$ the chain map
$\cc(t,A):\cc(Y,A)\to\cc(X,A)$ is a quasi-isomorphism
\cite[Proposition~1.5]{Spaltenstein88}. Assume that each complex
$X\in\cc$ has a right K\n-injective resolution $r_X:X\to Xi$, that is,
$r_X$ is a quasi-isomorphism and $Xi\in\Ob\cc$ is K\n-injective.
Moreover, if $X$ is K\n-injective, we assume that $Xi=X$ and $r_X=1_X$.
By definition, $\cc(r_X,A):s\cc(Xi,A)\to s\cc(X,A)$,
$fs\mapsto(r_Xf)s$ is a quasi-isomorphism. The assumption is
satisfied, when $\ca$ has enough injectives and $\cc=\Com^+(\ca)$, or
when $\ca=R\modul$,%
\footnote{If $R\in\sS{'}\fu\in\fu$, where $\sS{'}\fu$ is a smaller
universe, then $\ca=R\modul$ is a $\fu$\n-small $\sS{'}\fu$\n-category.
Is it possible to replace $\cb=\Acyc(\ca)$ with some
$\sS{'}\fu$\n-small category $\cb'\subset\cb$ to get a
$\sS{'}\fu$\n-category $\cd=\Dr(\cc|\cb')$ \ainf-equivalent to
$\fu$\n-small $\cd=\Dr(\cc|\cb)$?}
 or when $\co$ is a sheaf of rings on a topological space, and $\ca$ is
the category of sheaves of left $\co$\n-modules, see
\cite{Spaltenstein88}.

Assume now that $\kk$ is a field. Then for any chain complex of
$\kk$\n-modules of the form
$N=s\cc(X_0,X_1)\tens_\kk s\cc(X_1,X_2)\tens_\kk\dots\tens_\kk
s\cc(X_{n-1},X_n)$,
$n\ge0$, $X_i\in\Ob\cc$, for any quasi-isomorphism $r_X:X\to Y$ and
for any K\n-injective $A\in\cc,$ the following chain map
\begin{equation*}
u = \Hom(N,\cc(r_X,A)): \Hom_\kk^\bull(N,s\cc(Y,A)) \to
\Hom_\kk^\bull(N,s\cc(X,A)),
\end{equation*}
is a quasi-isomorphism (any $\kk$\n-module complex is K\n-projective).
Therefore, we may apply the results of
\secref{sec-constr-inj-res-ainf-fun} to the differential graded
category $\cc=\Com$ or $\Com^+$, and its full subcategories
$\cb=\Acyc(\ca)$ (resp. $\ci=\Inj(\ca)$, $\cj=\AcIn(\ca)$) of acyclic
(resp. K\n-injective, acyclic K\n-injective) complexes. Denote by
$e:\ci \rMono \cc$ the full embedding. Starting with the identity
functor $f=\id_\cc,$ we get the existence of $g=ie$ simultaneously with
the existence of a unital \ainf-functor $i:\cc\to\ci$ --
``K\n-injective resolution functor'' and a natural
\ainf-transformation $r:\id\to ie:\cc\to\cc$ (Propositions
\ref{prop-ainf-fr-g-tr-r}, \ref{prop-g-unital}). The said $i$, $r$ are
unique in the sense of Propositions \ref{pro-p-g-g-exists},
\ref{prop-p-unique} and \corref{cor-p-invertible}. Moreover, while
solving equations \eqref{eq-gb-bg-Tnb}--\eqref{eq-rb+br-TnB} we will
choose the solutions
\begin{gather*}
i_n =g_n =\id_n: s\cc(X_0,X_1)\tens_\kk\dots\tens_\kk s\cc(X_{n-1},X_n)
\to s\cc(X_0,X_n), \\
r_n =\uni^\cc_n: s\cc(X_0,X_1)\tens_\kk\dots\tens_\kk s\cc(X_{n-1},X_n)
\to s\cc(X_0,X_n),
\end{gather*}
if $X_0$, \dots, $X_n$ are K\n-injective (recall that $X_0i=X_0$,
$X_ni=X_n$).

Extending $e$, $i$ to \ainf-functors between the constructed
categories, we get a unital strict \ainf-embedding (actually, a
faithful differential graded functor)
$\overline{e}:\Dr(\ci|\cj)\to\Dr(\cc|\cb)$, which is injective on
objects, and a unital \ainf-functor $\iu:\Dr(\cc|\cb)\to\Dr(\ci|\cj)$.
Let us prove that these \ainf-functors are quasi-inverse to each other.
First of all, $ei=\id_\ci$ implies
$\overline{e}\,\iu=\id_{\Dr(\ci|\cj)}$. Secondly, there is a natural
\ainf-transformation
$\overline{r}:\id\to\iu\,\overline{e}:\Dr(\cc|\cb)\to\Dr(\cc|\cb)$. Let
us prove that it is invertible.

The 0\n-th component is
\[ \sS{_X}{\overline{r}}_0 = \bigl[ \kk \rTTo^{\sS{_X}r_0}
(s\cc)^{-1}(X,Xi) \rMono^{\ju_1} (s\Dr(\cc|\cb))^{-1}(X,Xi) \bigr]. \]
We have proved in \secref{sec-invert-qis} that since $r_X$ is a
quasi-isomorphism, the above element is invertible modulo boundary in
the sense of \secref{sec-Invertible-trans}: there exist $p_0$, $v_0$,
$w_0$ such that equations~\eqref{eq-rp-i-pr-i} hold. We conclude by
\propref{pro-Invertible-rpfgBC} that $\overline{r}$ is invertible,
hence, $\Dr(\cc|\cb)$ and $\Dr(\ci|\cj)$ are equivalent.

Each acyclic K\n-injective complex $X$ is contractible. Indeed,
$\Kht(\ca)(X,X)\simeq\Dc(\ca)(X,X)$\linebreak[1]${}=0$ by
\cite[Proposition~1.5]{Spaltenstein88}. Hence, $\cj$ is a contractible
subcategory of $\ci$. Thus, $\ju:\ci\to\Dr(\ci|\cj)$ is an equivalence.
We deduce that $\Dr(\cc|\cb)$ and $\ci$ are equivalent in
$A_\infty^u$. Taking $H^0$ we get equivalent categories
$H^0(\Dr(\cc|\cb))$ and $H^0(\ci)$. The latter is a full subcategory of
$\Kht(\ca)$, whose objects are K\n-injective complexes. It is
equivalent to the derived category $\Dc(\ca)$ (e.g. by
\cite[Proposition~1.6.5]{KAScHapIWARA}). Hence, $H^0(\Dr(\cc|\cb))$ is
equivalent to the derived category $\Dc(\ca)$. This result follows also
from Drinfeld's theory \cite{Drinf:DGquot}. It motivated our study of
\ainf-categories.

Let $F:\ca\to\cb$ be an additive $\kk$\n-linear functor between Abelian
categories. The standard recipe \cite{Spaltenstein88} of producing its
right derived functor can be formulated in terms of the K\n-injective
resolution \ainf-functor $i$ as follows. Apply $H^0$ to the
\ainf-functor
\[ \Dr(\Com(\ca)|\Acyc(\ca)) \rTTo^\iu \Dr(\Inj(\ca)|\AcIn(\ca))
\rTTo^{\Dr(F)} \Dr(\Com(\cb)|\Acyc(\cb))
\]
(when $F(\Ob\AcIn(\ca))\subset\Ob\Acyc(\cb)$). Some work is required to
identify the obtained functor \cite[Section~8.13]{Lyu-AinfCat}
\(H^0(\iu\Dr(F)):\Dc(\ca)\simeq H^0\bigl(\Dr(\Com(\ca)|\Acyc(\ca))\bigr)
\to H^0\bigl(\Dr(\Com(\cb)|\Acyc(\cb))\bigr)\simeq\Dc(\cb)\)
with $RF$; however, we shall not consider this topic here.

\begin{acknowledgement}
We are grateful to B.~L. Tsygan for drawing our attention to the
subject. We thank V.~G. Drinfeld for the possibility to become
acquainted with his work on quotients of differential graded categories
before it was completed. We are grateful to all the participants of the
\ainf-category seminar at the Institute of Mathematics, Kyiv,
especially to Yu.~Bespalov and O.~Manzyuk. We thank the referee for
useful suggestions, which allowed us to improve the exposition.
\end{acknowledgement}

%\bibliography{yuri}
\end{document}